\documentclass[11pt]{amsart}
\usepackage[top=2cm, bottom=2cm, left=3cm, right=3cm]{geometry}
\usepackage{setspace}
\usepackage{amsmath,amssymb,amsthm}
\usepackage{mathrsfs}  
\usepackage{bbm}        
\usepackage{esint}     
\usepackage{dsfont}
\usepackage[dvipsnames]{xcolor}
\usepackage{cancel}
\usepackage{todonotes}

\usepackage[utf8]{inputenc}
\usepackage[T1]{fontenc}
\usepackage[english]{babel}
\usepackage[noadjust]{cite}

\usepackage{calc}
\usepackage[inline]{enumitem}
\usepackage[normalem]{ulem}
\usepackage{url}

\usepackage[bookmarks, bookmarksnumbered, pdfstartview={XYZ null null 1.00}]{hyperref}

\usepackage{amsthm,amssymb,epsfig,mathtools}
\usepackage{mathrsfs}  
\usepackage{bbm}        
\usepackage{esint}     
\usepackage{dsfont}
\usepackage{color}
\usepackage{url}
\usepackage{algorithm2e}

\usepackage{calc}
\usepackage[inline]{enumitem}
\usepackage[normalem]{ulem}
\usepackage{url}

\usepackage{tikz}
\usetikzlibrary{shapes, arrows, arrows.meta, decorations.markings,
  patterns, decorations.pathreplacing}
\usepackage{pgfplots}
\pgfplotsset{compat=1.18}

\usepackage{todonotes}
\usepackage[ulem=normalem,draft]{changes}
\usepackage{cancel}

\makeindex

\numberwithin{equation}{section}

\newcommand{\theoname}{Theorem}
\newcommand{\lemmname}{Lemma}
\newcommand{\coroname}{Corollary}
\newcommand{\propname}{Proposition}
\newcommand{\definame}{Definition}

\newcommand{\remkname}{Remark}
\newcommand{\explname}{Example}

\theoremstyle{plain}
\newtheorem{theorem}{\theoname}[section]
\newtheorem{lemma}[theorem]{\lemmname}
\newtheorem{corollary}[theorem]{\coroname}
\newtheorem*{corollary*}{\coroname}
\newtheorem{proposition}[theorem]{\propname}

\theoremstyle{definition}

\newtheorem{remark}[theorem]{\remkname}

\newlist{hypothesis}{enumerate}{1}
\setlist[hypothesis]{label={\textup{(H\arabic*)}}, ref={(H\arabic*)}, leftmargin=*, widest*=10}

  {\list{}{\leftmargin=.5in\rightmargin=.5in}  \item[]  }%
  {\endlist}

\def\dd{{\rm d}}
\newcommand{\eqdef}{\ensuremath{\stackrel{\mbox{\upshape\tiny def.}}{=}}}
\newcommand{\norm}[1]{\left\lVert#1\right\rVert}

\def\inf{\mathop{\rm inf}}
\def\sup{\mathop{\rm sup}}
\def\min{\mathop{\rm min}}
\def\max{\mathop{\rm max}}

\newcommand{\mres}{\mathbin{\vrule height 1.6ex depth 0pt width
		0.13ex\vrule height 0.13ex depth 0pt width 1.3ex}}

\def\divv{{\mathop{\rm div}}}

\def\det{\mathop{\rm det}}

\def\dist{{\rm dist}}

\def\sign{\mathop{\rm sign}}

\def\Ker{\mathop{\rm Ker}}

\DeclareMathOperator{\id}{id}

\providecommand{\co}{\operatorname{co}}

\newcommand{\R}{\mathbb R}
\newcommand{\sym}{\mathbin{\triangle}}
\newcommand{\dPL}{\delta_{\rm PL}}
\newcommand{\dBM}{\delta_{\rm BM}}
\newcommand{\dvol}{\delta_{\rm vol}}
\newcommand{\dratio}{\delta_{\rm ratio}}
\newcommand{\did}{\delta_{\rm id}}

\newcommand{\dLSI}{\delta_{\mathsf{LSI}}}
\newcommand{\dTal}{\delta_{\mathsf{Tal}}}

\title[A height transport approach to quantitative PL]{A height transportation approach to the quantitative Prékopa--Leindler inequality
}

\author{João Miguel Machado}
\address{
CREST, ENSAE Paris, Institut Polytechnique de Paris\\
5 avenue Henry Le Chatelier\\
91120 Palaiseau, France}
\email{joaomiguel.machado@ensae.fr}

\author{João P. G. Ramos}
\address{Instituto Nacional de Matem\'atica Pura e Aplicada (IMPA) \\ 
Estrada Dona Castorina 110 \\
22460-320, Jardim Bot\^anico - Rio de Janeiro, RJ, Brasil}
\email{joao.ramo@impa.br}

\begin{document}

\begin{abstract} We give a new, simpler proof of sharp stability for the Pr\'ekopa--Leindler inequality. In contrast to the main arguments from Figalli, van Hintum and Tiba~\cite{FigalliVanHintumTibaPL}, our proof hinges mostly on the sharp stability theory for the Brunn--Minkowski problem, developed by the same authors \cite{figalli2024BM}, and on the idea of using a height transportation approach in order to match level sets. 

As a by-product of our techniques, we are able to prove the sharp dependency on $\tau = \min\{\lambda,1-\lambda\}$ in the Pr\'ekopa--Leindler stability constant whenever one of the inputs is logarithmically concave, producing a new, simplified proof of two of the main results from~\cite{FigalliVanHintumTibaPL}. 

Finally, in the same fashion as Krist\'aly and Pîrvuceanu \cite{kristaly2026quantitative}, we are able to obtain new stability results for the Bakry-Emery log-Sobolev inequality and for the Talagrand inequality. Due to the sharpness of the inputs, these results are optimal in terms of exponents. We moreover show that the $L^1$ nature of the distance in such results \emph{cannot} in general be improved to any other $L^p$ distance, where $p>1.$

\bigskip


\medskip

\end{abstract}

\setcounter{tocdepth}{1}
\maketitle
\tableofcontents

\section{Introduction}\label{sec.introduction}

In recent years there has been an increasing interest in the geometry of log-concave functions. These objects and their theory can be interpreted as a functional version of the Brunn-Minkowski theory developed for convex bodies~\cite{Gardner2002,schneider2013convex,rotem2022surface}, and play a pivotal role in theoretical guarantees of sampling algorithms~\cite{lovasz2007geometry}. One of the major structuring tools in this theory is the Prékopa-Leindler inequality, which gives a concavity principle to the natural generalization to the notion of log-concave volume in this geometry, as observed for instance in~\cite{cordero2015moment}. 

The Prékopa-Leindler inequality says that given three measurable functions $f,g,h : \mathbb{R}^n \to \mathbb{R}_+$ satisfying the Prékopa condition: for some $\lambda \in (0,1)$ it holds that
\begin{equation}\label{eq.prekopa_condition}
  h(\lambda x + (1-\lambda)y) \ge f(x)^\lambda g(y)^{1-\lambda}, \qquad 
  \text{ for all } x,y \in \mathbb{R}^n,
\end{equation} 
then this inequality can be lifted to their integrals, namely 
\begin{equation}\label{eq.prekopa_leindler}
  \int_{\mathbb{R}^n} h(x) \dd x 
  \ge 
  {\left(
    \int_{\mathbb{R}^n} f(x) \dd x
  \right)}^{\lambda}
  {\left(
    \int_{\mathbb{R}^n} g(x) \dd x
  \right)}^{1-\lambda}. 
\end{equation}

To emphasize the geometry of log-concave functions, it is often formulated in terms of two functions $f,g$ and their Asplund sum, defined as 
\begin{equation}\label{eq.asplund_sum}
    h_\lambda(z) 
    \eqdef 
    \sup_{z = \lambda x + (1-\lambda)y} 
    f(x)^\lambda g(y)^{1-\lambda}.
\end{equation}
Then of course, any function $h$ such that $f,g,h$ satisfy~\eqref{eq.prekopa_condition} is bigger than $h_\lambda$ and this one is the sharpest function for which~\eqref{eq.prekopa_leindler} may hold.

Since the recent contribution of Figalli-van-Hintum-Tiba~\cite{figalli2025sharp}, this inequality has a sharp quantitative version, which can be interpreted as a strong concavity principle for the natural notion of volume in this geometry. This has  found many applications to the stability of other functional and transport inequalities, as well as other variational problems that are central in log-concave geometry~\cite{machado2025quantitative,bonnet2026quantitative,kristaly2026quantitative}.

\begin{theorem}\label{thm.quantitative_PL}
For every $n\geq1$, there exists a dimensional constant $C_n>0$ such that the following holds. Let $f,g,h:\mathbb{R}^n\to\mathbb{R}_+$ be integrable functions satisfying the Prékopa condition~\eqref{eq.prekopa_condition} for some $\lambda\in(0,1)$, with $\|f\|_{L^1},\|g\|_{L^1}>0$. Assume that at least one of $f,g$ is log-concave such that almost equality holds
\[
  \int_{\mathbb{R}^n}h(x)\dd x
  \le 
  (1 + \dPL) 
  {
  \left(
    \int_{\mathbb{R}^n}f(x)\dd x
    \right)
    }^{\lambda}
    {
    \left(
    \int_{\mathbb{R}^n}g(x)\dd x
    \right)
    }^{1-\lambda},
\]
where $\tau\eqdef\min\{\lambda,1-\lambda\}$. Then we have that 
\begin{equation}\label{eq.quantitative_PL}
 \inf_{z\in\mathbb{R}^n}
 \norm{
   \frac{f}{\|f\|_{L^1}}
   -\frac{g(\cdot+z)}{\|g\|_{L^1}}
 }_{L^1(\mathbb{R}^n)}
 \leq C_n\sqrt{\frac{\dPL}{\tau}}.
\end{equation}
The exponent $1/2$ and the order of dependence on $\lambda$ are optimal.
Moreover, one may choose $C_n\leq e^{Cn}$ for an absolute constant $C$.
\end{theorem}

\begin{remark}\label{rem.dimension_comparison}
Keeping track of the dimension is useful here. The factor $Ce^n$ obtained when joining the translations is only one part of the stability constant: one must also include the polynomial factors from the height estimates and the dimensional constant from Brunn--Minkowski. Together, these give the bound $e^{Cn}$ in Theorem~\ref{thm.quantitative_PL}, both for two log-concave inputs and for a single log-concave input.

For comparison, Theorem~1.9 and Remark~1.10 in the second arXiv version of~\cite{FigalliVanHintumTibaPL} (6 October 2026), specialized to $p=0$, give a constant of order $n^{O(n)}$ multiplying the same factor $\sqrt{\dPL/\tau}$. Thus the dimensional upper bound obtained here improves the one stated there; we do not, however, claim that this specific dependence on $n$ is optimal. For general inputs, the constants and the smallness threshold in the Brunn--Minkowski input remain part of our estimate, as recorded in Section~\ref{sec.1logconcave_general_functions}.
\end{remark}

With minor modifications to our proof, the same result holds for arbitrary integrable functions $f,g$ and $h$, but with a non-optimal dependence on $\lambda$. We refer to Section~\ref{sec.1logconcave_general_functions} for a discussion of this case. On the other hand, as discussed below in this introduction, we believe that our proof is particularly well suited to the case of log-concave functions, and we expect that the optimal dependence on $\lambda$ is rightly simplified by the following argument due to~\cite{boroczky2021stability}. 

\begin{corollary}\label{corollary.quant_mid_point_implies_full}
    Under the assumptions of Theorem~\ref{thm.quantitative_PL}, assume in addition that both $f$ and $g$ are log-concave. Then the midpoint case of estimate~\eqref{eq.quantitative_PL} implies the full estimate for every $\lambda\in(0,1)$, up to a universal multiplicative change in $C_n$.
\end{corollary}
\begin{proof}
    Let $f,g$ be log-concave functions, and assume without loss of generality that $f,g$ are normalized so that
    \[
         \int_{\mathbb{R}^n} f(x) \dd x =  \int_{\mathbb{R}^n} g(x) \dd x = 1.
    \]
    Define their Asplund interpolation as 
    \[
      h_t(z) 
      \eqdef 
      \sup_{z = t x + (1-t)y} f(x)^t g(y)^{1-t}, 
      \text{ and set }
      \Phi(t) \eqdef \int_{\mathbb{R}^n} h_t(z) \dd z.
    \]
    As a result, the Prékopa-Leindler inequality implies that 
    \[
      t \mapsto \log \Phi(t) \text{ is concave, and } 
      \Phi(0) = \Phi(1) = 1. 
    \] 
    Up to swapping $f,g$, we can assume that $\lambda \le 1/2$, so that log concavity gives 
    \[
      \log \Phi(\lambda) = \log \Phi\left(
        2\lambda \cdot \frac{1}{2} + (1 - 2\lambda)\cdot 0
      \right)
      \ge 2\lambda \log\Phi(1/2). 
    \]

   Thanks to the chosen normalization, the deficit of the Prékopa-Leindler inequality with convexity parameter $t$ is given by
    \[
      t \mapsto 
      \delta_{\rm PL}(t) 
      \eqdef 
      \Phi(t) - 1, 
    \]
    hence assuming that the stability result holds for $t = 1/2$ with a universal stability constant $C_n$, we obtain the log-concavity of $\Phi$ that 
    \begin{align*}
      \frac{1}{C_n} 
      \inf_{z \in \mathbb{R}^n} 
      \norm{f - g(\cdot + z)}_{L^1(\mathbb{R}^n)}^2 
      \le \delta_{\rm PL}(1/2) 
      \le {\left(1 + \delta_{\rm PL}(\lambda)\right)}^{\frac{1}{2\lambda}} - 1.
    \end{align*}
    Since $\norm{f - g(\cdot + z)}_{L^1(\mathbb{R}^n)} \le 2$, up to increasing the constant $C_n$ it suffices to treat $\delta_{\rm PL}(\lambda)\leq\lambda$. In this case,
    \[
      {\left(1+\delta_{\rm PL}(\lambda)\right)}^{\frac{1}{2\lambda}}-1
      \leq \exp\left(\frac{\delta_{\rm PL}(\lambda)}{2\lambda}\right)-1
      \leq \frac{\delta_{\rm PL}(\lambda)}{\lambda},
    \]
    and the above estimates give
    \begin{align*}
      \inf_{z \in \mathbb{R}^n} 
      \norm{f - g(\cdot + z)}_{L^1(\mathbb{R}^n)} 
      \le
      C_n
      \sqrt{
        \frac{\delta_{\rm PL}(\lambda)}{\lambda}
      }.
    \end{align*}
    Since $\tau=\lambda$ and the original $h$ dominates $h_\lambda$, the result follows by mass normalization.
\end{proof}

\subsection*{Scheme of the proofs}

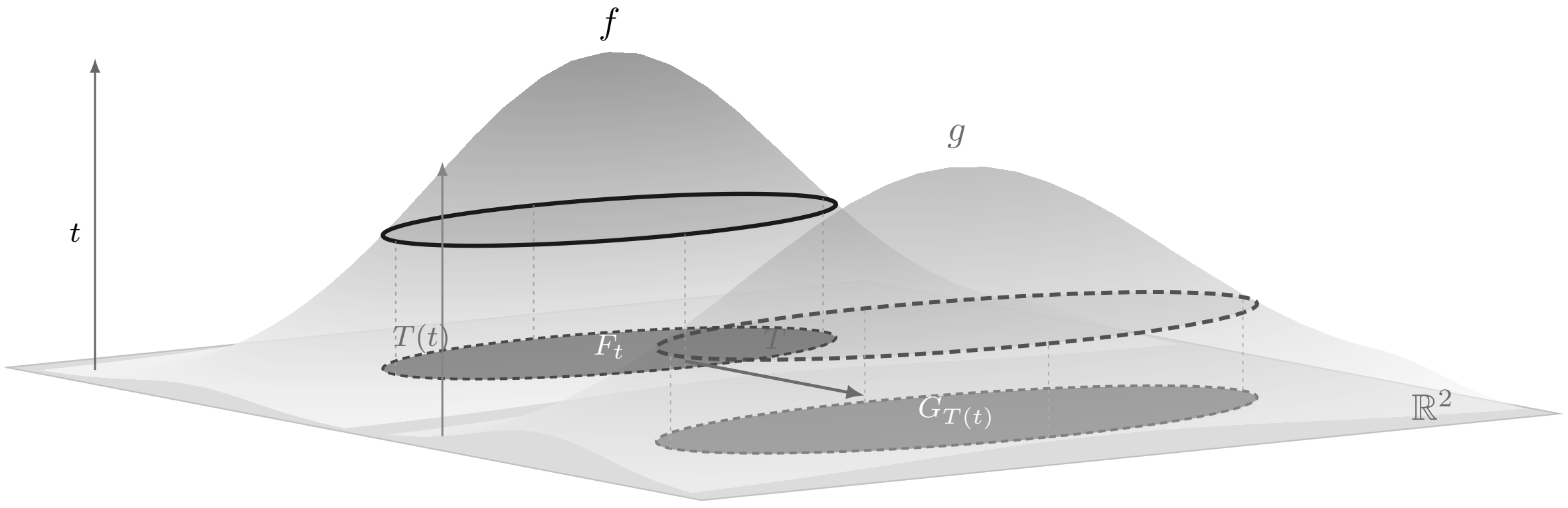
\begin{figure}[h!]
  \begin{center}
  \resizebox{\linewidth}{!}{
\begin{tikzpicture}
\begin{axis}[
  name=layeraxis,
  width=17.8cm,
  height=7.0cm,
  view={54}{27},
  axis lines=none,
  xmin=-4.1,xmax=4.1,
  ymin=-2.0,ymax=2.0,
  zmin=0,zmax=1.08,
  domain=-3.7:3.7,
  y domain=-1.7:1.7,
  samples=37,
  samples y=29,
  z buffer=sort,
  clip=false,
  colormap={surfacegray}{color(0cm)=(black!3); color(1cm)=(black!48)},
]

\addplot3[
  surf,
  shader=flat,
  draw=black!34,
  fill=black!19,
  opacity=.72,
  domain=-4:4,
  y domain=-1.75:1.75,
  samples=2,
  samples y=2
] {0};

\addplot3[
  surf,
  shader=interp,
  faceted color=black!16,
  opacity=.82,
  domain=-3.75:-.25,
  y domain=-1.7:1.7
] {exp(-1.05*((x+2)^2+y^2))};

\addplot3[
  surf,
  shader=interp,
  faceted color=black!22,
  opacity=.60,
  domain=.25:3.75,
  y domain=-1.7:1.7
] {.84*exp(-.86*((x-2)^2+.82*y^2))};

\path[fill=black!58,opacity=.15,draw=none]
  plot[domain=0:360,samples=60,variable=\a]
  (axis cs:{-2+.872*cos(\a)},{.872*sin(\a)},.45) -- cycle;
\draw[black!90,line width=1.15pt]
  plot[domain=0:360,samples=100,variable=\a]
  (axis cs:{-2+.872*cos(\a)},{.872*sin(\a)},.45) -- cycle;

\path[fill=black!48,opacity=.14,draw=none]
  plot[domain=0:360,samples=60,variable=\a]
  (axis cs:{2+1.058*cos(\a)},{1.168*sin(\a)},.32) -- cycle;
\draw[black!68,line width=1.1pt,dash pattern=on 3pt off 1.7pt]
  plot[domain=0:360,samples=100,variable=\a]
  (axis cs:{2+1.058*cos(\a)},{1.168*sin(\a)},.32) -- cycle;

\path[fill=black!70,opacity=.58,draw=none]
  plot[domain=0:360,samples=60,variable=\a]
  (axis cs:{-2+.872*cos(\a)},{.872*sin(\a)},.012) -- cycle;
\draw[black!72,line width=.75pt,dash pattern=on 2pt off 1.4pt]
  plot[domain=0:360,samples=100,variable=\a]
  (axis cs:{-2+.872*cos(\a)},{.872*sin(\a)},.014) -- cycle;

\path[fill=black!60,opacity=.54,draw=none]
  plot[domain=0:360,samples=60,variable=\a]
  (axis cs:{2+1.058*cos(\a)},{1.168*sin(\a)},.012) -- cycle;
\draw[black!50,line width=.75pt,dash pattern=on 2pt off 1.4pt]
  plot[domain=0:360,samples=100,variable=\a]
  (axis cs:{2+1.058*cos(\a)},{1.168*sin(\a)},.014) -- cycle;

\draw[black!38,line width=.35pt,dash pattern=on 1pt off 1.4pt]
  (axis cs:-2.872,0,.014) -- (axis cs:-2.872,0,.45);
\draw[black!38,line width=.35pt,dash pattern=on 1pt off 1.4pt]
  (axis cs:-1.128,0,.014) -- (axis cs:-1.128,0,.45);
\draw[black!38,line width=.35pt,dash pattern=on 1pt off 1.4pt]
  (axis cs:-2,.872,.014) -- (axis cs:-2,.872,.45);
\draw[black!38,line width=.35pt,dash pattern=on 1pt off 1.4pt]
  (axis cs:-2,-.872,.014) -- (axis cs:-2,-.872,.45);
\draw[black!32,line width=.35pt,dash pattern=on 1pt off 1.4pt]
  (axis cs:.942,0,.014) -- (axis cs:.942,0,.32);
\draw[black!32,line width=.35pt,dash pattern=on 1pt off 1.4pt]
  (axis cs:3.058,0,.014) -- (axis cs:3.058,0,.32);
\draw[black!32,line width=.35pt,dash pattern=on 1pt off 1.4pt]
  (axis cs:2,1.168,.014) -- (axis cs:2,1.168,.32);
\draw[black!32,line width=.35pt,dash pattern=on 1pt off 1.4pt]
  (axis cs:2,-1.168,.014) -- (axis cs:2,-1.168,.32);

\draw[-{Latex[length=4pt,width=3pt]},black!60,line width=.55pt]
  (axis cs:-3.55,-1.55,0) -- (axis cs:-3.55,-1.55,1.02);
\draw[-{Latex[length=4pt,width=3pt]},black!48,line width=.55pt]
  (axis cs:.45,-1.55,0) -- (axis cs:.45,-1.55,.90);

\node[font=\small] at (axis cs:-2,0,1.09) {$f$};
\node[font=\small,text=black!58] at (axis cs:2,0,.94) {$g$};
\node[font=\scriptsize,text=white]
  at (axis cs:-2,0,.035) {$F_t$};
\node[font=\scriptsize,text=white]
  at (axis cs:2,0,.035) {$G_{T(t)}$};
\node[font=\scriptsize] at (axis cs:-3.61,-1.61,.45) {$t$};
\node[font=\scriptsize,text=black!58] at (axis cs:.38,-1.61,.32) {$T(t)$};
\node[font=\small,text=black!58]
  at (axis cs:3.45,1.43,.025) {$\mathbb R^2$};

\draw[-{Latex[length=5pt,width=4pt]},black!58,line width=.85pt]
  (axis cs:-1.128,0,.035)
  -- node[midway,above=4pt,font=\scriptsize] {$T$}
  (axis cs:.942,0,.035);

\end{axis}
\end{tikzpicture}}
  \end{center}
  \caption{\emph{The qualitative mechanism: comparing suitable level sets via the Brunn--Minkowski inequality.}}
  \label{fig.height_transport}
\end{figure}

Our first contribution is a new \textit{height transport} proof of the classical Prékopa-Leindler inequality, which also yields a full characterization of the equality cases discovered by Dubuc~\cite{dubuc1977criteres}. The height transport method is based on the intuition that the Prékopa--Leindler inequality is the functional version of Brunn--Minkowski, and therefore it would be natural to expect a proof of the former by matching pairs of heights of the functions $f,g$ and applying the latter inequality to the corresponding matched level-sets, as illustrated in Figure~\ref{fig.height_transport}.

\begin{figure}[h!]
  \begin{center}
    \resizebox{.96\linewidth}{!}{
\begin{tikzpicture}[
  x=9.2cm,y=1.12cm,
  source/.style={draw=black!86,line width=1.05pt},
  target/.style={draw=black!55,line width=1.05pt,dash pattern=on 3pt off 1.8pt},
  guide/.style={draw=black!30,line width=.42pt,dash pattern=on 1.2pt off 1.7pt},
  maparrow/.style={{Latex[length=4.2pt,width=3.2pt]}-{Latex[length=4.2pt,width=3.2pt]},draw=black!58,line width=.75pt}
]
\path[use as bounding box] (-.06,-.53) rectangle (1.66,3.43);

\draw[black!72,line width=.5pt] (0,3.20) -- (0,0) -- (1.58,0);
\node[anchor=north,font=\scriptsize] at (0,-.08) {$0$};
\node[anchor=west,font=\scriptsize] at (1.58,0) {$t$};

\path[fill=black!12,draw=none]
  (0,0)
  -- plot[domain=0:.44,samples=70,smooth]
       (\x,{3*(1-\x)^2})
  -- (.44,0) -- cycle;
\path[pattern=north east lines,pattern color=black!28,draw=none]
  (0,0)
  -- plot[domain=0:.66,samples=90,smooth]
       (\x,{2*(1-\x/1.5)^2})
  -- (.66,0) -- cycle;

\draw[source]
  plot[domain=0:1,samples=150,smooth]
    (\x,{3*(1-\x)^2});
\draw[target]
  plot[domain=0:1.5,samples=190,smooth]
    (\x,{2*(1-\x/1.5)^2});

\draw[guide] (.44,0) -- (.44,{3*(1-.44)^2});
\draw[guide] (.66,0) -- (.66,{2*(1-.66/1.5)^2});
\draw[maparrow] (.44,-.22) -- (.66,-.22);
\node[font=\scriptsize] at (.55,-.40) {$T$};
\node[font=\scriptsize,anchor=north] at (.44,-.02) {$t$};
\node[font=\scriptsize,anchor=north] at (.66,-.02) {$T(t)$};

\node[font=\scriptsize,anchor=north] at (1,-.02) {$\|f\|_\infty$};
\node[font=\scriptsize,anchor=north,text=black!55]
  at (1.5,-.02) {$\|g\|_\infty$};

\draw[source] (1.02,2.95) -- (1.10,2.95);
\node[anchor=west,font=\scriptsize] at (1.12,2.95) {$|F_t|$};
\draw[target] (1.30,2.95) -- (1.38,2.95);
\node[anchor=west,font=\scriptsize,text=black!55] at (1.40,2.95) {$|G_t|$};
\node[font=\scriptsize,text=black!60] at (1.22,2.56)
  {equal shaded mass};
\end{tikzpicture}}
  \end{center}
  \caption{\emph{Heights are matched via the Brenier map transporting the height distribution of each function.}}
  \label{fig.height_transport_monotone_map}
\end{figure}
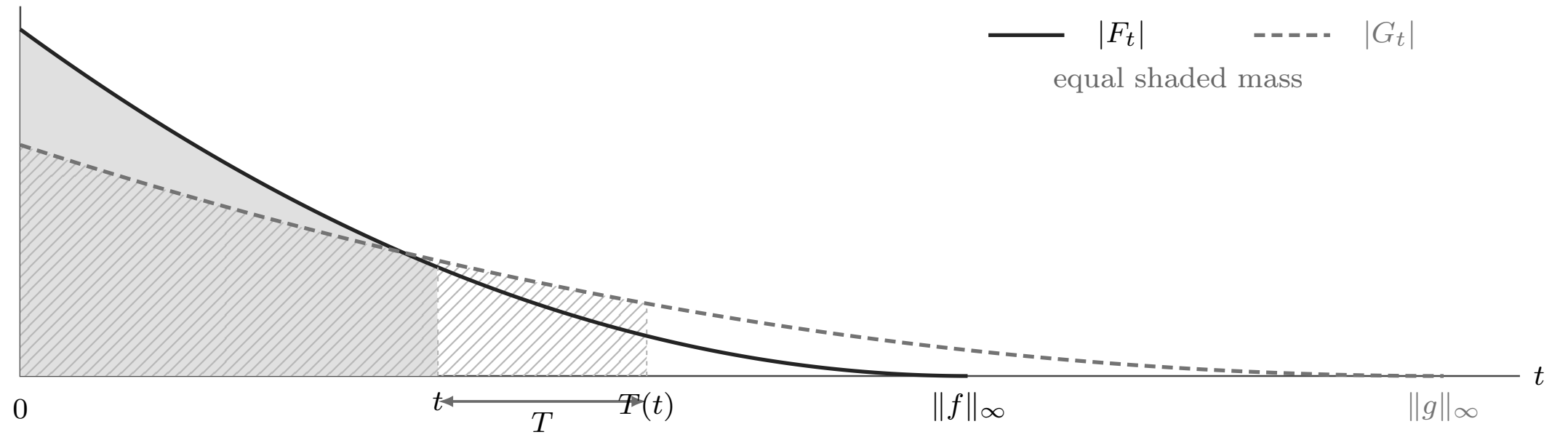

As a result, besides the Brunn-Minkowski theory enabling this argument, the backbone of the proof relies on how to match these level sets. The naive, but natural, strategy of matching sets of the same height quickly finds complications if the two functions $f,g$ have different maxima, since in this case many heights of one of them would become unmatched. Instead, this matching is performed by means of the \textit{Brenier map} transporting the heights distributions, defined as 
\[
    |F_t|\dd t, \quad |G_t|\dd t, \qquad \text{ where } \qquad 
    F_t = \{f > t\}, \ G_t = \{g > t\},
\] 
onto each other, see Figure~\ref{fig.height_transport_monotone_map}. This map is defined as the minimal map, whenever it exists, in the $L^2$-optimal transport problem between two measures $\mu,\nu \in \mathscr{P}(\mathbb{R}^d)$ by the map that minimizes the following problem whose value is the squared $2$-Wasserstein distance 
\[
  W_2^2(\mu, \nu) 
  \eqdef 
  \inf_{T_\sharp \mu = \nu} 
  \int_{\mathbb{R}^d} 
  |x - T(x)|^2 \dd \mu(x),
\]
where $T_\sharp \mu$ denotes the pushforward measure given by $T_\sharp \mu(B) = \mu(T^{-1}(B))$, for any Borel set $B \subset \mathbb{R}^d$. The quantity $W_2$ is known to define a distance between probability measures~\cite{ambrosio2021lectures,villani2009optimal}, and Brenier's theorem~\cite{brenier1991polar} asserts that whenever $\mu \ll \mathscr{L}^n$, this infimum is attained by an unique map $T = \nabla \varphi$, which can be proven to be the gradient of a convex function. 

In our case, where we consider the transport of height distributions, the map $T:\mathbb{R}\to \mathbb{R}$ coincides $\mathscr{L}^1$-a.e. with a monotone function, which will be used extensively throughout our arguments. In addition, this choice seems particularly fitting as it mirrors the use of the Brenier map in the proof of the Brunn-Minkowski inequality via Gromov's argument, see~\cite{figalli2009refined,figalli2010mass}, where the particular structure of the Brenier map has historically played a crucial role in proving sharp stability results. 

\begin{figure}[h!]
  \begin{center}
    \resizebox{\linewidth}{!}{\input{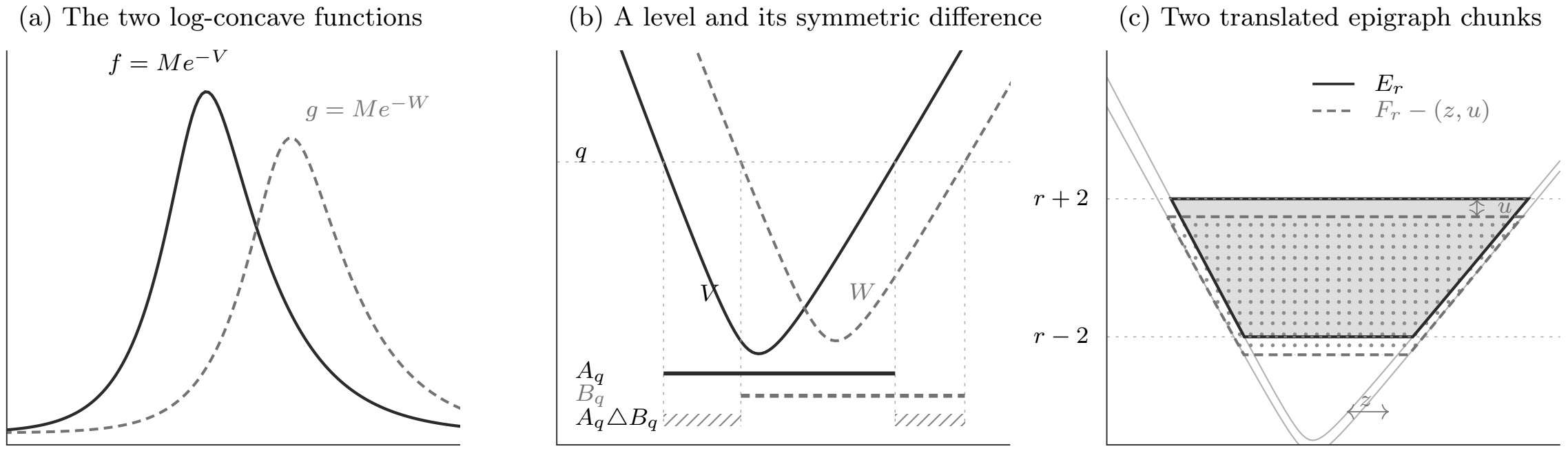}}
  \end{center}
  \caption{\emph{The quantitative mechanism: from levels to chunks.}}
  \label{fig.chunks_level_sets}
\end{figure}

This strategy of first matching heights via an optimal transport step and then comparing the resulting family of convex bodies via Brunn--Minkowski naturally suggests a quantitative extension. Indeed, as we shall see, this argument leads to an expression where \textit{the deficit of the Prékopa--Leindler inequality bounds an integral quantity of the Brunn-Minkowski deficit between matched level sets.} The issue is that each Brunn-Minkowski deficit naturally gives rise to a translation of a pair of matched level sets, yielding a continuum of such translations. To circumvent this problem, the idea is now to consider \textit{chunks of level-sets}, see Figure~\ref{fig.chunks_level_sets}. This forces us to pass to a suitable logarithmic scale, since chunks of level sets of log-concave functions are not necessarily convex bodies anymore. In fact, this change of scale is quite convenient as it allows us to work with the associated convex potentials and obtain many useful estimates, cf.~Lemma~\ref{lemma.tail_ab}. More importantly, we can then apply the quantitative version of Brunn-Minkowski to pairs of chunks, which require only countably many translations instead of a continuum. 

On the other hand, another problem arises with this simplification: if before the Brunn-Minkowski inequality would give rise to a translation in the ambient space, applying it for level set chunks gives rise to a translation in the natural space of the epigraph, meaning a translation $(z,u) \in \mathbb{R}^{n+1}$ composed of a horizontal translation and a height translation. The information on both sources of translations is contained in the height distributions, which in the logarithmic scale assume the form 
\[
  \dd \mu(q) = M e^{-q} a(q) \dd q, \quad 
  \dd \nu(q) = M e^{-q} b(q) \dd q,
\]
where $a(q) = |\{V < q\}|$ and $b(q) = |\{W < q\}|$, for convex functions $V,W$ such that $f = Me^{-V}$ and $g = M e^{-W}$.

In fact, all the information on both types of translations is encapsulated in a suitable distance between $\mu, \nu$, and we prove the following bounds: there is a dimensional constant $C_n$ such that 
\[
  \norm{\mu - \nu}_{L^1(\mathbb{R})}
  \le 
  C_n \sqrt{\dPL}, \qquad 
  W_2^2(\mu,\nu)
  \le 
  C_n \sqrt{\dPL}. 
\]
While the $L^1$ norm is a \textit{vertical distance}, shown to be the natural quantity to control the averaged vertical translations furnished by Brunn--Minkowski, the Wasserstein distance $W_2$ is a \textit{horizontal distance} and is shown to control the translations in the ambient space. 

\begin{figure}[h!]
  \begin{center}
    \resizebox{\linewidth}{!}{\input{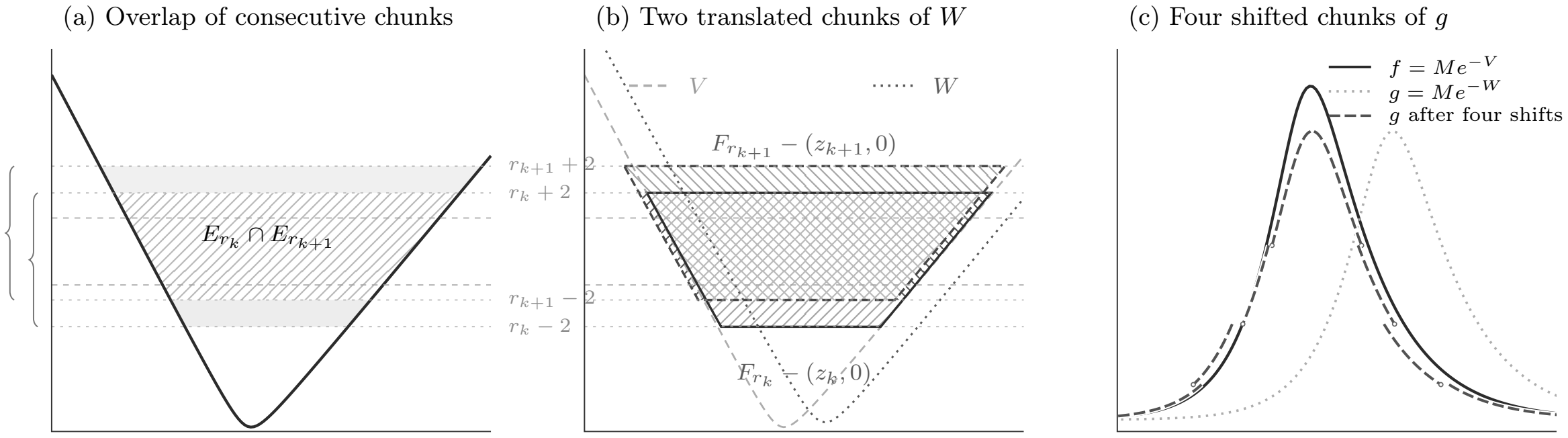}}
  \end{center}
  \caption{\emph{The quantitative mechanism: chunk overlap allows for a common translation.}}
  \label{fig.chunks_reconstruction}
\end{figure}

Our final argument, see the scheme in Figure~\ref{fig.chunks_reconstruction}, then goes as follows: first the $L^1$ bounds allow us to remove the height translations coming from the application of Brunn-Minkowski in the epigraph space, while only paying a factor of the order $\sqrt{\dPL}$, our final estimate. Finally, for the horizontal translations, we naturally cannot remove them all together, as it is a natural symmetry of the inequality. Taking however sufficiently large chunks, in such a way that they overlap, enables us to show that, since each chunk carries information of both its predecessor and successor, all horizontal translations cannot vary too much and might be chosen to be the same. 

\subsection*{One log-concave input and beyond}

The qualitative proof is completely independent of the log-concavity of the functions $f,g$. For the quantitative case, this property gives quantitative growth rates to the densities $\mu, \nu$ in logarithmic scale, which is useful in the two distance estimates described above. 

Without these assumptions, these estimates can still be proven if we adapt the way that we select level sets to form chunks. We show that the way of constructing level set chunks described above is equivalent to performing convolutions with a squared window kernel. We show then that, without the growth properties endowed by the log-concave geometry,  analogous Wasserstein estimates can be obtained if we replace this kernel with a smoother one, satisfying $\sqrt{\kappa} \in H^1(\mathbb{R})$. 

The only issue with this is that chunks must then be measured with a weightened measure, instead of the uniform Lebesgue measure over $\mathbb{R}^{n+1}$. The solution is to disintegrate this weighted measure into the uniform measure over $\mathbb{R}^{n+3}$ and unfold the level set chunks in the same space. 

This is where the real bottleneck of this proof appears: once these operations are performed, the quantitative result is as strong as the Brunn-Minkowski technology employed. To recover a stability version that presents at the same time the sharp exponent and the optimal dependence on $\lambda$, we need a stability version of the Brunn-Minkowski inequality that has the corresponding behavior. That is why in the case of one log-concave input, we can use the sharp result of Barchiesi--Julin~\cite{barchiesi2017robustness}, which gives precisely the optimal behavior under the assumption that only one of the sets is a convex body. This reasoning can be pushed even further, and we may recover a sharp result for general functions using the recent sharp quantitative Brunn--Minkowski theory ~\cite{figalli2024BM}.

\subsection*{Towards an \texorpdfstring{$L^1$}{L1} stability theory for functional and transport inequalities}

A large part of the quantitative theory of sharp functional inequalities is Hilbertian. Indeed, the question of Brezis and Lieb~\cite{brezis1985sobolev} on the stability of the Sobolev inequality was answered by Bianchi and Egnell~\cite{bianchi1991note} through a quadratic control of the distance to the Sobolev optimizers in the natural Sobolev norm.  As another example, we mention the use of spectral gap methods, as in~\cite{bonnefont2016spectral,andrade2026trudinger}. In these settings, the distance to the manifold of minimizers is controlled via a spectral gap condition for a suitable operator that arises in the Euler-Lagrange equations satisfied by optimizers. This is a recurrent, though by no means universal, source of stability phenomena in which a deficit controls an $L^2$-based energy distance.

The $L^1$ stability result for the Brascamp--Lieb inequality in~\cite{machado2025quantitative} suggests a complementary picture. Linearizing the Prékopa--Leindler inequality has shown to provide $L^1$ stability results with a constant that remains uniform over a large class of convex reference potentials. With that in mind, we shall employ some of our main results to obtain sharp stability estimates for the Bakry--Émery log-Sobolev and Talagrand inequalities. We will find that their density formulations admit sharp, potential-independent $L^1$ stability, whereas the analogous density estimates fail in every $L^p$ topology with $p>1$, even for the fixed Gaussian potential. Thus $L^1$ is not merely a convenient weak topology in this argument: it is the endpoint of a sharp uniform theory.

n~\cite{kristaly2026quantitative}, Krist\'aly and P{\^\i}rvuceanu first quantified the passage from Prékopa--Leindler to log-Sobolev and Talagrand stability. Their argument follows Maurey's method of linearizing the Prékopa--Leindler inequality by infimal convolution~\cite{maurey1991deviation} later expanded by Bobkov and Ledoux~\cite{bobkov2000brunn}. Their general exponent is $1/19$, reflecting the Prékopa--Leindler stability available at the time that presented the optimal dependence of the stability constant on the interpolation parameter~\cite{boroczky2021stability}. In the radial setting, they obtained the sharp exponent $1/2$. 

The quantitative stability result of Prékopa-Leindler, developed in the present manuscript, having simultaneously the sharp exponent and the optimal dependence of the constant on $\lambda$ for log-concave inputs, allows us to follow the arguments from~\cite{kristaly2026quantitative} to obtain sharp $L^1$ quantitative stability results for both log-Sobolev and Talagrand inequalities. In addition, since we may assume that a single input is log-concave, this yields results formulated for general probability
distributions with finite Fisher information, rather than only for the log-concave class considered in~\cite{kristaly2026quantitative}. 

While we do not claim to have a novel method for the sharp $L^1$-stability proofs for log-Sobolev and Talagrand inequalities, our interest is to shed new light into the universality of the stability constant w.r.t.~the reference measure's convex potential, and more importantly, how such behavior cannot hold in any stronger $L^p$ topologies with $p>1$. In particular, this shows that even for the Prékopa--Leindler inequality \textit{there is no estimate such as~\ref{thm.quantitative_PL} with a stability constant $C_n$ independent of $f,g$ controlling an $L^p$ distance with $p>1$.} This goes against the $L^2$ stability of Prékopa-Leindler suggested in~\cite{boroczky20262}. Indeed, if such stronger stability version of Prékopa--Leindler could be proved, this would imply a stronger $L^p$ stability result with $p>1$ for all the inequalities studied here: log-Sobolev and Talagrand, as well as the Brascamp-Lieb inequality. This would contradict Propositions~\ref{prop.app_Lp_LS_failure} and~\ref{prop.app_Lp_Tal_failure}. 

\subsection*{Alternative presentation and Lean formalization} In an effort to present the results of this paper in a more reader-friendly way, we have attempted to make an interactive version of the exposition, where one may see different parts of the proof with different levels of detail. In order to do so, we stored all results in this manuscript in the website 
\begin{center}
\href{https://w3.impa.br/~joao.ramos/height-transport-pl/}{\texttt{w3.impa.br/~joao.ramos/height-transport-pl/.}}
\end{center}
There, one may find a more detailed explanation of the importance of the Pr\'ekopa--Leindler inequality and its stability version, as well as the aforementioned interactive manuscript version, with a dependence tree of ideas used, expository videos targeted at explaining its different parts, and a dissected version of a Lean formalization of the results.

Regarding this last feature, we mention two caveats. The first is that we opted not to make any attempt to formalize the applications to the log-Sobolev and Talagrand inequalities, as these are only tangentially related to the main topic of the manuscript. The other one is that, although we formalized Theorem \ref{thm.quantitative_PL} in \emph{all} levels of generality in Lean, its formalization contains an explicit dependence on three results, depending on the respective generality level: the original sharp stability for Brunn--Minkowski for convex sets \cite{figalli2009refined}, the same result with the relaxed condition that only one of the convex sets is convex \cite{barchiesi2017robustness}, and the recent full stability with sharp exponent for the Brunn--Minkowski inequality \cite{figalli2024BM}. The choice of doing so was due to the length and effective cost that formalizing such theory-heavy results would bring about, and we would encourage the community (in particular, the authors of the respective manuscripts) to make a push towards formalizing them.  A more detailed comment on the Lean formalization is available in the website above,

\subsection*{Organization of this manuscript}
We start our analysis with the new proof of the Prékopa--Leindler inequality in Section~\ref{sec.qualitative_proof}. In Section~\ref{sec.quantitativePL} we give the proof in the case that both $f$ and $g$ are log-concave functions. In Section~\ref{sec.1logconcave_proof}, we show how the previous proof can be refined to the case of one log-concave function, still having the sharp exponent and the optimal dependence on the convexity parameter. In Section~\ref{sec.1logconcave_general_functions} we explain how these arguments also give the sharp exponent for general inputs. Finally, in Section~\ref{sec.applications}, we discuss applications to the log-Sobolev and Talagrand inequalities following the work of~\cite{kristaly2026quantitative}, showing that only the $L^1$ topology can yield such stability results. 

\section*{Acknowledgments}
J.M.M. thanks the hospitality of IMPA, during the summer program of January 2026. He was partially funded by the European Union (ERC, Optinfinite, 101201229). The views and opin-
ions expressed are, however, of the author(s) only and do not necessarily reflect those of the European Union or the European Research Council Executive Agency. Neither the European Union nor the granting authority can be held responsible for them.

J. P. G. R.  was supported by the FCT through project SHADE
(2023.17881.ICDT), DOI \allowbreak \texttt{10.54499/}
\texttt{2023.17881.ICDT}. He was also supported by FAPERJ through  JCNE grant
no.~SEI-260003/ \allowbreak 020475/2025, by Instituto Serrapilheira through
grant Serra-R-2510-59765, and by CNPq through Research Productivity Fellowship no. 303861/2026-6.

\section*{Statement of LLM usage}
The qualitative argument was first conceived in January of 2026 when the first author was a visiting researcher in the IMPA summer program. 

Although the quantitative arguments were conceptually clear to the authors, different rounds of interactions with Large Language Models helped catalyze the quantitative argument. Various versions of LLMs were used to explore different strategies for the quantitative proof, until it became clear to the authors that the conceptually clean idea that horizontal and vertical translations of epigraph chunks should be controlled with a Wasserstein and a $L^1$ distance, respectively. This led to the current path of proof with two log-concave inputs, which was then executed by the authors with AI assistence. 

The authors then realized that very few steps of the proof critically relied on the log-concave inputs. At this stage the authors used AI tools to identify all the steps in the proof that used this property. A frontier model then suggested the lift to $\mathbb{R}^{n+3}$ to avoid the estimates along Wasserstein geodesics with log-concave end-points. It then became clear to the authors that the bottleneck of the present argument is only the Brunn-Minkowski technology employed. In particular, that even the case of general inputs can be dealt with the present argument, but having a non-optimal dependence on the parameter $\lambda$. 

Although several technicalities were first checked with the help of LLMs, the manuscript's final version was in its majority human-written, completely thoroughly revised, and several arguments have been changed, simplified and improved by the authors. As the previous two paragraphs show, the main ideas in this manuscript are also completely due to both authors, with AI only helping expand, polish and verify them. The authors take full responsibility for the correctness of the mathematical contents of the paper. 



\section{The height transport method for the Prékopa-Leindler inequality}\label{sec.qualitative_proof}

\subsection{Gromov's argument, Brenier's map and multimarginal optimal transport}

The two backbones of our argument are the Brunn-Minkowski inequality and optimal transport theory. The former says that the euclidean volume, when see as a functional over the space of convex bodies or even measurable sets, is concave with respect to Minkowski addition: for measurable $A,B \subset \mathbb{R}^n$ it holds that
\begin{equation}\label{eq.BrunnMinkowski}
  |A + B|^{1/n} 
  \ge 
  |A|^{1/n} + |B|^{1/n}.
\end{equation}
The by now famous proof of Gromov\footnote{
    The original argument was formulated in~\cite[Appendix~I]{milman1986asymptotic} for the isoperimetric inequality.
} of the Brunn-Minkowski inequality is based on the idea of transporting the mass uniformly distributed over one set to the other. 

A summarized modern version of this argument is as follows: given two bounded measurable sets $A,B \subset \mathbb{R}^n$ of positive volume, consider the two probability measures uniformly distributed among them 
\[
    \dd\mu(x) \eqdef \frac{\mathbf{1}_A(x)}{|A|}\dd x, \quad 
    \dd\nu(x) \eqdef \frac{\mathbf{1}_B(x)}{|B|}\dd x. 
\]
Brenier's polar factorization theorem~\cite{brenier1991polar} states that there is a map $T = \nabla \varphi$, given by the gradient of a convex function such that 
\begin{equation}\label{eq.pushforward}
  T_\sharp \mu = \nu,
\end{equation}
where the push-forward measure is defined as $T_\sharp \mu(C) \eqdef \mu(T^{-1}(C))$, for any measurable $C \subset \mathbb{R}^n$. 

At least formally, equation~\eqref{eq.pushforward} can be seen as a change of variables formula, see~\cite[Lecture~6]{ambrosio2021lectures}. Exploiting the convexity of the potential $\varphi$, this is expressed as the following Monge--Ampère equation
\[
  \frac{|B|}{|A|} = \det DT(x),
  \qquad DT(x) = D^2\varphi(x) \ge 0,
  \quad \mu\text{-a.e.}
\]
which can be made rigorous since, as the gradient of a convex function, the Brenier map $T$ has bounded variation~{\color{black}\cite{alberti1999geometrical}} and is a monotone map. As a result, $\id+T$ is injective and, after discarding a null subset of $A$, satisfies $(\id+T)(A)\subseteq A+B$. Thus the area formula and the concavity of $\det^{1/n}$ yield
\begin{align*}
  |A+B|^{1/n}
  &\ge \left(\int_A \det(I+DT(x))\dd x\right)^{1/n}\\
  &\ge \left(\int_A \left(1+\det DT(x)^{1/n}\right)^n\dd x\right)^{1/n}\\
  &= |A|^{1/n}\left(1+\left(\frac{|B|}{|A|}\right)^{1/n}\right)
   = |A|^{1/n}+|B|^{1/n}.
\end{align*}
The reader is referred to~\cite{figalli2009refined} for more details. 

While Gromov's original argument can be carried out with any map performing the transportation~\eqref{eq.pushforward}, the fine structure of the Brenier map has been exploited in the seminal works~\cite{figalli2009refined,figalli2010mass} to prove sharp quantitative stability versions of the isoperimetric and Brunn--Minkowski inequalities. In the present work, we exploit the deeper variational and metric properties of the Brenier map, through the theory of optimal transport. 

More precisely, let $\mathscr P_2(\R^n)$ denote the probability measures with finite second moment. For $\mu,\nu\in\mathscr P_2(\R^n)$, let $\Pi(\mu,\nu)$ denote the probability measures on $\R^n\times\R^n$ with marginals $\mu$ and $\nu$. The quadratic Kantorovich problem defines the squared $2$-Wasserstein distance by
\begin{equation}\label{eq.wasserstein}
 \begin{aligned}
 W_2^2(\mu,\nu)&\eqdef\inf_{\pi\in\Pi(\mu,\nu)}
                       \int|x-y|^2\dd\pi(x,y)
 =\min_{T_\sharp\mu=\nu}\int_{\R^n}|x-T(x)|^2\dd\mu(x),
 \end{aligned}
\end{equation}
where the last equality holds if $\mu\ll\mathscr L^n$. In this case, Brenier's theorem identifies the unique optimal plan as $(\id,T)_\sharp\mu$, where $T=\nabla\varphi$ is the gradient of a convex function and the unique $\mu$-a.e. minimizing map. 

\subsection{The height transport method} 

We start with the height transport proof of the Prékopa-Leindler inequality. In spirit, it follows Gromov's proof of the isoperimetric and Brunn--Minkowski inequalities~\cite[Appendix~I]{milman1986asymptotic}, but at the level of heights. First we match the level sets of the inputs through their height distributions, then the matched sets are compared with the Brunn--Minkowski inequality.

\begin{theorem}\label{thm.qualitative_PL}
    Let $f,g,h : \mathbb{R}^n \to \mathbb{R}_+$ be measurable functions satisfying the Prékopa condition~\eqref{eq.prekopa_condition} for some $\lambda \in (0,1)$, then the Prékopa-Leindler inequality~\eqref{eq.prekopa_leindler} holds. In addition, equality holds if and only if there exists a log-concave function $\varrho : \mathbb{R}^n \mapsto \mathbb{R}_+$ and a vector $z \in \mathbb{R}^n$ such that 
    \begin{equation}\label{eq.equality_cases_thm}
      h = \varrho, \quad 
      f = a^{1-\lambda} \varrho(\cdot + (1-\lambda)z), \quad 
      g = a^{-\lambda} \varrho(\cdot + -\lambda z), \quad 
      a \eqdef \frac{\norm{f}_{L^1}}{\norm{g}_{L^1}}.
    \end{equation} 
\end{theorem}
\begin{proof}
  Without loss of generality we can assume that $h$ is given by the Asplund sum of $f,g$, that is 
  \[
    h(z) = 
    \sup_{
      \lambda x + (1-\lambda)y = z
    }
    f(x)^\lambda g(y)^{1-\lambda}. 
  \]
  In addition, up to replacing $f,g,h$ by 
  \[
    \bar f \eqdef \frac{f}{\norm{f}_{L^1}}, \qquad
    \bar g \eqdef \frac{g}{\norm{g}_{L^1}}, \qquad 
    \bar h \eqdef \frac{h}{
      \norm{f}_{L^1}^\lambda
      \norm{g}_{L^1}^{1-\lambda}
      },
  \]
  we can assume that $\displaystyle \int_{\mathbb{R}^n} f = \int_{\mathbb{R}^n} g = 1$, so that it suffices for us to show that $\displaystyle \int_{\mathbb{R}^n} h \ge 1$.

  Now define the following sets
  \[
    F_t \eqdef \{ f > t\}, \qquad
    G_t \eqdef \{ g > t\}, \qquad
    H_t \eqdef \{ h > t\},
  \]
  so that the normalization condition on $f,g$ and the layer-cake formula {\color{black}$\int_{\R^n}f\,\dd x=\int_0^\infty|F_t|\,\dd t$} imply that 
  \[
      \mu_f \eqdef |F_t|\dd t, \qquad 
      \mu_g \eqdef |G_t|\dd t
  \]
  are both probability measures over $\mathscr{P}(\mathbb{R})$. On the other hand, the Prékopa condition implies that 
  \begin{equation}\label{eq.PL_level_set_inclusion}
    \lambda F_r + (1-\lambda) G_s 
    \subseteq 
    H_{r^\lambda s^{1-\lambda}}.
  \end{equation}

  Since both measures $\mu_f, \mu_g$ are atomless, we can consider the Brenier map $T$ performing the $L^2$ optimal transport $T_\sharp \mu_f = \mu_g$. As a result, it solves the following Monge-Ampère equation 
  \begin{equation}\label{eq.MA_levelset_transport}
      |F_t| = |G_{T(t)}| T'(t), \quad \mu_f-a.e.
  \end{equation}
  Now, the Brunn-Minkowski inequality applied to~\eqref{eq.PL_level_set_inclusion} with $r = t$ and $s = T(t)$, alongside the Monge-Ampère equation above gives 
  \begin{align*}
    |H_{t^\lambda T(t)^{1-\lambda}}|
    &\ge 
    {
      \left(
        \lambda |F_t|^{1/n} + (1-\lambda)|G_{T(t)}|^{1/n}
      \right)
    }^n\\
    &= 
    |F_t|
    {
      \left(
        \lambda + (1-\lambda){T'(t)}^{-1/n}
      \right)
    }^n. 
  \end{align*}

  Focusing now on the integral of $h$, in order to profit from the estimate above, it is convenient to define the following change of variables induced by the level-set transportation via the one-dimensional Brenier map
  \begin{equation}\label{eq.change_variables_Brenier_map}
    u(t) \eqdef t^\lambda T(t)^{1-\lambda},
  \end{equation}
  as well as the Jacobian that it induces by change of variables
  \begin{equation}\label{eq.change_variables_jacobian_Brenier_map}
    \begin{aligned}
      u'(t)
      &= 
      \lambda t^{\lambda -1}T(t)^{1-\lambda} + 
      (1-\lambda)t^\lambda T(t)^{-\lambda}T'(t),\\ 
      &=
      \frac{u(t)}{t}\left[
        \lambda + (1-\lambda)\frac{t T'(t)}{T(t)}
      \right].
    \end{aligned}
  \end{equation}

  Indeed, the layer-cake combined with the change of variables formula imply that 
  \begin{align*}
    \int_{\mathbb{R}^n} h 
    &= 
    \int_0^{+\infty} |H_t|\dd t
    = 
    \int_0^{+\infty} |H_{u(t)}||u'(t)| \dd t \\
    &\ge 
    \int_0^{+\infty} 
    |F_t|
    {
      \left(
        \lambda + (1-\lambda){T'(t)}^{-1/n}
      \right)
    }^n
    |u'(t)| \dd t.
  \end{align*}

  We can equivalently rewrite the integrand of the right hand-side above as 
  \begin{align*}
    &{
      \left(
        \lambda + (1-\lambda){T'(t)}^{-1/n}
      \right)
    }^n
    u'(t)\\
    &\quad
    =  
    {
      \left(
        \lambda + (1-\lambda){T'(t)}^{-1/n}
      \right)
    }^n
    \frac{u(t)}{t}\left[
        \lambda + (1-\lambda)\frac{t T'(t)}{T(t)}
      \right]\\
    &\quad =
     {
      \left(
        \lambda {T'(t)}^{\frac{1-\lambda}{n}} 
        + 
        (1-\lambda){T'(t)}^{-\frac{\lambda}{n}}
      \right)
    }^n 
    \left(
        \lambda
        {
          \left(
            \frac{T(t)}{tT'(t)}
          \right)
        }^{1-\lambda} 
        + 
        (1-\lambda){  
          \left(
            \frac{tT'(t)}{T(t)}
          \right)
        }^{\lambda}
      \right)
  \end{align*}

  It is useful to write the two factors in the last equality under the following notation 
  \begin{equation}
    \begin{aligned}
      \dratio(t) 
      &\eqdef 
      {
        \left(
          \lambda {T'(t)}^{\frac{1-\lambda}{n}} 
          + 
          (1-\lambda){T'(t)}^{-\frac{\lambda}{n}}
        \right)
      }^n - 1,\\
      \did(t) 
      &\eqdef 
      \left(
          \lambda
          {\left(
              \frac{T(t)}{tT'(t)}
            \right)
          }^{1-\lambda} 
          + 
          (1-\lambda){  
            \left(
              \frac{tT'(t)}{T(t)}
            \right)
          }^{\lambda}
      \right)-1.
    \end{aligned}
  \end{equation}
  From the weighted arithmetic/geometric means inequality, we have for all $t$ that $\dratio(t), \did(t) \ge 0$. On the other hand, equality holds for $\dratio(t) = 0$ if and only if we have 
  \[
    T'(t)^{\frac{1-\lambda}{n}} = T'(t)^{\frac{-\lambda}{n}}
    \Longleftrightarrow 
    T'(t) = 1,
  \]
  while for the second we have $\did(t)= 0$ if and only if 
  \[
    {\left(
      \frac{T(t)}{tT'(t)}
    \right)
    }^{1-\lambda} 
    = 
    {  
      \left(
        \frac{tT'(t)}{T(t)}
      \right)
    }^{\lambda}
    \Longleftrightarrow 
    T(t) = tT'(t). 
  \]
  As a result, the first quantity is a pointwise measure of how far from being a small translation the transport of the $t$ level set of $f$ to $g$, while in conjunction with the second deficit, it gives a measure of how far the transportation is from not moving, that is $T(t) = t$.

  So coming back to the estimate, we have that 
  \begin{align*}
    \int_{\mathbb{R}^n} h 
    &\ge 
    \int_0^{+\infty} 
    |F_t|\left(
      1 + \dratio(t)
    \right)\left(
      1 + \did(t)
    \right) \dd t
    \ge
    \int_0^{+\infty} 
    |F_t| \dd t
    = 
    \int_{\mathbb{R}^n} f = 1.
  \end{align*}
This finishes the proof. 
\end{proof}

\subsection{The equality cases}\label{sec.qualitative_equality_cases}

The qualitative argument also allows us to recover the equality cases identified by Dubuc~\cite{dubuc1977criteres}. For this, we first prove a lemma showing that the translations of the level sets furnished by Brunn--Minkowski can be chosen independently of the height. It also serves as a warm-up to the quantitative proof.

\begin{lemma}\label{lemma.fixed_translation}
    Fix $\lambda \in (0,1)$, let $I \subseteq (0,+\infty)$ be an open interval, and consider a triplet of families of Borel sets of finite volume $A_t, B_t, C_t \subset \mathbb{R}^n$, which are increasing in $t \in I$ and such that 
    \[
      \lambda A_r + (1 - \lambda) B_s 
      \subseteq C_{\lambda r + (1 - \lambda)s} 
      \quad 
      \text{ for all }
      r,s \in I. 
    \]
    If for a.e. $t \in I$ we have 
    \[
      |A_t| = |B_t| = |C_t| > 0,
    \]
    then there is a vector $z \in \mathbb{R}^n$ such that 
    \[
      B_t = A_t + z, \quad 
      C_t = A_t + (1-\lambda)z. 
    \]
    In addition, the family can be chosen so that 
    \[
      \lambda A_r + (1-\lambda)A_s 
      \subseteq A_{\lambda r + (1-\lambda)s}.
    \]
\end{lemma}
\begin{proof}
    We start by fixing the following set of full measure in $I$
    \begin{equation}\label{eq.full_measure_set}
      I_0 \eqdef 
      \left\{
        t \in I: \quad 
        \substack{
          \displaystyle
          \lambda A_t + (1 - \lambda) B_t \subseteq C_{t}\\ 
          \displaystyle
          0 < |A_t| = |B_t| = |C_t| \eqdef m_t < +\infty
        }
      \right\}.
    \end{equation}
    Then from the Brunn--Minkowski inequality for $t \in I_0$ we have that 
    \[
        m_t = |C_t| \ge |\lambda A_t + (1 - \lambda) B_t| 
        \ge 
        {\left(
          \lambda {|A_t|}^{1/n} + (1-\lambda) {|B_t|}^{1/n}
        \right)}^n
        = m_t. 
    \]
    From the equality case of the Brunn--Minkowski inequality, then all three sets are convex bodies and there exists a vector $z_t \in \mathbb{R}^n$ such that 
    \begin{equation}\label{eq.variable_translation}
      B_t = z_t + A_t 
      \text{ and }
      C_t = \lambda A_t + (1 - \lambda)B_t = A_t + (1 - \lambda)z_t,
    \end{equation}
    where the second equality follows from all sets having all the same volume and the inclusion $\lambda A_t + (1 - \lambda) B_t \subseteq C_{t}$. As a result, to finish the proof it suffices to prove that the mapping $t \mapsto z_t$ is constant for a.e. $t \in I_0$. 

    Fixed $\theta \in \mathbb{S}^{n-1}$, we define the following functions, whose dependence on $\theta$ is omitted:
    \[
      a(t) \eqdef \sup_{ x \in A_t} \theta\cdot x, \qquad 
      b(t) \eqdef \sup_{ x \in B_t} \theta\cdot x, \qquad 
      c(t) \eqdef \sup_{ x \in C_t} \theta\cdot x.
    \]
    These functions coincide with the support function of their corresponding convex bodies, evaluated at the vector $\theta$. Thanks to~\eqref{eq.variable_translation} it holds that 
    \[
      a(t) = c(t) - (1-\lambda) \theta\cdot z_t, \quad 
      b(t) = c(t) + \lambda \theta\cdot z_t.
    \]
    On the other hand, for every $r, s \in I_0$ it follows from the set inclusion $\lambda A_r + (1 - \lambda) B_s \subseteq C_{\lambda r + (1 - \lambda)s}$ that 
    \begin{equation}\label{eq.inclusion_translation_support_inequalities}
      \begin{aligned}
        \lambda c(s) + (1 - \lambda)c(r) + 
        \lambda(1 - \lambda)\theta\cdot(z_r - z_s)
        &\le 
        c(\lambda s + (1-\lambda)r),\\
        (1 - \lambda) c(s) + \lambda c(r) + 
        \lambda(1 - \lambda)\theta\cdot(z_s - z_r)
        &\le 
        c((1-\lambda) s + \lambda r). 
      \end{aligned}
    \end{equation}
    Summing both of these inequalities, we get that 
    \begin{equation}\label{eq.almost_concavity_estimate}
      c(r) + c(s) \le 
      c((1-\lambda) s + \lambda r) + c(\lambda s + (1-\lambda)r).
    \end{equation}

    The above estimate suggests that $t \mapsto c(t)$ is actually a concave function, and our goal is to show this, which will require a somewhat delicate argument. For this, it suffices to show that $c$ is midpoint concave. In other words, that the usual concavity inequality holds with the convex interpolation parameter $1/2$. 

    Consider the mapping 
    \[
      I_0 \times I_0 \ni (s,r) \mapsto 
      (s', r') \eqdef (\lambda s + (1 - \lambda r), (1- \lambda)s + \lambda r),
    \]
    and define, for now at a formal level, the following sequence by induction 
    \[
      (s_{k +1}, r_{k +1}) 
      \eqdef 
      (s_k', r_k'), \text{ and }
      (s_{0}, r_{0}) 
      \eqdef 
      (s, r),
    \]
    for some fixed $(s,r) \in I_0\times I_0$. An elementary argument shows that
    \[
      \lim_{k \to \infty} s_k = \lim_{k \to \infty} r_k
      = 
      \frac{s + r}{2}. 
    \]
    
    As a result, if the sequence ${(s_k, r_k)}_{k \in \mathbb{N}}$ is contained in $I_0 \times I_0$, we can apply estimate~\eqref{eq.almost_concavity_estimate} for every $k$ and obtain that
    \begin{align*}
      c(s) + c(r) 
      &= 
      c(s_0) + c(r_0) 
      \le 
      \lim_{k \to \infty}
      c(s_k) + c(r_k) \\
      &\le 
      \lim_{k \to \infty}
      c(s_{k+1}) + c(r_{k+1}) 
      = 
      2 c\left(
        \frac{s + r}{2}
      \right),
    \end{align*}
    whenever $c$ is continuous over $\displaystyle \frac{s + r}{2}$. To obtain the concavity of $c$ then, we start by choosing a subset of full measure of $I_0$ with good properties, since then the concavity will follow with a continuous extension. 

    First, notice that $t \mapsto c(t)$ is a monotone function, because of the set inclusion $C_s \subset C_t$ for $s \le t$, so it is almost everywhere continuous. Consider now 
    \[
      I'_0 \times I_0' \eqdef 
      \left\{
        (s,r) \in I_0 \times I_0 : 
        \substack{
          \displaystyle
          (s', r') \in I_0 \times I_0, \\ 
          \displaystyle
          c \text{ is continuous at } \frac{s + r}{2}
        } 
      \right\}.
    \]
    Since $\lambda \in (0,1)$, the operation $(s,r) \mapsto (s',r')$ is Lipschitz, so the set $I_0'\times I_0'$ has full measure in $I_0\times I_0$. Applying this operation countably many times, we obtain a set of full measure  
    \[
       J_\star \subset I_0 \times I_0, \quad 
       |I_0 \times I_0 \setminus J_\star| = 0,
       \text{ and }
       J_\star \subset 
      \left\{
        (s,r) \in I_0 \times I_0 : 
        \substack{
          \displaystyle
          {(s_k, r_k)}_{k \in \mathbb{N}} \subset I_0 \times I_0, \\ 
          \displaystyle
          c \text{ is continuous at } \frac{s + r}{2}
        }
      \right\}.
    \]
    From the symmetry of the above operations, we have that $J_\star$ is of the form $I_\star \times I_\star$. Considering first $(s,r) \in I_\star \times I_\star$, followed by a continuous extension argument, shows that $c$ is midpoint concave, and hence concave. 

    Now coming back to the estimates~\eqref{eq.inclusion_translation_support_inequalities}, we notice that we can estimate the distance between $z_t - z_s$ by the Jensen divergence of $c$, from which we can now obtain more information using concavity. Indeed, using the notation above of $(s',r')$, and defining the Jensen divergence of $c$ as
    \[
      J_\lambda(s,r) 
      \eqdef 
      c(s') - \lambda c(s) - (1-\lambda)c(r),
    \]
    we obtain from the concavity of $c$ for almost every $s,r$ that $c$ is differentiable at $s,r$ and then 
    \[
      c(s') = c(r) + \lambda (s-r) c'(r) + O(|s-r|^2), \quad 
      c(s) = c(r) + (s-r) c'(r) + O(|s-r|^2). 
    \]
    Therefore the zero-th and first order terms cancel each other in $J_\lambda(s,r)$ and $J_\lambda(r,s)$, and both are of the order $O(|s-r|^2)$. Since $J_\lambda(s,r), J_\lambda(r,s) \ge 0$ from concavity, we obtain from~\eqref{eq.inclusion_translation_support_inequalities} that 
    \[
      \lambda(1-\lambda)|\theta \cdot (z_r - z_s)| 
      \le 
      \max\{
        J_\lambda(s,r), J_\lambda(r,s)
      \}
      = O(|s-r|^2). 
    \]
    Since the right-hand side above is independent of $\theta$, and it was an arbitrary element of $\mathbb{S}^{n-1}$, this can only be the case if $t \mapsto z_t$ is constant.
\end{proof}

\begin{corollary}\label{cor.equality}
Let $\lambda\in(0,1)$ and let $f,g,h:\R^n\to[0,\infty)$
be measurable functions satisfying the Prékopa condition. Assume that
\[
  m_f\eqdef\int_{\R^n}f \dd x \in(0,\infty),\qquad
  m_g\eqdef\int_{\R^n}g \dd x \in(0,\infty),
\]
and set $m_\lambda\eqdef m_f^\lambda m_g^{1-\lambda}$.
Then equality holds in the Prékopa--Leindler inequality,
namely
\[
  \int_{\R^n}h \dd x = m_\lambda,
\]
if and only if there exist a nonzero integrable log-concave
function $\varrho$ and a vector $z\in\R^n$ such that,
almost everywhere,
\begin{equation}\label{eq.equality_cases_cor}
  h=\varrho,\qquad
  f=\frac{m_f}{m_\lambda}
    \varrho\bigl(\cdot+(1-\lambda)z\bigr),\qquad
  g=\frac{m_g}{m_\lambda}
    \varrho(\cdot-\lambda z).
\end{equation}
\end{corollary}
\begin{proof}
To prove the equality cases we assume once again for simplicity that the functions are normalized. If we have $\norm{h}_{L^1} = 1$, then the previous estimates from Theorem~\ref{thm.qualitative_PL} give that 
  \[
    \int_0^{+\infty} 
    |F_t|
    \left[
      \left(
      1 + \dratio(t)\did(t) + \did(t)
    \right)\left(
      1 + \did(t)
    \right) - 1
    \right]
    \dd t = 0.
  \]
  As a result, since both $\dratio$ and $\did$ are non-negative, this can only be the case if 
  \[
    \dratio(t) = \did(t) = 0
  \]
  for $\mu_f$-almost every $t > 0$. This also yields that 
  \[
    0 = \int_{0}^\infty (|H_t| - |F_t|) \dd t,
  \] 
  and both of these together imply that for almost every $t$ we have $T = \id$, so  $\mu_f = \mu_g$ and we conclude that 
  \[
    |H_t| = |F_t| = |G_t|.
  \]
  
  Now passing to the logarithmic change of variables $t(q) = e^{-q}$, so that $q = -\log t$ and defining 
  \[
    E_q \eqdef E_{t(q)}, \quad \text{ for } E \in \{F,G,H\},
  \]
  we obtain from the Prékopa condition that 
  \[
    \lambda F_p + (1-\lambda) G_q \subseteq H_{(1-\lambda)q + \lambda p}
  \]
  for every $p,q \in \mathbb{R}\cup\{+\infty\}$. As a result, the assumptions of Lemma~\ref{lemma.fixed_translation} are satisfied, so there exists a vector $z \in \mathbb{R}^n$ such that 
  \[
    G_q = z + F_q, \quad 
    H_q = (1-\lambda)z + F_q,
  \]
  and all these level sets are convex bodies. In addition, the final condition of Lemma~\ref{lemma.fixed_translation} applied in the log-scale implies that $f,g,h$ are log-concave. 
  
  To verify that they are given by translations of the same log-concave function, notice that for any smooth test function $\varphi$ the layer-cake formula gives 
  \[
    \int_{\mathbb{R}^n} \varphi h \dd x 
    = 
    \int_{\mathbb{R}} \int_{H_t} \varphi(x) \dd x \dd t
    = 
    \int_{\mathbb{R}} \int_{F_t - (1-\lambda)z} \varphi(x) \dd x \dd t
    = 
    \int_{\mathbb{R}^n} \varphi f(\cdot - (1-\lambda)z) \dd x
  \]
  which gives the desired translation identities. 
\end{proof}

\section{The quantitative height transport proof}\label{sec.quantitativePL}

Now we move on to making the level set proof from Section~\ref{sec.qualitative_proof} quantitative. Recall that from Corollary~\ref{corollary.quant_mid_point_implies_full}, it suffices to consider the midpoint case $\lambda = 1/2$. Therefore, throughout this section we consider $f,g,h$ satisfying the Prékopa condition, and we may also assume that 
\begin{equation}\label{eq.quant_normalization}
 \int f=\int g=1,\qquad
 h((x+y)/2)\geq\sqrt{f(x)g(y)}.
\end{equation} 
Another possible simplification is replacing $h$ by the Asplund midpoint sum of $f,g$ that is 
\[
  h(z) = \sup_{
    (x + y)/2 = z
  }
  \sqrt{f(x)g(y)}. 
\]
This replacement will only decrease the Prékopa-Leindler deficit defined below. In addition the resulting $h$ is log-concave, and
$\|h\|_\infty=\sqrt{\|f\|_\infty\|g\|_\infty}$.

As a result, the Prékopa-Leindler inequality is equivalent to the non-negativity of the following deficit 
\begin{equation}\label{eq.PL_deficit}
  \dPL = 
  \dPL(f,g)
  \eqdef \int_{\mathbb{R}^n} h -1,
\end{equation}
and we recall that our goal is to show that there exists a dimensional constant $C_n$ such that for all $f,g$ log-concave it holds that 
\begin{equation}
    \inf_{z \in \mathbb{R}^n}\norm{f - g(\cdot + z)}^2_{L^1(\mathbb{R}^n)} 
    \le 
    C_n \dPL(f,g). 
\end{equation}
The zero-deficit case was proved in Section~\ref{sec.qualitative_equality_cases}.
Since the distance between probability densities is at most two, up to increasing $C_n$ it suffices to treat $0<\dPL\leq1$. Throughout the quantitative proofs, $C$ denotes an absolute constant, which may change from line to line; the dimensional constants from Brunn--Minkowski will be kept explicit.

\subsection{The change of variables to logarithmic scale}\label{sec.change_variables}
Our strategy towards this quantitative statement is to interpret the epigraphs as convex bodies and apply the by now well understood stability of the Brunn-Minkowski inequality. Of course, in Euclidean coordinates the epigraphs of $f,g,h$ are not convex, so, for that reason, before proceeding we establish some notation to properly compare them in a logarithmic scale. 

At this scale we work with the convex potentials associated with $f,g,h$. In addition, interchanging $f,g$ if necessary, assume that $\|f\|_\infty\geq\|g\|_\infty$, and 
we set the following notation to a common reference height:
\begin{equation}\label{eq.common_height}
 \begin{gathered}
 M=\|g\|_\infty=\min\{\|f\|_\infty,\|g\|_\infty\},\qquad f=Me^{-V},\quad g=Me^{-W},\quad h=Me^{-U},\\
 \min V=-m,\qquad \min W=0,\qquad m=\log\frac{\|f\|_\infty}{\|g\|_\infty},
 \qquad \min U=-m/2.
 \end{gathered}
\end{equation}
Since $V,W,U$ are convex potentials such that $0 < \int e^{-V}, \int e^{-W}, \int e^{-U} < +\infty$, they are coercive, have bounded sublevels and
attain their minima. For $q\in\R$ let
\begin{equation}\label{eq.log_levels}
 \begin{aligned}
 A_q&=\{f>Me^{-q}\}=\{V<q\},\\
 B_q&=\{g>Me^{-q}\}=\{W<q\},\\
 C_q&=\{h>Me^{-q}\}=\{U<q\}.
 \end{aligned}
\end{equation}
Write $a(q)=|A_q|$, $b(q)=|B_q|$. These two volume profiles recur throughout the proof; in particular $a=0$ on $(-\infty,-m]$ and $b=0$ on $(-\infty,0]$. In fact, these tail of these distributions inherit the geometry of the log-concave functions defining them, which is discussed in the sequel, where we prove a few fundamental properties of the height distributions. 
\begin{lemma}\label{lemma.tail_ab}
    For $i \in \{0,1\}$, let $\varrho_i \in \mathscr{P}(\mathbb{R})$ be such that
    \[
      \dd \varrho_i = M a_i(q) e^{-q} \dd q, 
    \]
    so that $\varrho_0 = \mu$ and $a_0 = a$, while $\varrho_1 = \nu$ and $a_1 = b$. Then, setting $m_0=-m$ and $m_1=0$, the minimal values of the potentials $V, W$, associated with $\varrho_0, \varrho_1$ respectively, the following properties hold.
    \begin{itemize}
      \item[(a)] For $m_i < t < q$ we have
      \[
          a_i(q) \le {\left(
            \frac{q - m_i}{t - m_i}
          \right)}^n a_i(t),
      \]   
      for both $i = 0,1$.
      \item[(b)] Both measures $\varrho_0, \varrho_1 \in \mathscr{P}(\mathbb{R})$ are log-concave, hence their displacement interpolation given by
      \[
        \varrho_t \eqdef (\mathcal{T}_t)_{\sharp} \mu, \quad 
        \mathcal{T}_t \eqdef 
        (1-t)\id + t \mathcal{T}
      \]
      remains log concave for all $t \in [0,1]$. 
      \item[(c)] Identifying $\varrho_t$ with its density and writing
      $m_t \eqdef (1-t)m_0+tm_1 = -(1-t)m$, then for a.e. $q>m_t$
      \[
        -1\leq(\log\varrho_t)'(q)
        \leq\frac{n+1}{q-m_t}.
      \]
      Consequently, for every $q>m_t$,
      \begin{equation}\label{eq.lemma_interp_1st_identity}
        \varrho_t(q+s)\geq e^{-s}\varrho_t(q),
        \qquad s\geq0,
      \end{equation}
      and
      \begin{equation}\label{eq.lemma_interp_2nd_identity}
        \varrho_t(q-s)\geq
        \left(\frac{q-s-m_t}{q-m_t}\right)^{n+1}
        \varrho_t(q),
        \qquad 0\leq s<q-m_t.
      \end{equation}
    \end{itemize} 
\end{lemma}
\begin{proof}
  For $i \in \{0,1\}$, let $V_i$ be the corresponding convex potential generating $a_i$ in the above construction, along with the level sets $A_{i,q} \eqdef \{V_i(x) < q\}$. So we have $V_0 = V$ and $A_{0,q} = A_q$, while $V_1 = W$ and $A_{1,q} = B_q$ in the above notation. 

  We start with item (a). Let $x_\star$ be a minimizer of $V_i$, then we can show that for $m_i<t<q$ the following inclusion holds 
  \[
    x_\star + \frac{t - m_i}{q - m_i}(A_{i,q} - x_\star) \subset A_{i,t}. 
  \]
  Indeed, considering some $x \in A_{i,q}$ we have from the convexity of $V_i$ that
  \begin{align*}
    V_i\left(
      \frac{q-t}{q-m_i} x_\star + \frac{t - m_i}{q - m_i} x
    \right)
    \le 
    \frac{q-t}{q-m_i} V_i(x_\star) + \frac{t - m_i}{q - m_i} V_i(x) 
    < \frac{q-t}{q-m_i} m_i + \frac{t - m_i}{q - m_i} q = t. 
  \end{align*}
  Now using the monotonicity of the Lebesgue measure on the above inclusion, it follows that 
  \[
    {\left(
      \frac{t - m_i}{q - m_i}
    \right)}^{n} a_i(q) 
    \le a_i(t).
  \]

  For item (b), first notice that 
  \[
    q \mapsto a_i(q)^{1/n} \text{ is a concave function }
  \]
  thanks to the Brunn-Minkowski inequality. As a result, we can rewrite 
  \[
    \log \varrho_i = 
    \log M - q + n \log a_i^{1/n},
  \]
  which corresponds to a composition of concave functions, and hence concave. Although the log-concavity of displacement interpolations is not generally preserved along the interpolation, even if the marginals enjoy such geometry, see~\cite{santambrogio2016convexity} for many counter-examples, it is nonetheless known to be true for one-dimensional distributions, but for completeness we give a proof here that is based on the particular structure of geodesics in $\mathbb{R}$ as well as a criterion of log-concavity by Bobkov~\cite{bobkov1996extremal} in this setting. Given $\varrho \in \mathscr{P}(\mathbb{R})$ we define its distribution and quantile functions as 
  \[
    G_\varrho(q) 
    \eqdef 
    \varrho((-\infty,q]), \quad 
    Q_\varrho(u) 
    \eqdef 
    G_\varrho^{-1}(u).  
  \]
  As a result, we naturally have $G_\varrho \circ Q_\varrho = \id$ and the chain rule implies that 
  \[
    I_\varrho(u) 
    \eqdef 
    \varrho(Q_\varrho(u)) 
    = 
    \frac{1}{Q'_\varrho(u)}. 
  \]
  It was then proven by Bobkov, see also~\cite[Prop.~2.22]{saumard2014logconcavity_review}, that $\varrho$ is log-concave if and only if $I_\varrho$ is concave over $(0,1)$. 

  On the other hand, we may represent any one-dimensional distribution as 
  \[
    \varrho = {(Q_\varrho)}_\sharp \mathscr{L}^1\mres [0,1]. 
  \]
  Hence a routine computation for the geodesics 
  \[
  \begin{aligned}
    \varrho_t 
    &= 
    {
      ((1-t)\id + t\mathcal{T})
    }_\sharp \varrho_0 
    = 
    {
      ((1-t)Q_0 + t\mathcal{T}\circ Q_0)
    }_\sharp \mathscr{L}^1\mres [0,1] \\
    &= 
    {
      ((1-t)Q_0 + t Q_1)
    }_\sharp \mathscr{L}^1\mres [0,1].
  \end{aligned}
  \]
  implies that the quantile of geodesic interpolations is the interpolation of quantiles $Q_t = (1-t) Q_0 + t Q_1$. As a consequence for Bobkov's function $I_t$ associated to the interpolation we have 
  \[
    I_t(u) = \frac{1}{Q'_t(u)} 
    = 
    {
      \left(
        (1-t) Q_0'(u) + t Q_1'(u)
      \right)
    }^{-1}
    = 
    {
      \left(
        \frac{1-t}{I_0(u)} + \frac{t}{I_1(u)}
      \right)
    }^{-1}.
  \]
  Since the marginals are both log-concave, the necessity of Bobkov's criterion gives that $I_0, I_1$ are concave, while composing with the harmonic mean gives the concavity of $I_t$ and implies the log-concavity of $\varrho_t$ thanks to the sufficiency part of the same statement. 

  We finally prove item (c). First, since $a_i$ is non-decreasing, for every $q > m_i$ it holds that
  \[
    1-G_{\varrho_i}(q)
    =M\int_q^{+\infty}e^{-s}a_i(s)\dd s
    \geq M a_i(q)\int_q^{+\infty}e^{-s}\dd s
    =\varrho_i(q).
  \]
  On the other hand, item (a) and the inequality $e^{-s}\geq e^{-q}$ for $m_i<s<q$ give
  \begin{align*}
    G_{\varrho_i}(q)
    &=M\int_{m_i}^q e^{-s}a_i(s)\dd s\\
    &\geq Me^{-q}a_i(q)
    \int_{m_i}^q
    {\left(
      \frac{s-m_i}{q-m_i}
    \right)}^n\dd s
    =\frac{q-m_i}{n+1}\varrho_i(q).
  \end{align*}
  Evaluating these two inequalities at $q=Q_i(u)$, and recalling that $I_i(u)=\varrho_i(Q_i(u))$, we obtain
  \[
    I_i(u)\leq1-u,\qquad
    \frac{u}{I_i(u)}\geq\frac{Q_i(u)-m_i}{n+1}.
  \]
  Both estimates are preserved by the displacement interpolation. Indeed, the formula for $I_t$ obtained above gives
  \[
    I_t(u)\leq1-u
  \]
  and
  \begin{align*}
    \frac{u}{I_t(u)}
    &=(1-t)\frac{u}{I_0(u)}+t\frac{u}{I_1(u)}\\
    &\geq
    \frac{(1-t)(Q_0(u)-m_0)+t(Q_1(u)-m_1)}{n+1}
    =\frac{Q_t(u)-m_t}{n+1}.
  \end{align*}
  Taking $u=G_{\varrho_t}(q)$ in the previous inequalities, we conclude that
  \begin{equation}\label{eq.interpolated_distribution_bounds}
    \frac{\varrho_t(q)}{1-G_{\varrho_t}(q)}\leq1,
    \qquad
    \frac{\varrho_t(q)}{G_{\varrho_t}(q)}
    \leq\frac{n+1}{q-m_t}.
  \end{equation}

  We now use the log-concavity proved in item (b). At every differentiability point $q>m_t$ where $(\log\varrho_t)'(q)<0$, concavity gives
  \[
    \varrho_t(s)\leq
    \varrho_t(q)e^{(\log\varrho_t)'(q)(s-q)}
    \qquad \text{ for all } s\geq q.
  \]
  Integrating this inequality and using the first estimate in~\eqref{eq.interpolated_distribution_bounds}, we obtain
  \[
    \varrho_t(q)\leq1-G_{\varrho_t}(q)
    \leq\frac{\varrho_t(q)}{-(\log\varrho_t)'(q)},
  \]
  which implies $(\log\varrho_t)'(q)\geq-1$. Similarly, if $(\log\varrho_t)'(q)>0$, concavity gives
  \[
    G_{\varrho_t}(q)\leq
    \frac{\varrho_t(q)}{(\log\varrho_t)'(q)}.
  \]
  Combining this inequality with the second estimate in~\eqref{eq.interpolated_distribution_bounds}, we obtain
  \[
    (\log\varrho_t)'(q)\leq\frac{n+1}{q-m_t}.
  \]
  The case $(\log\varrho_t)'(q)=0$ is immediate, and this proves the first conclusion of item (c).

  Finally, integrating the lower bound for $(\log\varrho_t)'$ between $q$ and $q+s$ gives
  \[
    \varrho_t(q+s)\geq e^{-s}\varrho_t(q).
  \]
  Likewise, integrating the upper bound between $q-s$ and $q$, for $0\leq s<q-m_t$, gives
  \[
    \log\frac{\varrho_t(q)}{\varrho_t(q-s)}
    \leq(n+1)\log\frac{q-m_t}{q-s-m_t},
  \]
  which is equivalent to the last conclusion of item (c).
\end{proof}
As seen in the qualitative proof, another advantage of the log-scale is that the Asplund sum of log-concave functions is translated into a Minkowski sum for all level sets. Indeed, we have that
\begin{equation}\label{eq.cross_levels_midpoint}
 \frac{A_r+B_s}{2}\subset C_{(r+s)/2},\qquad
 U((x+y)/2)\leq\frac{V(x)+W(y)}2,
\end{equation}
while the change of variables for the total mass of $f,g$ gives
\begin{equation}
  M\int_\R e^{-q}a(q)\dd q =M\int_\R e^{-q}b(q)\dd q=1,\label{eq.common_masses}
\end{equation}
and Cavalieri's principle give the following formula for representing 
\begin{align}
 \|f-g(\cdot+z)\|_1
 =M\int_\R e^{-q}|A_q\sym(B_q-z)|\dd q.
 \label{eq.log_layer_cake}
\end{align}

Finally, the central objects on the qualitative proof of Theorem~\ref{thm.qualitative_PL} are the height distributions and the Brenier map $T$ transporting them, denoted in Euclidean heights as 
\[
 \mu \eqdef |F_t|\dd t, \qquad 
  \nu \eqdef |G_t|\dd t, \qquad 
 T_\sharp \mu=\nu. 
\]
Performing the change of variables $t = t(q) = Me^{-q}$ above, the height distribution in the common logarithmic coordinate become
\begin{equation}\label{eq.log_transport_measures}
 \dd\mu(q)=Me^{-q}a(q)\dd q,\qquad
 \dd\nu(q)=Me^{-q}b(q)\dd q.
\end{equation}
While the same map in log-scale can be written as $\mathcal T$, having the property
$F_t=A_q$ and $G_{T(t)}=B_{\mathcal T(q)}$ when $t=Me^{-q}$, becomes 
\begin{equation}\label{eq.log_transport_definition}
 \mathcal T(q) \eqdef 
 t^{-1}\circ T(t(q)) 
 =
 -\log\frac{T(Me^{-q})}{M},\qquad
 T(Me^{-q})=Me^{-\mathcal T(q)}.
\end{equation}
The common change of variables also gives $\mathcal T_\sharp \mu=\nu$.
Both changes of coordinate reverse orientation, so $\mathcal T$ is
increasing and therefore corresponds to the Brenier map transporting $\mu$ to $\nu$, see~\cite[Chap.~2]{santambrogio2015optimal}. Its derivative and displacement satisfy, for $t=Me^{-q}$,
\begin{equation}\label{eq.log_transport_identities}
 \mathcal T'(q)=\frac{tT'(t)}{T(t)},\qquad
 \mathcal T(q)-q=-\log\frac{T(t)}t.
\end{equation}

These identities for the transport map in logarithmic scale allow us to better understand the deficits that naturally appeared in the qualitative proof~\ref{thm.qualitative_PL}. Specializing to the case $\lambda = 1/2$, we obtain from this proof that 
\begin{equation}\label{eq.quanli2quanti}
  \begin{aligned}
    \dPL 
    &\ge 
    \int_{0}^\infty 
    |F_t|\left(
        \dratio(t) + \did(t)
    \right) \dd t 
    = 
    \int_{\mathbb R} 
    \left(
        \dratio(q) + \did(q)
    \right) \dd \mu(q),
  \end{aligned}
\end{equation}
where we employ the slight abuse of notation $\dratio(q) \eqdef \dratio(t(q))$ and $\did(q) \eqdef \did(t(q))$, and the deficits can be written in our setting on Euclidean heights as  
\[
  \dratio(t) ={\left[
    \frac{T'(t)^{\frac{1}{2n}}}{2} + \frac{T'(t)^{-\frac{1}{2n}}}{2}
  \right]}^{n} - 1, \quad 
  \did(t) = 
  \frac{1}{2} 
  \left[
    {\left(
      \frac{T(t)}{t T'(t)}
    \right)}^{\frac{1}{2}} + {\left(
      \frac{T(t)}{t T'(t)}
    \right)}^{-\frac{1}{2}} 
  \right] - 1.
\]
Replacing the expressions of $T'(t)$ by the equivalent one $T'(t(q)) = T'(q) T(t(q))/t(q)$ obtained in~\eqref{eq.log_transport_identities} as well as $\frac{T(t(q))}{t(q)} = e^{q - \mathcal{T}(q)}$ in the deficits above, they can be rewritten in log scale as 
\begin{equation}\label{eq.log_scalar_deficits}
 \begin{aligned}
 \dratio(q)
 &=\cosh\!\left(\frac{\log \mathcal{T}'(q) + (q - \mathcal{T}(q))}{2n}\right)^n-1,\\
 \did(q)
 &=\cosh\!\left(\frac{\log \mathcal{T}'(q)}{2}\right)-1.
 \end{aligned}
\end{equation}
Not only these coordinates allow us to see quite clearly how the deficit $\dPL$ penalizes the deviation of $\mathcal{T}'$ to $1$ and $\mathcal{T}$ from the identity, but they will also be useful while estimating the distance between the distributions $\mu, \nu$.

\subsection{Brunn--Minkowski deficit and the translation problem}\label{sec.brunn_minkowski}

The change of variables to logarithmic heights described above allows us to use the technology of the quantitative stability of Brunn-Minkowski inequality, due to the fact that this change of variables turns the Prékopa condition into a Minkowski sum. 

For convex bodies $E,F\subset\R^d$, we define the \emph{additive volume deficit}
\begin{equation}\label{eq.BM_deficit}
 \dBM(E,F)=\left|\frac{E+F}{2}\right|
       -\left(\frac{|E|^{1/d}+|F|^{1/d}}2\right)^d.
\end{equation}
The ambient dimension is the dimension of the bodies; it will be $n$
for individual levels, $n+1$ for the windows in this section, and $n+3$ for the lifted windows introduced later.

The stability principle of Figalli-Maggi-Pratelli~\cite{figalli2009refined,figalli2010mass} implies that this quantity controls the Fraenkel asymmetry between two convex bodies, namely 
\begin{equation}\label{eq.FMP_asymmetry}
 \mathcal A(E,F)^2
 \leq C_d
 \sigma_{E,F}
 \left[\frac{|(E+F)/2|^{1/d}}{(|E|^{1/d}+|F|^{1/d})/2}-1\right],
\end{equation}
where 
\begin{equation}
  \mathcal A(E,F)=\inf_{v \in \mathbb{R}^d}
  \left\{
    \frac{|E\sym(v+\rho F)|}{|E|} : \  
    |\rho F| = |E|
  \right\} \text{ and }
  \sigma_{E,F}
  \eqdef 
  \max\left\{\frac{|E|}{|F|},\frac{|F|}{|E|}\right\}^{1/d}. 
\end{equation}

The dilation $\rho$ is essential since a homothetic pair can attain equality
in Brunn--Minkowski, so the deficit $\dBM$ alone cannot force the unscaled volumes of level sets to agree. To distinguish a change of shape from a change of volume, define also
\begin{equation}\label{eq.scalar_volume_gap}
 \dvol(E,F)=\frac{|E|+|F|}{2}
       -\left(\frac{|E|^{1/d}+|F|^{1/d}}2\right)^d,
 \qquad d\geq2.
\end{equation}
Since the scalar volume deficit above depends only on $|E|, |F|$, we can use a slight abuse of notation of writing $\dvol(|E|,|F|)$ to be a function of scalar arguments for which the following equivalent bounds hold. This estimate is elementary, but it will be useful in comparing the height distributions $\mu, \nu$ to the Hellinger distance. 

\begin{lemma}\label{lemma.scalar_Hellinger}
For $p, q > 0$, the following bounds hold
\begin{equation}\label{eq.dvol_bounds}
 \frac14(\sqrt p-\sqrt q)^2\leq\dvol(p,q)
 \leq\frac12(\sqrt p-\sqrt q)^2.
\end{equation}  
\end{lemma}
\begin{proof}
  Starting with the lower bound, we recall the power-mean inequality, which is a fundamental consequence of Jensen's inequality. Since $d\geq2$, it gives
  \[
    \sqrt{pq}\leq
    \left(\frac{p^{1/d}+q^{1/d}}2\right)^d
    \leq\left(\frac{\sqrt p+\sqrt q}{2}\right)^2.
  \]
  Subtracting these inequalities from $(p+q)/2$ gives both bounds.
\end{proof}

As a result $\dBM$ controls a shape mismatch after volume normalization, whereas
$\dvol$ handles any dilation that might be introduced by a mismatch of volumes. We record in the following theorem, as a consequence of~\eqref{eq.FMP_asymmetry} that these two deficits control the $L^1$ norm between convex bodies, in a way that is well adapted to handle the distance between level sets. 
\begin{theorem}[Convex Brunn--Minkowski stability]
\label{thm.BM_stability}
For convex bodies $E,F\subset\R^d$ with $d\geq 2$,
\begin{equation}\label{eq.BM_stability}
 \inf_{v \in \mathbb{R}^d}
 |E\sym(F-v)|^2
 \leq B_d^{\rm conv}(|E|+|F|)
       [\dBM(E,F)+\dvol(E,F)].
\end{equation}
The constant can be chosen with $B_d^{\rm conv}\leq e^{Cd}$.
\end{theorem}
\begin{proof}
Our starting point is the relative asymmetry estimate of
Figalli--Maggi--Pratelli~\cite{figalli2009refined}; denote its constant
in~\eqref{eq.FMP_asymmetry} by $A_d$, and set $p=|E|$, $q=|F|$
and $\rho=(p/q)^{1/d}$. First we suppose that $p/q\in[1/2,2]$.
After translating $F$ so that it contains the origin, its dilates are
nested. Choosing a translation $v$ attaining the relative asymmetry gives
\[
 \inf_z|E\sym(z+F)|
 \leq p\mathcal A(E,F)+|p-q|.
\]
Writing $b=((p^{1/d}+q^{1/d})/2)^d\geq p/2$, the concavity of
$t^{1/d}$ gives
\[
 \frac{|(E+F)/2|^{1/d}}{b^{1/d}}-1
 \leq\frac{\dBM(E,F)}{db}.
\]
Combining this with~\eqref{eq.FMP_asymmetry} and
$|p-q|^2\leq8(p+q)\dvol(p,q)$ proves the desired estimate with
$B_d^{\rm conv}\leq C(1+A_d/d)$ in this case.

If the volumes are not comparable, the lower bound in
Lemma~\ref{lemma.scalar_Hellinger} gives
$\dvol(p,q)\geq c(p+q)$, and the bound by $p+q$ finishes the proof.
Finally, the constants in the proof of~\cite{figalli2009refined} satisfy
$A_d\leq e^{Cd}$, so the same is true of $B_d^{\rm conv}$.
\end{proof}

\subsection{Two estimates on the heights}

Our goal is to apply the quantitative version of Brunn-Minkowski discussed in the previous paragraph to suitable sections of the epigraphs of $f,g$. This naturally introduces translations, which are desirable to capture the equality cases of the Prékopa-Leindler inequality. These correspond to horizontal translations of the epigraphs, but simply applying Theorem~\ref{thm.BM_stability} also yields an undesired vertical translation, and hence we must be able to control both. 

In this section we will derive two estimates on the distances between the height distributions $\mu$ and $\nu$. The notion of distance being employed will be chosen to capture each one of these sources of translations; indeed, we will prove later on that the Wasserstein distance $W_2^2(\mu,\nu)$ controls horizontal translations of the epigraph, while an $L^1$ estimate controls height translations. This is in accordance with the intuition that optimal transport gives rise to \textit{horizontal distances}, while in the $L^p$ space we measure vertical deviations, see the discussion on the start of Chapter 5 of~\cite{santambrogio2015optimal}. 

\begin{lemma}\label{lemma.height_transport}
Let $\mu,\nu \in \mathscr{P}(\mathbb{R})$ be defined as in  \eqref{eq.log_transport_measures} and the map $\mathcal T$ defined in~\eqref{eq.log_transport_definition}. Then the following bound holds
\begin{equation}\label{eq.height_full_cost}
 \int_\R
 \left[|\mathcal T(q)-q|^2+|\log\mathcal T'(q)|^2\right]
       \dd\mu(q)\leq16(n+1)\dPL.
\end{equation}
In particular, since $\mathcal{T}$ is the Brenier map transporting $\mu$ onto $\nu$ we obtain that
\[
  W_2^2(\mu, \nu) \le16(n+1)\dPL. 
\]

In addition, setting $c(q)\eqdef|C_q|$, $\dd\mu_h(q)\eqdef Me^{-q}c(q)\dd q$, and defining the new transportation map 
\[
  \mathcal{S} \eqdef \frac{\id+\mathcal T}{2},
\]
there is a nonnegative measure $\xi$ such that
\begin{equation}\label{eq.output_height_measure}
 \mu_h=\mathcal{S}_\sharp\mu+\xi,\qquad \xi(\R)=\dPL.
\end{equation}
\end{lemma}
\begin{proof}
Recalling the expression for the deficits~\eqref{eq.log_scalar_deficits} obtained in section~\ref{sec.change_variables}
\begin{equation*}
  \begin{aligned}
 \dratio(q)
 &=\cosh\!\left(\frac{\log \mathcal{T}'(q) + (q - \mathcal{T}(q))}{2n}\right)^n-1,\\
 \did(q)
 &=\cosh\!\left(\frac{\log \mathcal{T}'(q)}{2}\right)-1,
 \end{aligned}
\end{equation*}
we have obtained by carefully tracking the deficits in the qualitative proof the expression in~\eqref{eq.quanli2quanti} which is rewritten here as
\begin{equation}\label{eq.log_scalar_budget}
 \int\left[
 \cosh\!\left(\frac{\log \mathcal{T}'(q) + (q - \mathcal{T}(q))}{2n}\right)^n -1
 + 
 \cosh\!\left(\frac{\log \mathcal{T}'(q)}2\right)-1\right]\dd\mu
 \leq\dPL.
\end{equation}
Using the elementary inequalities $\cosh v\geq1+v^2/2$ and $(1+x)^n\geq1+nx$ for $x\geq0$, coming for instance from the Taylor expansion of $\cosh$, we have that 
\begin{equation}\label{eq.log_derivative_control}
 \int\left[\log \mathcal{T}'(q) + (q - \mathcal{T}(q))\right]^2\dd\mu
 \leq8n\dPL,
\end{equation}
while $\int|\log\mathcal T'|^2\,\dd\mu\leq8\dPL$. More precisely,
the sum of these two squared quantities divided by $8n$ and $8$,
respectively, is at most $\dPL$. The elementary inequality
$|a-b|^2+|b|^2\leq2|a|^2+3|b|^2$ therefore proves
\eqref{eq.height_full_cost}.

To prove \eqref{eq.output_height_measure}, we start by going back to the pointwise estimate done in the proof of Theorem~\ref{thm.qualitative_PL}. Indeed, we recall the notation 
\[
    F_t \eqdef \{ f > t\}, \qquad
    G_t \eqdef \{ g > t\}, \qquad
    H_t \eqdef \{ h > t\},
\]
so that the Prékopa condition allows us to apply the Brunn-Minkowski inequality for any convex combination of level sets, in particular we have 
\begin{equation*}
    \left(\frac{|F_t|^{1/n}+|G_{T(t)}|^{1/n}}2\right)^n 
    \le 
    |H_{u(t)}|,
\end{equation*}
where $u(t) = t^{1/2}{T(t)}^{1/2}$. Following once again the arguments from Theorem~\ref{thm.qualitative_PL}, we multiply both sides by $u'(t)$, and we express the left-hand side in terms only of $|F_t|$ and the map $T$ by means of the Monge-Ampère equation associated with $T$. These computations then yield 
\[
  |H_{u(t)}|u'(t) \ge 
  |F_t| \left(1 + \dratio(t)\right) \left(1 + \did(t)\right). 
\] 
After the change of variations to the logarithmic scale from Section~\ref{sec.change_variables}, we retain the pointwise bound
\[
 Me^{-\mathcal{S}(q)}c(\mathcal{S}(q))\mathcal{S}'(q)
 \geq Me^{-q}a(q)
       [1+\dratio(Me^{-q})][1+\did(Me^{-q})]
 \geq Me^{-q}a(q),
\]
where $c(q) \eqdef |\{U < q\}| = |\{h > t(q)\}|$, and we set $\dd \mu_h(q) \eqdef M c(q) e^{-q} \dd q$. 

From the layer cake formula, we then have that 
\[
  \begin{aligned}
    1 + \dPL 
    &= \int_{\mathbb{R}^n} h(z) \dd z 
    = 
    \mu_h(\mathbb{R})
    = 
    M\int_{\mathbb{R}}
    e^{-\mathcal{S}(q)}c(\mathcal{S}(q))\mathcal{S}'(q) \dd q \\
    & \ge
    M\int_{\mathbb{R}}
    e^{-q}a(q) \dd q
    = 
    \int_{\mathbb{R}^n} f(x) \dd x = 1.  
  \end{aligned}
\]
So defining $\xi \eqdef \mu_h - {(\mathcal{S})}_\sharp \mu$, it follows directly that $\xi(\mathbb{R}) = \dPL$, and we must only show that it is a nonnegative measure. For this, given an arbitrary smooth nonnegative function $\varphi$ it suffices to integrate the above pointwise estimate against $\varphi(\mathcal{S})$ to get that 
\[
  \begin{aligned}
    \int_{\mathbb{R}} \varphi \dd \mu_h
    &= 
    M\int_{\mathbb{R}}
    \varphi(\mathcal{S}(q))
    e^{-\mathcal{S}(q)}c(\mathcal{S}(q))\mathcal{S}'(q) \dd q \\
    & \ge
    M\int_{\mathbb{R}}
    \varphi(\mathcal{S}(q))
    e^{-q}a(q) \dd q
    = 
    \int_{\mathbb{R}^n} \varphi \dd ({\mathcal{S}}_\sharp \mu). 
  \end{aligned}
\]
Hence $\xi \ge 0$ in the sense of distribution, and the result follows. 
\end{proof}

Now we move on to the vertical distance estimate between $\mu, \nu$ by means of the global deficit $\dPL$. Once again this will be a direct consequence of carefully keeping track of the remainder deficits coming from the qualitative proof. 

\begin{lemma}\label{lem.profiles}
For the measures $\mu,\nu$ defined in~\eqref{eq.log_transport_measures}, we have
\begin{equation}\label{eq.profile_comparison}
\norm{\mu-\nu}_{L^1(\mathbb R)}
=M\int_\R e^{-q}|a(q)-b(q)|\dd q
\leq4\sqrt{2n\dPL}.
\end{equation}
In particular, since $B_q=\varnothing$ for $q\leq0$,
\begin{equation}\label{eq.two_negative_heights}
 M\int_{-\infty}^0e^{-q}|A_q\sym(B_q-z)|\dd q
 =\mu(( -\infty,0])\leq4\sqrt{2n\dPL}
 \qquad\text{for every }z\in\R^n.
\end{equation}
\end{lemma}
\begin{proof}
  Carrying the deficits in Euclidean height in the qualitative proof yields the following inequality 
  \[
    \int_0^\infty 
    |F_t|
    \left[
      {\left(
        \frac{
          {T'(t)}^{\frac{1}{2n}} + {T'(t)}^{-\frac{1}{2n}}
        }{2}
      \right)}^n -1
    \right]\dd t
    \le \dPL. 
  \]
  Using $\cosh v\geq1+v^2/2$ and $(1+x)^n\geq1+nx$, this gives
  \begin{equation}\label{eq.deficit_estimate}
    \int_0^\infty |F_t|\,|\log T'(t)|^2\dd t
    \leq8n\dPL.
  \end{equation}

  Now, to compare the two distributions, it is useful to recall the Monge-Ampère equation matching them 
  \[
    |F_t| = |G_{T(t)}|T'(t)
  \]
  and define the following change of variables 
  \[
    S(t) = \int_0^t \min\{1, T'(r)\} \dd r. 
  \]
  As a result, we clearly have $S(t) \le t$ and $S(t) \le T(t)$. We then define the profiles 
  \[
    \alpha(s) \eqdef |F_{S^{-1}(s)}|, \qquad 
    \beta(s)  \eqdef |G_{T\circ S^{-1}(s)}|
  \]
  meaning that $\alpha(S(t)) = |F_t|$ and $\beta(S(t)) = |G_{T(t)}|$. Since $t \mapsto |F_t|$ is nonincreasing, by definition, as well as $t \mapsto |G_{T(t)}|$ thanks to the monotonicity of $T$, it follows that 
  \[
    \alpha(s) \le |F_s|, \qquad \beta(s) \le |G_s|
  \]
  for all $s \ge 0$. 

  Now, turning our attention to our quantity of interest, we have from our previous change of variables and the triangle inequality that
  \[
    \norm{\mu - \nu}_{L^1}
    = 
    \int_0^\infty
    ||F_t| - |G_t|| \dd t 
    \le
    \underbrace{
      \int_0^\infty
    ||F_t| - \alpha(t)| \dd t
    }_{I_1}
    + 
    \underbrace{
      \int_0^\infty
    ||G_t| - \beta(t)| \dd t
    }_{I_2}
    + 
    \underbrace{
          \int_0^\infty
    |\alpha(t) - \beta(t)| \dd t
    }_{I_3},
  \]
  while each of these integrals is given by 
  \[
  \begin{aligned}
    I_1 &= \int_0^\infty 
    |F_t|\left(
      1 - \min\{1,T'\}  
    \right)\dd t, 
    \quad
    I_2 = \int_0^\infty 
    |F_t|\left(
      1 - \frac{\min\{1,T'\}}{T'} 
    \right)\dd t, \\
    I_3 &= \int_0^\infty 
    |F_t|
    \min\{1,T'\}
    \left|
      1 - \frac{1}{T'} 
    \right|\dd t.
  \end{aligned}
  \]
  Summing all the integrands above we obtain
  \[
    2\frac{|T'-1|}{\max\{1,T'\}}
    \leq2|\log T'|.
  \]
  Combining this with~\eqref{eq.deficit_estimate} and Cauchy--Schwarz
  gives~\eqref{eq.profile_comparison}. The last assertion follows
  from $b(q)=0$ for $q\leq0$.
\end{proof}

\subsection{Overlapping epigraph chunks, their volume, and the lift}\label{sec.epigraph_chunks}

The biggest virtue of the qualitative proof of the Prékopa-Leindler inequality from {\color{black}Theorem~\ref{thm.qualitative_PL}} is the confirmation of our understanding that it can be seen as a functional version of Brunn-Minkowski. Therefore, it is conceptually clear that the stability of Prékopa-Leindler could, in principle, be attained by applying the stability of Brunn-Minkowski over its level sets and gluing everything back together via the Cavalieri principle to obtain an estimate on the $L^1$ distance between the log-concave measures $f,g$. The issue is that from Brunn-Minkowski, we obtain a continuum of translations, potentially a different one for each level set, which cannot be easily glue back together. 

A very natural way of circumventing this is to make quantitative the philosophy behind Lemma~\ref{lemma.fixed_translation}, where the translation $z_t$ of the level set $t$ is proved to behave in the order $|t-s|^2$, being therefore constant. Instead, our approach is to consider \textit{chunks of level-sets from the epigraphs} of the functions $f,g$. 

For this, we shall work with the epigraphs in the log-scale, where the Prékopa condition becomes a Minkowski sum, and consider chunks of level sets so that we end up with countably many translations, so having a common one for each chunk, instead of a continuum. To implement this idea,
consider a log-concave function $e^{-P}$ associated with a convex potential $P : \mathbb{R}^n \to \mathbb{R}\cup \{+\infty\}$. One can define a chunk of the epigraph of $e^{-P}$ around the level $r$ as 
\[
  \mathcal{Q}_r(P)
  \eqdef 
  \left\{
    (x,q) \in \mathbb{R}^n \times \mathbb{R}: 
    P(x) < q, \quad |q - r| < 2
  \right\}.
\]
Since $P$ is convex and $e^{-P}$ is integrable, its sublevels are bounded, and therefore $\mathcal Q_r(P)$ is a bounded convex set of $\mathbb{R}^{n+1}$. 

A simple application of Fubini's Theorem shows that the volume of $\mathcal{Q}_r(P)$ is written as
\[
  |\mathcal{Q}_r(P)| 
  = 
  \int_{\{
    P(x) < q, \ |q-r| < 2
  \}}
  1 \dd x \dd q
  = 
  \int_{r - 2}^{r + 2} 
  |\{P(x) < q\}| \dd q.
\]
Hence, recalling $\dd \mu = M e^{-q} a(q) \dd q$, where $a(q) = |\{V < q\}|$, we can express the weighted volume as a convolution with a squared window 
\begin{equation}\label{eq.square_window_convolution}
  M e^{-r}
  |\mathcal{Q}_r(V)| 
  = 
  M e^{-r}
  \int_{r - 2}^{r + 2} 
  a(q) \dd q
  = 
  \int_{-\infty}^{+\infty}
  e^{-(r-q)} \mathbf{1}_{|r-q| < 2} \dd \mu 
  = 
  \mu \star \kappa_0(r),
\end{equation}
where the square window is defined as $\kappa_0 \eqdef e^{-s} \mathbf{1}_{|s| < 2}$. The analogous identity holds for $\nu$. 

In the sequel, we shall consider the following chunks: for $r\in\R$ define
\begin{equation}\label{eq.chunks_def}
    E_r \eqdef \mathcal{Q}_r(V), \quad 
    F_r \eqdef \mathcal{Q}_r(W), \quad 
    H_r \eqdef \mathcal{Q}_r(U),
\end{equation}
so that we have 
\begin{equation}\label{eq.chunks_volume_identity}
  Me^{-r}|E_r| = \mu \star \kappa_0(r), \quad 
  Me^{-r}|F_r| = \nu \star \kappa_0(r). 
\end{equation}
As a result, from Fubini and fundamental properties of the convolution operation (see for instance~\cite[Chap.~2]{ambrosio2000functions}), setting $K_0 \eqdef \norm{\kappa_0}_{L^1(\mathbb{R})}$ it holds that 
\begin{equation}
  \begin{aligned}
    M \int_{\mathbb{R}} e^{-r}|E_r| \dd r 
    &= K_0 \norm{f}_{L^1(\mathbb{R}^n)} = K_0, \\
    M \int_{\mathbb{R}} e^{-r}|F_r| \dd r 
    &= K_0 \norm{g}_{L^1(\mathbb{R}^n)} = K_0, \\
    M \int_{\mathbb{R}} e^{-r}|H_r| \dd r 
    &= K_0 \norm{h}_{L^1(\mathbb{R}^n)} = K_0(1 + \dPL).
  \end{aligned}
\end{equation}
Notice that the integrals must be taken over all the real line, which will be slight technical complication afterward. 

This simple computation shows that the integrals of these quantities carry exactly the information of $\dPL$. As a result, using this construction to estimate the desired distance with the deficit depends on being able to control the terms $\dBM(E_r, F_r)$ $\dvol(|E_r|, |F_r|)$, which we do in the following Lemma. 

\begin{lemma}\label{lemma.integral_deficits_control}
  For $r > 0$, define the quantity 
  \begin{equation}
    \delta(r) \eqdef 
    |H_r| - 
    {
      \left(
        \frac{|E_r|^{1/d} + |F_r|^{1/d}}{2}
      \right)
    }^{d}
    \text{ with } d = n + 1. 
  \end{equation}
  Then for $r > 1/4$ we have that both level set chunks $E_r, F_r$ have positive volume, so that $0 \le \dBM(E_r, F_r) \le \delta(r)$, and there exists an absolute constant $C$ for which it holds that 
  \begin{align}
    \label{eq.estimate_dvol_exp}
    M\int_{1/2}^\infty
      e^{-r} \dvol(E_r, F_r) \dd r & \le C(n+1)^2 \dPL, \\
    \label{eq.estimate_dBM_exp}
    M\int_{1/2}^\infty
      e^{-r} \dBM(E_r, F_r) \dd r & \le C(n+1)^2 \dPL.
  \end{align}
\end{lemma}
\begin{proof}
  Let us start with estimate~\eqref{eq.estimate_dvol_exp}, which will later be used in the proof of~\eqref{eq.estimate_dBM_exp}. From Lemma~\ref{lemma.scalar_Hellinger}, the chunk volume identity from~\eqref{eq.chunks_volume_identity}, and the fact that $\mathcal{T}_\sharp \mu = \nu$ is optimal alongside the transport estimates from $\mu$ to $\nu$ via the deficit $\dPL$ obtained in Lemma~\ref{lemma.height_transport}, it suffices to show that 
  \begin{equation}\label{eq.hellinger_kantorovitch_bound}
    \begin{aligned}
    \int_{1/4}^\infty 
    {
      \left(
        \sqrt{\mu \star \kappa_0}(r) - \sqrt{\nu \star \kappa_0}(r)
      \right)
    }^2 \dd r
    &\le C(n+1) 
    \int_{\mathbb{R}} |q - \mathcal{T}(q)|^2 \dd \mu(q)
    \\
    &=C(n+1)W_2^2(\mu,\nu)
    \le C(n+1)^2\dPL.
    \end{aligned}
  \end{equation}

 Let then ${(\varrho_t)}_{t \in [0,1]}$ be the Wasserstein geodesic connecting $\mu, \nu$, which are defined as 
  \[
    \varrho_t \eqdef (\mathcal{T}_t)_{\sharp} \mu, \quad 
    \mathcal{T}_t \eqdef 
    (1-t)\id + t \mathcal{T},
  \]
  as considered in Lemma~\ref{lemma.tail_ab}. The Benamou-Brenier formula states that this interpolation solves the continuity equation 
  \[
    \partial_t \varrho_t + \divv (\varrho_t v_t) = 0, 
    \text{  where  } 
    v_t \eqdef (\mathcal{T} - \id) \circ \mathcal{T}_t^{-1},
  \]
  or equivalently, if we work with momentum variables $M_t \eqdef \varrho_t v_t$ so that $\partial_t \varrho_t + \divv M_t = 0$ and 
  \[
      \int_{0}^1 \int_\mathbb{R} 
      \frac{|M_t|^2}{\varrho_t} \dd q \dd t 
      = 
      W_2^2(\mu, \nu) \le16(n+1)\dPL. 
  \] 

  Now we want to turn this into a suitable interpolation for the convolutions $\mu \star \kappa_0, \nu \star \kappa_0$. This can easily be done as 
  \[
    \eta_t(r) \eqdef \varrho_t \star \kappa_0(r) 
    = 
    \int_{r-2}^{r+2} e^{q-r} \varrho_t(q) \dd q, \quad \text{ for all } t \in [0,1]. 
  \]
  Going back to the very first quantity we wanted to estimate in this proof, the fundamental theorem of calculus gives 
  \begin{align*}
    \int_{1/4}^\infty 
    {
      \left(
        \sqrt{\mu \star \kappa_0}(r) - \sqrt{\nu \star \kappa_0}(r)
      \right)
    }^2 \dd r 
    \le 
    \int_{1/4}^\infty 
    {\left(
      \int_0^1 \frac{\partial_t \eta_t}{\sqrt{\eta_t}} \dd t
    \right)
    }^2 \dd r
    \le 
    \int_{1/4}^\infty \int_0^1
     \frac{|\partial_t \eta_t|^2}{\eta_t} \dd t \dd r. 
  \end{align*}

  In the sequel, we can compute $\partial_t \eta_t$ through the continuity equation solved by $\varrho_t$ and elementary properties of the convolution operator. Indeed, since the Wasserstein geodesic $\varrho_t$ is an absolutely continuous curve we can apply differentiation under the integral sign to obtain that 
  \begin{equation*}
      \begin{aligned}
          \partial_t \eta_t(r) 
          &= 
          \int_{\mathbb{R}} 
          \partial_t \varrho_t (q) \kappa_0(r - q)\dd q
          = 
          -\int_{r - 2}^{r + 2} 
            \partial_q M_t(q)
          e^{-(r-q)} \dd q \\
          &= 
          e^{-2} M_t(r-2) 
          - 
          e^2 M_t(r + 2) 
          + 
          \int_{r-2}^{r+2} e^{q-r} M_t(q)\dd q. 
      \end{aligned}
  \end{equation*}

  As a result, back to the previous estimate on the Hellinger distance, our goal is to estimate 
  \[
    \frac{|\partial_t \eta_t|^2}{\eta_t} 
    \le 
    C\left[
      \frac{M_t(r-2)^2}{\eta_t} + \frac{M_t(r+2)^2}{\eta_t} 
      + 
      \frac{1}{\eta_t} 
      {
        \left(
          \int_{r-2}^{r+2} 
          e^{q-r} M_t(q) \dd q  
        \right)
      }^{2}
    \right],
  \]
  and it suffices to estimate each term. The first one uses~\eqref{eq.lemma_interp_1st_identity} the first identity from item (c) of Lemma~\ref{lemma.tail_ab}. Indeed, using this mentioned identity for each $q$ in the integration interval being bounded from below with the value at $r-2$ we get 
  \[
    \eta_t(r) = 
    \int_{r-2}^{r+2} e^{q-r} \varrho_t(q) \dd q 
    \ge 
    \int_{r-2}^{r+2} e^{-2}\varrho_t(r-2) \dd q 
    = 
    4e^{-2} \varrho_t(r-2). 
  \]
  As result, using that $M_t = \varrho_t v_t$, we get 
  \[
    \frac{M_t(r-2)^2}{\eta_t} 
    \le 
    C\varrho_t(r-2)|v_t(r-2)|^2. 
  \]
  Similarly for the second term, write $Q=r+2$ and
  $\varepsilon=1/(n+1)$. Since $m_t\leq0$ and $r>1/4$, we have
  $Q-m_t>2$. For $0\leq s\leq\varepsilon$, the second bound
  \eqref{eq.lemma_interp_2nd_identity} gives
  \[
    \varrho_t(Q-s)\geq
    \left(1-\frac{1}{2(n+1)}\right)^{n+1}\varrho_t(Q)
    \geq\frac12\varrho_t(Q).
  \]
  Integrating only on $[Q-\varepsilon,Q]$, we obtain
  $\eta_t(r)\geq c\varrho_t(r+2)/(n+1)$, and hence
  \[
    \frac{M_t(r+2)^2}{\eta_t(r)}
    \leq C(n+1)\varrho_t(r+2)|v_t(r+2)|^2.
  \]
  As for the last term, we use Cauchy--Schwarz:
  \[
    \frac{1}{\eta_t(r)}
    \left(\int_{r-2}^{r+2}e^{q-r}M_t(q)\dd q\right)^2
    \leq\int_{r-2}^{r+2}e^{q-r}\varrho_t(q)|v_t(q)|^2\dd q.
  \]
  Its integral in $r$ is at most $K_0\int\varrho_t|v_t|^2$ by Fubini.
  Integrating all these estimates over $r > 1/4$ we obtain the desired bound~\eqref{eq.hellinger_kantorovitch_bound}. Finally, using the upper bound on~\eqref{eq.dvol_bounds} from Lemma~\ref{lemma.scalar_Hellinger}, the first identity~\eqref{eq.estimate_dvol_exp} follows.

Now, for the second identity~\eqref{eq.estimate_dBM_exp}, set $K_0=\int_\R\kappa_0$. Fubini gives the exact global mass identities
\begin{equation}\label{eq.window_masses}
 M\int_\R e^{-r}|E_r|\dd r=M\int_\R e^{-r}|F_r|\dd r=K_0,
 \qquad M\int_\R e^{-r}|H_r|\dd r=K_0(1+\dPL),
\end{equation}
Notice that these identities require an integral over $r\in\R$, since the overlap of level set chunks will forcibly give a positive mass to some negative centers. Hence instead of the previous integrals of fixed volume, we use a selecting kernel 
$\beta\in C^2(\R)$ with $0\leq\beta\leq1$, $\beta=0$ on
$(-\infty,1/4]$, and $\beta=1$ on $[1/2,\infty)$, with bounded
first and second derivatives. This induces the following test function
\[
 \Psi(q)=\int_{-2}^2 e^{-t}\beta(q+t)\dd t,
\]
which satisfies $0\leq\Psi\leq K_0$ and
$\|\Psi''\|_\infty\leq K_0\|\beta''\|_\infty$. Now applying Fubini to the weighted mass integral gives
\[
 \begin{aligned}
 M\int_\R\beta(r)e^{-r}|E_r|\dd r
 &=\int_\R\Psi(q)\dd\mu(q),\\
 M\int_\R\beta(r)e^{-r}|F_r|\dd r
 &=\int_\R\Psi(q)\dd\nu(q),\\
 M\int_\R\beta(r)e^{-r}|H_r|\dd r
 &=\int_\R\Psi(q)\dd\mu_h(q).
 \end{aligned}
\]
We now use \eqref{eq.output_height_measure}; namely
$\mu_h=\mathcal S_\sharp\mu+\xi$ with
$\mathcal S=(\id+\mathcal T)/2$, together with
$\nu=\mathcal T_\sharp\mu$.  Subtracting the first two identities from the
third gives
\begin{equation}\label{eq.signed_balance}
 \begin{aligned}
 &M\int_\R\beta(r)e^{-r}
       \left[|H_r|-\frac{|E_r|+|F_r|}{2}\right]\dd r\\
 &\quad=\int_{\mathbb{R}}\Psi\,\dd\left(\mu_h-\frac{\mu+\nu}{2}\right)\\
 &\quad=\int_{\mathbb{R}}\Psi\,\dd\xi+
 \int_{\mathbb{R}}\left[\Psi\!\left(\frac{q+\mathcal T(q)}2\right)
            -\frac{\Psi(q)+\Psi(\mathcal T(q))}{2}\right]\dd\mu(q)\\
 &\quad\leq K_0\dPL+
       \frac{\|\Psi''\|_\infty}{8}
        \int_{\mathbb{R}}|\mathcal T(q)-q|^2\dd\mu(q)
 \leq C(n+1)\dPL, 
 \end{aligned}
\end{equation}
where we have used that the nonnegative measure $\xi$ has mass
$\dPL$ from Lemma~\ref{lemma.height_transport}.

By definition,
\[
 \delta(r)=|H_r|-\frac{|E_r|+|F_r|}{2}+\dvol(E_r,F_r).
\]
The scalar estimate already proved on $r\geq1/4$, together with
\eqref{eq.signed_balance}, therefore yields
\[
 M\int_\R\beta(r)e^{-r}\delta(r)\dd r\leq C(n+1)^2\dPL.
\]
On the other hand, if $r\geq1/4$ both $E_r$ and $F_r$ have positive volume so that $0 \le \dBM(E_r, F_r) \le \delta(r)$ and we get that
\[
 M\int_{1/2}^\infty e^{-r}\dBM(E_r,F_r)\dd r
 \leq M\int_{1/2}^\infty e^{-r}\delta(r)\dd r
 \leq C(n+1)^2\dPL,
\]
and assertion~\eqref{eq.estimate_dBM_exp} follows. 
\end{proof}

\subsection{Proof of the quantitative Pr\'ekopa--Leindler estimate, Theorem~\ref{thm.quantitative_PL}}
\label{sec.control_translations}

We now combine the estimates obtained above and prove
\eqref{eq.quantitative_PL}. Recall that we are first working at the
midpoint, under the normalization \eqref{eq.quant_normalization}, and that $0<\dPL\leq1$. Our choice of the common reference height
ensures that both $E_r$ and $F_r$ have positive volume for every
$r\geq1/2$, without a smallness condition on $m$.

We start by applying Theorem~\ref{thm.BM_stability} to each pair of
chunks $E_r,F_r\subset\R^{n+1}$. Combining the conclusion of Lemma~\ref{lemma.scalar_Hellinger}, identity \eqref{eq.window_masses}, and an easy application of the Cauchy--Schwarz inequality give
\begingroup\small
\begin{equation}\label{eq.all_window_translations}
  \begin{aligned}
    &M\int_{1/2}^\infty e^{-r}
       \inf_{v\in\R^{n+1}}|E_r\sym(F_r-v)|\dd r \\
    &\quad\leq\sqrt{B_{n+1}^{\rm conv}}
    \left(M\int_{1/2}^\infty e^{-r}(|E_r|+|F_r|)\dd r\right)^{1/2}
    \left(M\int_{1/2}^\infty e^{-r}
          [\dBM(E_r,F_r)+\dvol(E_r,F_r)]\dd r\right)^{1/2}\\
    &\quad\leq C(n+1)\sqrt{B_{n+1}^{\rm conv}\dPL}.
  \end{aligned}
\end{equation}
\endgroup
For each $r$, the translation furnished by Brunn--Minkowski belongs to
$\R^{n+1}$ and may depend on the chunk's center. We write it as
$v=(z,u)$, where $z\in\R^n$ is the horizontal component and $u\in\R$
is the vertical, or height, component. Only the $z$ component corresponds to a translation of the original functions. We must therefore remove $u$ and
then show that the remaining horizontal translations can all be replaced
by a single horizontal translation.

\vspace{2mm}

\emph{Step 1: Removing the vertical translations.}
For $r\geq1/2$, we recall that both $E_r, F_r$ are non-empty, and we set
\[
 \chi(r)=\int_{r-2}^{r+2}|a(q)-b(q)|\dd q.
\]
We first show that, for every $v=(z,u)$,
\begin{equation}\label{eq.pinning}
 |E_r\sym(F_r-(z,0))|
 \leq2|E_r\sym(F_r-v)|+\chi(r).
\end{equation}
To prove this, consider the height marginals
\[
 m_E(q)=a(q)\mathbf1_{(r-2,r+2)}(q),\qquad
 m_F(q)=b(q)\mathbf1_{(r-2,r+2)}(q).
\]
The sections of $F_r$ and $F_r-(0,u)$ at each height are nested
level sets, or one of them is empty. Hence
\[
 |F_r\sym(F_r-(0,u))|=\|m_F-m_F(\cdot+u)\|_1.
\]
Taking marginals decreases the $L^1$ distance, and therefore
\[
 \|m_E-m_F(\cdot+u)\|_1\leq|E_r\sym(F_r-v)|,
 \qquad\|m_E-m_F\|_1=\chi(r).
\]
Two applications of the triangle inequality prove, more precisely,
\[
 |E_r\sym(F_r-(z,0))|
 \leq2|E_r\sym(F_r-v)|+\chi(r).
\]

Define the error after allowing only horizontal translations of chunks by
\begin{equation}\label{eq.local_error}
 \mathcal E(r)=\inf_{z\in\R^n}
       \int_{r-2}^{r+2}|A_q\sym(B_q-z)|\dd q.
\end{equation}
Now, once height translations are removed, epigraph translations are reflected back to the original functions via the identity 
\[
  |E_r\sym(F_r-(z,0))|
 =\int_{r-2}^{r+2}|A_q\sym(B_q-z)|\dd q,
\]
which can be combined with~\eqref{eq.pinning} yielding
\[
 \mathcal E(r)\leq2\left(
   \inf_{v\in\R^{n+1}}|E_r\sym(F_r-v)|+\chi(r)\right).
\]
The first term is controlled after integration by
\eqref{eq.all_window_translations}. For the second one, Fubini gives
\begin{align*}
 M\int_{1/2}^\infty e^{-r}\chi(r)\dd r
 &\leq K_0M\int_\R e^{-q}|a(q)-b(q)|\dd q\\
 &=K_0\|\mu-\nu\|_{L^1(\R)}
 \leq C\sqrt{n\dPL},
\end{align*}
where the last estimate is \eqref{eq.profile_comparison}. Consequently,
\begin{equation}\label{eq.local_integrated}
 M\int_{1/2}^\infty e^{-r}\mathcal E(r)\dd r
 \leq C(n+1)(1+\sqrt{B_{n+1}^{\rm conv}})\sqrt{\dPL}.
\end{equation}

\vspace{2mm}

\emph{Step 2: Joining the horizontal translations.}
Estimate \eqref{eq.local_integrated} still allows a different horizontal
translation for each chunk, for a continuum of these, which yields a continuum of translations. So first we bound this integral from below with a series, which selects only countably many chunks, and we exploit the overlap of consecutive ones to replace their corresponding sequence of horizontal translations by a single vector.

For each integer $k\geq0$, first notice that for $r\in [k+1/2,k+3/4]$ we have $e^{-k}\le e^{3/4} e^{-r}$, hence applying the mean value theorem we can choose $r_k$ in this interval so that
\begin{equation}\label{eq.sample_centers}
 \sum_{k\geq0}e^{-k}\mathcal E(r_k)
 \leq4e^{3/4}\int_{1/2}^\infty e^{-r}\mathcal E(r)\dd r.
\end{equation}

Let $z_k$ be a minimizing translation for $\mathcal E(r_k)$. Both
intervals $[r_k-2,r_k+2]$ and $[r_{k+1}-2,r_{k+1}+2]$ contain
$[k+1,k+2]$. Hence, up to another averaging argument, there exists
$\rho_k\in[k+1,k+2]$ such that
\begin{equation}\label{eq.overlap_increment}
 |B_{\rho_k}\sym(B_{\rho_k}+z_{k+1}-z_k)|
 \leq\mathcal E(r_k)+\mathcal E(r_{k+1}).
\end{equation}

Estimate~\eqref{eq.overlap_increment} controls the effect of the
translation $z_{k+1}-z_k$ at the single level $\rho_k$. To join all
translations, we need to control the same effect at every other level,
rather than the Euclidean distance $|z_{k+1}-z_k|$. The convexity of
$W$ and $\min W=0$ give, for $0<r\leq s$ and $v\in\R^n$,
\begin{equation}\label{eq.propagation}
 |B_r\sym(B_r+v)|
 \leq|B_s\sym(B_s+v)|
 \leq(s/r)^n|B_r\sym(B_r+v)|.
\end{equation}
Indeed, integration along lines parallel to $v$ writes the symmetric
difference of a convex body $K$ as twice the integral of
$\min\{|v|,\text{length of the chord of }K\}$. This quantity is
monotone under inclusion, and increases by at most $t^n$ under a
dilation by $t\geq1$. Apply these facts to the inclusions from
Lemma~\ref{lemma.tail_ab}(a).

Now we fix the translation $z_n$, associated with a chunk near height
$n$. This choice keeps track of the dimensional loss in the following
summation. If $q\in[k,k+1]$, the triangle inequality gives
\begin{equation}\label{eq.global_from_local}
 |A_q\sym(B_q-z_n)|
 \leq |A_q\sym(B_q-z_k)|
 +\sum_{j=\min\{k,n\}}^{\max\{k,n\}-1}
     |B_q\sym(B_q+z_{j+1}-z_j)|.
\end{equation}
The local terms integrate to at most
$\sum_{k\geq0}e^{-k}\mathcal E(r_k)$, because
$[k,k+1]\subset(r_k-2,r_k+2)$.

For $j<n$, the increment is used only at heights $q\leq j+1\leq\rho_j$.
The first inequality in~\eqref{eq.propagation} therefore bounds all
these contributions by
\[
 \sum_{j=0}^{n-1}[\mathcal E(r_j)+\mathcal E(r_{j+1})]
 \leq2e^n\sum_{k\geq0}e^{-k}\mathcal E(r_k).
\]
For $j\geq n$, the increment is used only at $q\geq j+1$. Combining
both inequalities in~\eqref{eq.propagation} and $\rho_j\geq j+1$
gives
\begin{equation}\label{eq.propagated_increment}
 |B_q\sym(B_q+z_{j+1}-z_j)|
 \leq\left(\frac q{j+1}\right)^n
       [\mathcal E(r_j)+\mathcal E(r_{j+1})],\qquad q\geq j+1.
\end{equation}
The corresponding tail satisfies
\[
 \begin{aligned}
 \int_{j+1}^{\infty}e^{-q}\left(\frac q{j+1}\right)^n\dd q
 &\leq e^{-j-1}\int_0^\infty
       e^{-s(1-n/(j+1))}\dd s\\
 &=e^{-j-1}\frac{j+1}{j+1-n}
 \leq(n+1)e^{-j-1}.
 \end{aligned}
\]
Reversing the order of summation in~\eqref{eq.global_from_local},
and then using~\eqref{eq.sample_centers}, proves
\begin{equation}\label{eq.gluing}
 \int_0^\infty e^{-q}|A_q\sym(B_q-z_n)|\dd q
 \leq Ce^n\int_{1/2}^\infty e^{-r}\mathcal E(r)\dd r.
\end{equation}
Notice that the argument controls the symmetric difference produced by
the changes $z_{j+1}-z_j$; no estimate of their Euclidean lengths is
required. Anchoring at $z_n$ avoids propagating the first increment
from a height of order one up to heights of order $n$.

Finally, multiply~\eqref{eq.gluing} by $M$ and use
\eqref{eq.local_integrated}. The negative heights are controlled by
\eqref{eq.two_negative_heights}. By~\eqref{eq.log_layer_cake}, we obtain
\begin{equation}\label{eq.midpoint_conclusion}
 \|f-g(\cdot+z_n)\|_1
 \leq Ce^n(n+1)(1+\sqrt{B_{n+1}^{\rm conv}})\sqrt{\dPL}
 \leq e^{Cn}\sqrt{\dPL}.
\end{equation}
This proves the midpoint estimate for $0<\dPL\leq1$. The reduction in
Corollary~\ref{corollary.quant_mid_point_implies_full}, followed by the
mass normalization, proves \eqref{eq.quantitative_PL} for every
$\lambda\in(0,1)$.

\section{Once again, but faster: one log-concave input and beyond}\label{sec.quantitativePL_1logconcave}

A careful inspection of the previous proof reveals that the assumption of $f,g$ being log-concave is only actively essential on the growth estimates proved in Lemma~\ref{lemma.tail_ab}, which in their turn are crucial for estimating the Hellinger-type distance 
\[
    \norm{
        \sqrt{\mu \star \kappa_0}
        - 
        \sqrt{\nu \star \kappa_0}
    }_{L^2}^2
    \le 
    C_n W_2^2(\mu, \nu),
\]
with the Wasserstein distance between the height distributions. If the kernel $\kappa_0$ was smoother instead of the squared window kernel 
$\kappa_0 \eqdef e^{ - s}\mathbf{1}_{\{|s| < 2\}}$,
this estimate could have been much more easily established. Of course, another kernel would mean another way of selecting chunks that is much less natural, but let us entertain this idea for a while. 

For this reason, in this section we consider the following kernel generated by the trapezoidal window $\phi$ below
\begin{equation}\label{eq.smooth_kernel}
    \kappa (s) \eqdef \pi e^{-s} \phi(s)^2, 
    \text{ where }
    \phi(s)\eqdef
    \begin{cases}
      1, & |s| \le 1,\\
      2-|s|, & 1\leq |s|\leq 2,\\
      0, & \text{otherwise}.
    \end{cases}
\end{equation}
The same argument would work with other kernels, but this one has the virtue of being conceptually close to the square window kernel $\kappa_0$. The middle part being constant and equal $1$, has the same function from before regarding the overlap of chunks, while the linear decay of $\phi$ and squaring it in the definition of $\kappa$ ensures that $\sqrt{\kappa} \in H^1(\mathbb{R})$.

Using this improved regularity, we can establish the Wasserstein estimate with a very simple argument. The issue is that if $f$ is not log-concave, then its corresponding height distribution $\mu$ might not have finite second moments in the logarithmic scale, so the Wasserstein distance between them might not be finite. 

\begin{lemma}\label{lemma.smooth_height}
Let $\mu,\nu$ be the logarithmic height distributions of two normalized
inputs, defined as in~\eqref{eq.log_transport_measures}, and let
$\mathcal T$ be their monotone transport. No log-concavity is assumed.
Then we have
\begin{equation}\label{eq.smooth_height_truncated_cost}
 \int_\R\bigl|\sqrt{\mu\star\kappa(r)}-\sqrt{\nu\star\kappa(r)}\bigr|^2\dd r
 \leq 
 C_\kappa\int_\R \min\{|q-\mathcal{T}(q)|^2,1\}\dd\mu(q).
\end{equation}
Here $C_\kappa=\max\{\|(\sqrt\kappa)'\|_2^2,4\|\kappa\|_1\}$, which is absolute for the fixed kernel above.
If in addition $\mu,\nu\in\mathscr P_2(\R)$, then
\begin{equation}\label{eq.smooth_height_W2}
 \int_\R\bigl|\sqrt{\mu\star\kappa(r)}-\sqrt{\nu\star\kappa(r)}\bigr|^2\dd r
 \leq \|\bigl(\sqrt{\kappa}\bigr)'\|_2^2 W_2^2(\mu,\nu).
\end{equation}
\end{lemma}
\begin{proof}
Recall that $\mathcal{T}$ is the monotone map transporting $\mu$ onto $\nu$;
its definition as a monotone rearrangement does not require finite second
moments. Define the plan $\pi_0 = {(\id, \mathcal{T})}_\sharp \mu$, whose
first and second marginals are $\mu$ and $\nu$. Then
\[
 \mu\star\kappa(r)=\int \kappa(r-q)\dd\pi_0(q,s),
 \qquad
 \nu\star\kappa(r)=\int \kappa(r-s)\dd\pi_0(q,s).
\]
Since $\kappa\geq0$ is bounded and $\pi_0$ is a probability measure, for every fixed $r\in\R$ the functions
\[
 \varphi_r(q,s)\eqdef\sqrt{\kappa(r-q)},\qquad
 \psi_r(q,s)\eqdef\sqrt{\kappa(r-s)}
\]
belong to $L^2(\pi_0)$. As the first and second marginals of $\pi_0$ are $\mu$ and $\nu$, the two identities above read
\[
 \sqrt{\mu\star\kappa(r)}=\|\varphi_r\|_{L^2(\pi_0)},
 \qquad
 \sqrt{\nu\star\kappa(r)}=\|\psi_r\|_{L^2(\pi_0)}.
\]
The reverse triangle inequality for the norm of $L^2(\pi_0)$,
$$\bigl|\|\varphi_r\|_{L^2(\pi_0)}-\|\psi_r\|_{L^2(\pi_0)}\bigr|\leq\|\varphi_r-\psi_r\|_{L^2(\pi_0)},$$ 
therefore gives, after squaring, for every $r$,
\[
 \bigl|\sqrt{\mu\star\kappa(r)}-\sqrt{\nu\star\kappa(r)}\bigr|^2
 \leq\int\bigl|\sqrt{\kappa(r-q)}
                 -\sqrt{\kappa(r-s)}\bigr|^2\dd\pi_0(q,s).
\]
Integrating in $r$ and using Fubini yields
\begin{equation}\label{eq.smooth_height_coupling}
 \int_\R\bigl|\sqrt{\mu\star\kappa}-\sqrt{\nu\star\kappa}\bigr|^2 \dd r
 \leq\int_{\R^2}
       \|\sqrt{\kappa(\cdot-q)}
          -\sqrt{\kappa(\cdot-s)}\|_2^2\dd\pi_0(q,s).
\end{equation}
For the chosen kernel it holds that $\sqrt{\kappa} \in H^1(\R)$, so that the fundamental theorem of calculus for translations in $H^1$ gives
\[
 \|\sqrt{\kappa(\cdot-q)}-\sqrt{\kappa(\cdot-s)}\|_2^2
 \leq\|\bigl(\sqrt{\kappa}\bigr)'\|_2^2|q-s|^2.
\]
The triangle inequality also gives the uniform bound
\[
 \|\sqrt{\kappa(\cdot-q)}-\sqrt{\kappa(\cdot-s)}\|_2^2
 \leq4\|\sqrt{\kappa}\|_2^2.
\]
Taking the better of these two estimates in
\eqref{eq.smooth_height_coupling} proves
\eqref{eq.smooth_height_truncated_cost}, for instance with
\[
 C_\kappa=\max\left\{
   \|\bigl(\sqrt{\kappa}\bigr)'\|_2^2,
   4\|\sqrt{\kappa}\|_2^2
 \right\}.
\]
If $\mu,\nu\in\mathscr P_2$, use
only the first translation estimate and the optimality of $\mathcal{T}$ gives~\eqref{eq.smooth_height_W2}.
\end{proof}

On the one hand, this argument has the clear advantage of not using log-concavity, but the interpretation of the convolution kernel as the volume of a chunk of the epigraphs is lost, instead we obtain a weighted volume of epigraph chunks. The issue is that to the best of our knowledge, there is no quantitative Brunn-Minkowski theory that can deal with the weighted measure that appears naturally. As a result, we unfold the epigraph in $\mathbb{R}^{n+3}$ which implies a disintegration of the weighted measure over $\mathbb{R}^{n+1}$ into the Lebesgue measure over $\mathbb{R}^{n+3}$. After this lift, the remaining arguments remain almost unchanged. 

\subsection{Review of notation}\label{sec.1logconcave_review_notation}

Throughout this section, $f$ is an arbitrary nonnegative measurable function and $g$ is log-concave. Here the reduction to the midpoint argument from Corollary~\ref{corollary.quant_mid_point_implies_full} does not hold, since we need both inputs to be log-concave. As a result, we keep the dependence on $\lambda$. We can still assume the normalization $\int f=\int g=1$, but the Prékopa condition~\eqref{eq.prekopa_condition} for $h$ now reads 
\[
    h(\lambda x + (1-\lambda) y) 
    \ge 
    f(x)^{\lambda} g(y)^{1-\lambda}, 
\]
so we set 
\[
 \tau\eqdef\min\{\lambda,1-\lambda\},\qquad
 \dPL=\dPL(\lambda)\eqdef\int_{\R^n}h-1.
\]
Our goal is to obtain an estimate of order $\sqrt{\dPL/\tau}$, while keeping track of the dimensional constant. We assume $\dPL>0$; no dimensional smallness threshold will be needed.

We keep the notation from Section~\ref{sec.change_variables}, taking the smaller of the two maximum heights as the common reference:
\begin{equation}\label{eq.1logconcave_height}
 \begin{gathered}
 M=\min\{\mathop{\rm ess\,sup}f,\|g\|_\infty\},\qquad
 f=Me^{-V},\quad g=Me^{-W},\quad h=Me^{-U},\\
 m_g\eqdef\min W=\log\frac{M}{\|g\|_\infty}\leq0.
 \end{gathered}
\end{equation}
where $M$ is finite and positive even if $f$ is unbounded. Only $W$ is convex; the potentials may take negative values, and $V,U$ need not be convex. Both inputs have positive-volume level sets at every $q>0$, while at least one has level sets of zero volume for $q\leq0$. 

The change of variables to the log-scale works similarly to the definition of the sets $A_q,B_q,C_q$ from~\eqref{eq.log_levels}; their masses $a(q), b(q), c(q)$, and the layer-cake formulas connecting these quantities to $f,g$ from~\ref{eq.log_layer_cake} remain the same. The central objects to the height transport method are still the one-dimensional distributions
\begin{equation}\label{eq.1logconcave_measures}
    \dd\mu(q) \eqdef M a(q) e^{-q} \dd q, 
    \quad 
    \dd\nu(q) \eqdef M b(q) e^{-q} \dd q,
\end{equation}
as well as the Brenier map transporting them, which in log-scale is represented by $\mathcal{T}$ and still satisfies~\ref{eq.log_transport_definition},\eqref{eq.log_transport_identities}. 

These identities require no moment assumption on $\mu$, which is no longer log-concave in general, and Lemma~\ref{lemma.tail_ab} only applies to $\nu$, which is log-concave. In particular, we still have that if $x_\star$ minimizes $W$ then
\begin{equation}\label{eq.1logconcave_growth}
 \begin{gathered}
 B_r\subset B_s\subset x_\star+\frac{s-m_g}{r-m_g}(B_r-x_\star),\\
 b(s)\leq\left(\frac{s-m_g}{r-m_g}\right)^n b(r)
 \leq(s/r)^n b(r),\qquad 0<r\leq s.
 \end{gathered}
\end{equation}
Since $s\mapsto b(m_g+s)/s^n$ is nonincreasing, comparing its exponentially weighted integrals on $(0,s)$ and $(s,\infty)$ gives
\begin{equation}\label{eq.1logconcave_top_mass}
 \nu((m_g,m_g+s))\geq\frac1{n!}\int_0^s e^{-q}q^n\dd q,\qquad s>0.
\end{equation}

Passing to the Brunn--Minkowski stability that we will employ here, we shall use the technology developed by Barchiesi and Julin~\cite[Theorem~2]{barchiesi2017robustness} since it simultaneously yields the sharp exponent $1/2$ and the optimal dependence on $\lambda$, while considering only one convex body and a general measurable set, differently from previous approaches such as~\cite{figalli2009refined}. Their result is formulated as a stability estimate for the Euclidean concentration inequality, as follows.

\begin{theorem}[{\cite[Theorem~2]{barchiesi2017robustness}}]\label{thm.BJ}
Let $d\geq2$ and let $K\subset\R^d$ be an open, bounded and convex set containing the origin. There exists a constant $c_d>0$, depending only on $d$, such that for every $r>0$ and every measurable set $E\subset\R^d$ with $|E|=|K|$,
\[
 |E+rK|-|(1+r)K|
 \ \geq\
 c_d\max\{r^{d-1},r\}\,
 \frac{1}{|E|}\inf_{x\in\R^d}\bigl|(E+x)\sym K\bigr|^2 .
\]
\end{theorem}

In terms of the deficit $\dBM$ from~\eqref{eq.1logconcave_BM_deficits}, Theorem~\ref{thm.BJ} reads as follows for sets of equal volume.

\begin{corollary}\label{cor.BJ_lambda}
Let $d\geq2$, let $E\subset\R^d$ be measurable and let $F\subset\R^d$ be a convex body with $|E|=|F|>0$. Then
\[
 \inf_{v\in\R^d}|E\sym(F-v)|^2
 \ \leq\
 \frac{2^{d-1}}{c_d}\,\frac{|E|}{\tau}\,\dBM(E,F),
\]
where $c_d$ is the constant of Theorem~\ref{thm.BJ}.
\end{corollary}
\begin{proof}
Replacing $F$ by $F+y$ translates $\lambda E+(1-\lambda)F$ by $(1-\lambda)y$ and only shifts the variable $v$, so neither side of the estimate changes, and we may assume that the origin lies in the interior of $F$. Let $K$ be the interior of $F$. Since the boundary of a convex body is a null set, $|K|=|F|=|E|$ and $|E\sym(K-v)|=|E\sym(F-v)|$ for every $v$; moreover $\lambda E+(1-\lambda)K\subseteq\lambda E+(1-\lambda)F$, so that $\dBM(E,K)\leq\dBM(E,F)$. It therefore suffices to prove the estimate for $K$.

Set $r\eqdef(1-\lambda)/\lambda$. Then $\lambda E+(1-\lambda)K=\lambda(E+rK)$ and, as $K$ is convex, $(1+r)K=\lambda^{-1}K$. Since $|E|=|K|$, we have $(\lambda|E|^{1/d}+(1-\lambda)|K|^{1/d})^d=|K|$, and hence
\[
 \dBM(E,K)=\lambda^d|E+rK|-|K|
 =\lambda^d\bigl(|E+rK|-|(1+r)K|\bigr).
\]
By Theorem~\ref{thm.BJ}, and since $|(E+x)\sym K|=|E\sym(K-x)|$,
\[
 \dBM(E,K)\ \geq\ c_d\,\lambda^d\max\{r^{d-1},r\}\,
 \frac{1}{|E|}\inf_{v\in\R^d}|E\sym(K-v)|^2 .
\]
Finally, $\lambda^d r^{d-1}=\lambda(1-\lambda)^{d-1}$ and $\lambda^d r=\lambda^{d-1}(1-\lambda)$. If $\lambda\leq1/2$, the first of these is at least $2^{1-d}\lambda$; if $\lambda\geq1/2$, the second is at least $2^{1-d}(1-\lambda)$. In both cases $\lambda^d\max\{r^{d-1},r\}\geq2^{1-d}\tau$, which concludes the proof.
\end{proof}

To establish notation, once again we keep the parameter $\lambda$ in both deficits:
\begin{equation}\label{eq.1logconcave_BM_deficits}
 \begin{aligned}
 \dBM(E,F)&\eqdef|\lambda E+(1-\lambda)F|
       -(\lambda|E|^{1/d}+(1-\lambda)|F|^{1/d})^d,\\
 \dvol(E,F)&\eqdef\lambda|E|+(1-\lambda)|F|
       -(\lambda|E|^{1/d}+(1-\lambda)|F|^{1/d})^d.
 \end{aligned}
\end{equation}
As before, $\dvol(p,q)$ denotes the scalar expression with volumes $p,q$ and we have the following bounds
\begin{equation}\label{eq.1logconcave_scalar_volume}
 \frac\tau2(\sqrt p-\sqrt q)^2
 \leq\dvol(p,q)\leq(d-1)\tau(\sqrt p-\sqrt q)^2,
 \qquad p,q\geq0,\quad d\geq2, 
\end{equation}
where the lower bound follows by comparing the weighted power means of orders $1/d$ and $1/2$. For the upper bound, use the tangent inequality for $t^d$ at the endpoint with the larger weight, together with
\[
 x^d-y^d-dy^{d-1}(x-y)
 \leq(d-1)(x^{d/2}-y^{d/2})^2,\qquad x,y\geq0.
\]
The latter inequality is the weighted arithmetic--geometric mean inequality after expansion. Similarly, using the results of Barchiesi--Julin~\cite{barchiesi2017robustness}, we will use the following reformulation of the Brunn-Minkowski stability principle that also measures volume fluctuations between the sets being compared. 

\begin{theorem}[Brunn--Minkowski stability with one convex summand]
\label{thm.1logconcave_BM}
Let $E\subset\R^d$ be measurable with positive finite volume and let
$F\subset\R^d$ be a convex body, where $d\geq2$. Then
\begin{equation}\label{eq.1logconcave_BM_stability}
 \inf_{v\in\R^d}|E\sym(F-v)|^2
 \leq\frac{B_d^{\rm one}}{\tau}(|E|+|F|)
       [\dBM(E,F)+\dvol(E,F)],
\end{equation}
where $B_d^{\rm one}$ depends only on $d$, and one may choose
\[
 B_d^{\rm one}\leq e^{Cd}
\]
for an absolute $C$.
\end{theorem}
\begin{proof}
For equal-volume inputs, the concentration constant in~\cite[Theorem~2 and Remark~2]{barchiesi2017robustness} is bounded below by $c9^{-d}$. Its conversion to a weighted Minkowski sum introduces the factor $\tau(1-\tau)^{d-1}\geq2^{1-d}\tau$, as in the proof of Corollary~\ref{cor.BJ_lambda}. Thus the equal-volume version of~\eqref{eq.1logconcave_BM_stability} has a constant at most $e^{Cd}$.

Write $p=|E|$, $q=|F|$. If $p/q\notin[1/2,2]$, the lower bound in~\eqref{eq.1logconcave_scalar_volume} gives $\dvol(E,F)\geq c\tau(p+q)$, and the conclusion is immediate. Otherwise normalize $E,F$ by the factors $p^{-1/d},q^{-1/d}$. The interpolation parameter becomes
\[
 t=\frac{\lambda p^{1/d}}{\lambda p^{1/d}+(1-\lambda)q^{1/d}},
 \qquad \min\{t,1-t\}\geq c\tau.
\]
The relative deficit is unchanged. Rescaling the equal-volume estimate gives comparison with a translate of $\rho F$, where $\rho=(p/q)^{1/d}$. After placing the origin in $F$, the dilates $F,\rho F$ are nested and their symmetric difference has volume $|p-q|$. Finally,
\[
 |p-q|^2\leq2(p+q)(\sqrt p-\sqrt q)^2
 \leq\frac4\tau(p+q)\dvol(p,q),
\]
which absorbs this volume correction without changing the exponential order of the dimensional constant.
\end{proof}

We can now give the geometric interpretation of the smoother kernel
from~\eqref{eq.smooth_kernel}. Introduce two auxiliary coordinates
$u\in\R^2$ and the convex height window
\[
 K_r\eqdef\left\{(q,u)\in\R\times\R^2:
 |u|<1,\quad |q-r|+|u|<2\right\}.
\]
For every $r\in\R$, the chunks $E_r,F_r,H_r$ are defined in $\R^{n+3}$ by
\begin{equation}\label{eq.1logconcave_chunks}
 \begin{aligned}
 E_r&\eqdef\{(x,q,u):x\in A_q,\ (q,u)\in K_r\},\\
 F_r&\eqdef\{(x,q,u):x\in B_q,\ (q,u)\in K_r\},\\
 H_r&\eqdef\{(x,q,u):x\in C_q,\ (q,u)\in K_r\}.
 \end{aligned}
\end{equation}

These definitions were crafted in order to recover the volume interpretation of the kernel convolution from~\eqref{eq.chunks_volume_identity}. 
\begin{lemma}\label{lemma.1logconcave_volume_Rn+3}
Concerning the previous construction: 
    \begin{enumerate}
        \item for all $r\in\R$ it holds that
        $ \displaystyle
             |E_r| =
             \pi\int_{\mathbb{R}}
             \phi(q-r)^2 a(q)\dd q$, and similarly for $F_r$, $H_r$. In particular, this gives \begin{equation}\label{eq.1logconcave_chunk_convolution}
             Me^{-r}|E_r|=\mu\star\kappa(r),\qquad
             Me^{-r}|F_r|=\nu\star\kappa(r),\qquad
             Me^{-r}|H_r|=\mu_h\star\kappa(r),
        \end{equation}
        \item Since $g$ is log-concave, i.e. $W$ is convex, $F_r$ is convex; it has positive volume and is a convex body whenever $r>1/4$. 
    \end{enumerate}
\end{lemma}
\begin{proof}
To prove item (1), first notice that a simple application of Fubini's theorem gives 
\[
 \begin{aligned}
 |E_r|
 &=\int_\R |A_q|
   \left|\left\{u\in\R^2:(q,u)\in K_r\right\}\right|\dd q
 \end{aligned}
\]
From the definition of $K_r$, it follows that 
\[
    \left|
        \left\{u\in\R^2:(q,u)\in K_r\right\}
    \right|
    = 
    \left|
        \left\{
            |u| < 2 - |q-r|
        \right\}
        \cap
        B_{\mathbb{R}^2}(0,1)
    \right|. 
\]
Now, if $2 - |q-r| < 1$, or equivalently $|q-r|>1$, the intersection above is given by the ball over $\mathbb{R}^2$ of radius $2 - |q-r|$. Otherwise, the intersection is just $B_{\mathbb{R}^2}(0,1)$. As a result, the volume becomes 
\[
    \left|
        \left\{u\in\R^2:(q,u)\in K_r\right\}
    \right| 
    = 
    \begin{cases}
        \pi,& \quad |q-r| \le 1 \\
        \pi(2-|q-r|)^2,& \quad 1 \le |q-r| \le 2 \\
        0,& \text{ otherwise.}
    \end{cases}
\]
In other words, the intersection is always a ball of radius $\phi(|q-r|)$, hence the volume $\left|\left\{u\in\R^2:(q,u)\in K_r\right\}\right| = \pi\phi(|q-r|)^2$. 

Finally, this gives the desired formula
\[
    |E_r|=\pi\int_\R\phi(q-r)^2a(q)\dd q, 
\]
and similarly for $F_r, H_r$, which can also be written in the form of a convolution giving 
\[
 \begin{aligned}
 \mu\star\kappa(r)
 &=\int_\R \kappa(r-q)\dd\mu(q)
 =\pi M e^{-r}\int_\R\phi(r-q)^2a(q)\dd q
 =Me^{-r}|E_r|.
 \end{aligned}
\]
The same computation for $F_r,H_r$ yields the conclusion. Item (2) follows from a direct inspection of the definition, by means of the convexity of $W$. 
\end{proof}

Recalling that $\mu(\R)=\nu(\R)=1$ and
$\mu_h(\R)=1+\dPL$, the same relation between integral volumes of $E_r, F_r, H_r$ hold, now with the new kernel
\begin{equation}\label{eq.1logconcave_chunk_masses}
 \begin{aligned}
 M\int_\R e^{-r}|E_r|\dd r
 &=M\int_\R e^{-r}|F_r|\dd r=\|\kappa\|_{L^1},\\
 M\int_\R e^{-r}|H_r|\dd r&=(1+\dPL)\|\kappa\|_{L^1},
 \end{aligned}
\end{equation}
which makes explicit the links with the Prékopa-Leindler deficit, while the Prékopa condition and convexity of $K_r$ imply
\begin{equation}\label{eq.1logconcave_chunk_inclusion}
 \lambda E_r+(1-\lambda)F_s\subset H_{\lambda r + (1-\lambda) s} 
 \text{ for all } r,s \in\R. 
\end{equation}
Thus we can apply the quantitative stability of Brunn-Minkowski inequality from Theorem~\ref{thm.1logconcave_BM} in dimension $d=n+3$, even though $E_r$ and $H_r$ need not be convex.

\subsection{The two distance estimates}\label{sec.1logconcave_2_distance_estimates}

One important insight from the proof developed in Section~\ref{sec.quantitativePL} is that the strategy of applying Brunn-Minkowski to level-set chunks and then controlling the translations requires two types of control. While the $L^1$ estimate controls vertical translations of level-set chunks, the $W_2$ estimate controls horizontal translations, which is coherent with how the Wasserstein geometry measures distances. But under the hood of that proof is that a Hellinger distance is being bounded with $W_2$, which does not hold in general.

Lemma~\ref{lemma.1logconcave_distances} below does precisely this step, but using the elementary estimate from~\ref{lemma.smooth_height}, we bypass the Wasserstein estimate all together. As already discussed, the price to pay for this simplification was the lift to $\mathbb{R}^{n+3}$ above. 

\begin{lemma}\label{lemma.1logconcave_distances}
For the height distributions in~\eqref{eq.1logconcave_measures},
\begin{align}
 \|\mu-\nu\|_{L^1(\R)}
 =M\int_\R e^{-q}|a(q)-b(q)|\dd q
 &\leq C\sqrt{\frac{n\dPL}{\tau}},
 \label{eq.1logconcave_profile_distance}\\
 \int_\R\bigl|\sqrt{\mu\star\kappa(r)}-\sqrt{\nu\star\kappa(r)}\bigr|^2\dd r
 &\leq Cn\frac{\dPL}{\tau}.
 \label{eq.1logconcave_convolution_distance}
\end{align}
In addition, with $\mathcal S=\lambda\id+(1-\lambda)\mathcal T$, we have
\begin{equation}\label{eq.1logconcave_output_measure}
 \mu_h=\mathcal S_\sharp\mu+\xi,\qquad
 \xi\geq0,\qquad \xi(\R)=\dPL.
\end{equation}
These conclusions do not require log-concavity of either input.
\end{lemma}
\begin{proof}
For the $L^1$ estimate, it suffices to work in Euclidean heights. Indeed,
recalling $F_t=\{f>t\}$ and $G_t=\{g>t\}$, the change of variables
$t=Me^{-q}$ gives
\[
 \|\mu-\nu\|_{L^1(\R)}
 =M\int_\R e^{-q}|a(q)-b(q)|\dd q
 =\int_0^\infty\bigl||F_t|-|G_t|\bigr|\dd t.
\]
The comparison in Lemma~\ref{lem.profiles} uses the level sets $t\mapsto|F_t|$ and $t\mapsto|G_t|$ are nonincreasing, which holds for arbitrary measurable inputs. The qualitative proof of Theorem~\ref{thm.qualitative_PL} gives
\[
 \int_0^\infty |F_t|\dratio(t)\dd t\leq\dPL,
 \qquad
 \dratio(t)=\left[
   \lambda T'(t)^{(1-\lambda)/n}
   +(1-\lambda)T'(t)^{-\lambda/n}
 \right]^n-1.
\]
To estimate this remainder, we use the following elementary estimate
\begin{equation}\label{eq.1logconcave_elementary_remainder}
 \lambda e^{(1-\lambda)s}+(1-\lambda)e^{-\lambda s}-1
 \geq\frac{\tau}{4e}\min\{|s|^2,|s|\},\qquad s\in\R.
\end{equation}
Applying~\eqref{eq.1logconcave_elementary_remainder} with
$s=\log T'(t)/n$, and using $(1+x)^n\geq1+nx$ for $x\geq0$, gives
\begin{equation}\label{eq.1logconcave_euclidean_scalar_cost}
 \int_0^\infty |F_t|
 \min\{|\log T'(t)|^2,|\log T'(t)|\}\dd t
 \leq Cn\frac{\dPL}{\tau}.
\end{equation}
Here $n\min\{|s/n|^2,|s/n|\}\geq n^{-1}\min\{|s|^2,|s|\}$ accounts for the factor $n$.

The same argument found in Lemma~\ref{lem.profiles} gives
\[
 \begin{aligned}
 \|\mu-\nu\|_{L^1(\R)}
 &\leq2\left(\int_0^\infty|F_t|
       \min\{|\log T'(t)|^2,|\log T'(t)|\}\dd t\right)^{1/2}
 \leq C\sqrt{\frac{n\dPL}{\tau}},
 \end{aligned}
\]
which proves~\eqref{eq.1logconcave_profile_distance}. The only dependence on $\lambda$ comes from the scalar estimate~\eqref{eq.1logconcave_euclidean_scalar_cost}.

It remains to prove the convolution estimate, which is the only new point.
The Hellinger distance is more easily estimated in logarithmic heights. For $t=Me^{-q}$,
identities~\eqref{eq.log_transport_identities} give
\[
 \log T'(t)=\log\mathcal T'(q)+q-\mathcal T(q),
 \qquad
 \log\frac{T(t)}{tT'(t)}=-\log\mathcal T'(q).
\]
Writing $\dratio(q)=\dratio(Me^{-q})$ and
$\did(q)=\did(Me^{-q})$, estimate~\eqref{eq.1logconcave_elementary_remainder}
applied to the two scalar deficits yields
\begin{equation}\label{eq.1logconcave_scalar_cost}
 \begin{aligned}
 &\min\left\{
   |\log\mathcal T'(q)+q-\mathcal T(q)|^2,
   |\log\mathcal T'(q)+q-\mathcal T(q)|
 \right\}
 \leq\frac{Cn}{\tau}\dratio(q),\\
 &\min\{|\log\mathcal T'(q)|^2,|\log\mathcal T'(q)|\}
 \leq\frac{C}{\tau}\did(q).
 \end{aligned}
\end{equation}
Moreover, writing $ \mathcal T(q)-q =\log\mathcal T'(q) -[\log\mathcal T'(q)+q-\mathcal T(q)]$ and using also that 
$\min\{|a-b|^2,|a-b|\}
 \leq4\bigl[\min\{|a|^2,|a|\}+\min\{|b|^2,|b|\}\bigr]
$
along with the scalar deficit integral estimate from~\eqref{eq.quanli2quanti}, we obtain
\begin{equation}\label{eq.1logconcave_transport_cost}
 \int_\R\min\{|\mathcal T(q)-q|^2,|\mathcal T(q)-q|\}\dd\mu(q)
 \leq Cn\frac{\dPL}{\tau}.
\end{equation}
Since $\min\{|s|^2,1\}\leq\min\{|s|^2,|s|\}$,
Lemma~\ref{lemma.smooth_height} proves
\eqref{eq.1logconcave_convolution_distance}.

Finally, the proof of~\eqref{eq.1logconcave_output_measure} is the same as
the proof of~\eqref{eq.output_height_measure}, with
 $\mathcal S=\lambda\id+(1-\lambda)\mathcal T$. 
\end{proof}

Since at least one of $A_q,B_q$ has zero volume for $q\leq0$, the part omitted by restricting the layer-cake integral to positive heights satisfies, for every $z\in\R^n$,
\begin{equation}\label{eq.1logconcave_negative_heights}
 \begin{aligned}
 M\int_{-\infty}^0e^{-q}|A_q\sym(B_q-z)|\dd q
 &=M\int_{-\infty}^0e^{-q}|a(q)-b(q)|\dd q\\
 &\leq\|\mu-\nu\|_{L^1}
 \leq C\sqrt{\frac{n\dPL}{\tau}}.
 \end{aligned}
\end{equation}
Both chunks have positive volume for every center $r>1/4$ by the choice of $M$; no smallness assumption on $\dPL/\tau$ is needed.

We now turn these two distance estimates into a bound for the integrated
chunk deficits. This is the counterpart of
Lemma~\ref{lemma.integral_deficits_control}, with the regularity of the
kernel replacing the displacement interpolation argument.

\begin{lemma}\label{lemma.1logconcave_integrated_deficits}
For the lifted chunks defined in~\eqref{eq.1logconcave_chunks}, we have
\begin{align}
 M\int_\R e^{-r}\dvol(E_r,F_r)\dd r&\leq Cn^2\dPL,
 \label{eq.1logconcave_integrated_volume}\\
 M\int_{1/2}^\infty e^{-r}\dBM(E_r,F_r)\dd r&\leq Cn^2\dPL.
 \label{eq.1logconcave_integrated_BM}
\end{align}
\end{lemma}
\begin{proof}
For the scalar volume deficit, use~\eqref{eq.1logconcave_scalar_volume}
in dimension $d=n+3$ and the convolution identities from Lemma~\ref{lemma.1logconcave_volume_Rn+3} to obtain
\[
 \begin{aligned}
 M\int_\R e^{-r}\dvol(E_r,F_r)\dd r
 &\leq (n+2)\tau\int_\R
       \bigl|\sqrt{\mu\star\kappa(r)}-\sqrt{\nu\star\kappa(r)}\bigr|^2\dd r
 \leq Cn^2\dPL.
 \end{aligned}
\]
In particular, the factor $\tau$ from the scalar deficit cancels the
factor $\tau^{-1}$ in~\eqref{eq.1logconcave_convolution_distance}.

For the Brunn--Minkowski deficit, keep the cutoff $\beta$ used in
Lemma~\ref{lemma.integral_deficits_control}: it takes values in $[0,1]$,
vanishes on $(-\infty,1/4]$, equals one on $[1/2,\infty)$, and has
bounded first two derivatives. The corresponding test function is now
\[
 \Psi(q)\eqdef\int_\R\beta(q+s)\kappa(s)\dd s.
\]
Fubini and~\eqref{eq.1logconcave_output_measure} give the same signed
balance, with its full dependence on $\lambda$:
\begin{equation}\label{eq.1logconcave_signed_balance}
 \begin{aligned}
 &M\int_\R\beta(r)e^{-r}
       [|H_r|-\lambda|E_r|-(1-\lambda)|F_r|]\dd r\\
 &=\int_\R\Psi\dd\xi+
   \int_\R[\Psi(\mathcal S(q))-\lambda\Psi(q)
                      -(1-\lambda)\Psi(\mathcal T(q))]\dd\mu(q).
 \end{aligned}
\end{equation}
For $|q-s|\leq1$, the second-order interpolation remainder is bounded by
$\frac12\|\Psi''\|_\infty\lambda(1-\lambda)|q-s|^2$.
For $|q-s|>1$, the weighted triangle inequality and the Lipschitz bound give
$2\|\Psi'\|_\infty\lambda(1-\lambda)|q-s|$.
Thus, in both cases,
\[
 \begin{aligned}
 &|\Psi(\lambda q+(1-\lambda)s)-\lambda\Psi(q)-(1-\lambda)\Psi(s)|\\
 &\qquad\leq C\tau\min\{|q-s|^2,|q-s|\}.
 \end{aligned}
\]
Together with $\xi(\R)=\dPL$ and~\eqref{eq.1logconcave_transport_cost},
this bounds the right-hand side of~\eqref{eq.1logconcave_signed_balance}
by $Cn\dPL$. On the support of $\beta$, the chunks have positive volume,
so~\eqref{eq.1logconcave_chunk_inclusion} yields
\[
 0\leq\dBM(E_r,F_r)
 \leq |H_r|-\lambda|E_r|-(1-\lambda)|F_r|+\dvol(E_r,F_r).
\]
Integration with $M\beta(r)e^{-r}$ and the scalar deficit estimate just
proved give~\eqref{eq.1logconcave_integrated_BM}.
\end{proof}

\subsection{Proof of Theorem~\ref{thm.quantitative_PL} with one log-concave input}\label{sec.1logconcave_proof}

We can now follow the final argument of Section~\ref{sec.control_translations}.
Theorem~\ref{thm.1logconcave_BM}, Cauchy--Schwarz, the mass identities
in~\eqref{eq.1logconcave_chunk_masses}, and
Lemma~\ref{lemma.1logconcave_integrated_deficits} give
\begin{equation}\label{eq.1logconcave_all_translations}
 M\int_{1/2}^\infty e^{-r}
       \inf_{v\in\R^{n+3}}|E_r\sym(F_r-v)|\dd r
 \leq Cn\sqrt{B_{n+3}^{\rm one}}\sqrt{\frac{\dPL}{\tau}}.
\end{equation}
The only change at this stage is that the translations have three
components outside the ambient space of the original functions: the
height component and the two auxiliary coordinates of the lift.

To remove all three at once, write $v=(z,w)\in\R^n\times\R^3$ and
consider the marginals of the chunks in $y=(q,u)$:
\[
 m_E(y)=a(q)\mathbf1_{K_r}(y),\qquad
 m_F(y)=b(q)\mathbf1_{K_r}(y),\qquad
 \chi(r)\eqdef\pi\int_\R\phi(q-r)^2|a(q)-b(q)|\dd q.
\]
The spatial sections of $F_r$ and $F_r-(0,w)$ at any fixed $y$ are
two sublevel sets of $W$, or one of them is empty. They are therefore
nested, and Fubini gives
\[
 |F_r\sym(F_r-(0,w))|=\|m_F-m_F(\cdot+w)\|_1.
\]
Taking marginals decreases the $L^1$ distance, so
\[
 \|m_E-m_F(\cdot+w)\|_1\leq|E_r\sym(F_r-(z,w))|,
 \qquad \|m_E-m_F\|_1=\chi(r).
\]
Two applications of the triangle inequality consequently yield
\begin{equation}\label{eq.1logconcave_pinning}
 |E_r\sym(F_r-(z,0))|
 \leq2|E_r\sym(F_r-(z,w))|+\chi(r).
\end{equation}
This is the same marginal argument as in Section~\ref{sec.control_translations}, now with three coordinates outside $\R^n$. It only uses the nesting of the level sets. Moreover, Fubini gives
\[
 M\int_\R e^{-r}\chi(r)\dd r
 =\|\kappa\|_{L^1}\|\mu-\nu\|_{L^1}
 \leq C\sqrt{\frac{n\dPL}{\tau}}.
\]

On the smaller window $|q-r|<1$, we have $\phi(q-r)=1$. Hence we
can return to the same unweighted local error as before,
\[
 \mathcal E(r)\eqdef\inf_{z\in\R^n}
       \int_{r-1}^{r+1}|A_q\sym(B_q-z)|\dd q,
\]
since
\[
 |E_r\sym(F_r-(z,0))|
 =\pi\int_\R\phi(q-r)^2|A_q\sym(B_q-z)|\dd q
 \geq\pi\mathcal E(r).
\]
Combining~\eqref{eq.1logconcave_all_translations} and
\eqref{eq.1logconcave_pinning} therefore gives
\begin{equation}\label{eq.1logconcave_local_error}
 M\int_{1/2}^\infty e^{-r}\mathcal E(r)\dd r
 \leq Cn\bigl(1+\sqrt{B_{n+3}^{\rm one}}\bigr)
       \sqrt{\frac{\dPL}{\tau}}.
\end{equation}

It remains to join the horizontal translations. Choose
$r_k\in[k+1/2,k+3/4]$ and minimizing translations $z_k$ as in
\eqref{eq.sample_centers}. Such minimizers exist because the overlap of
the two finite-volume unweighted chunks is continuous in the horizontal
translation and tends to zero at infinity. Notice that
$[k,k+1]\subset(r_k-1,r_k+1)$, while two consecutive windows both contain
$[k+1,k+3/2]$. The overlap argument supplies
$\rho_k\in[k+1,k+3/2]$ such that
\[
 |B_{\rho_k}\sym(B_{\rho_k}+z_{k+1}-z_k)|
 \leq2[\mathcal E(r_k)+\mathcal E(r_{k+1})].
\]
Neither this selection nor the triangle inequality uses convexity of
$A_q$. The propagation argument uses only the convexity and growth of
$B_q$, which remain available in~\eqref{eq.1logconcave_growth}.
More precisely, the chord formula stated in
Section~\ref{sec.1logconcave_general_functions} below, applied to the inclusions
in~\eqref{eq.1logconcave_growth}, gives, for $0<r\leq s$,
\[
 |B_r\sym(B_r+v)|\leq|B_s\sym(B_s+v)|
 \leq(s/r)^n|B_r\sym(B_r+v)|.
\]
We retain the countably many chunks, but anchor their translations at $z_n$.
This avoids propagating every error from a height of order one to heights of order $n$.
For the increments with $j<n$, only heights $q\leq j+1\leq\rho_j$ occur in the comparison with $z_n$, so monotonicity gives a total contribution bounded by
\[
 2\sum_{j=0}^{n-1}[\mathcal E(r_j)+\mathcal E(r_{j+1})]
 \leq4e^n\sum_{k\geq0}e^{-k}\mathcal E(r_k).
\]
For $j\geq n$, the increments occur only at $q\geq j+1$, and
\[
 \begin{aligned}
 \int_{j+1}^\infty e^{-q}\left(\frac q{j+1}\right)^n\dd q
 &\leq e^{-j-1}\int_0^\infty
       e^{-t+nt/(j+1)}\dd t\\
 &=e^{-j-1}\frac{j+1}{j+1-n}
 \leq(n+1)e^{-j-1}.
 \end{aligned}
\]
Summing these bounds and using~\eqref{eq.sample_centers}, exactly as in the proof of~\eqref{eq.gluing}, gives
\[
 M\int_0^\infty e^{-q}|A_q\sym(B_q-z_n)|\dd q
 \leq Ce^n M\int_{1/2}^\infty e^{-r}\mathcal E(r)\dd r.
\]
Finally, the negative heights are bounded by~\eqref{eq.1logconcave_negative_heights}.
Combining this with~\eqref{eq.1logconcave_local_error}, we conclude that
\begin{equation}\label{eq.1logconcave_conclusion}
 \begin{aligned}
 \|f-g(\cdot+z_n)\|_{L^1(\R^n)}
 &\leq Ce^n n\bigl(1+\sqrt{B_{n+3}^{\rm one}}\bigr)
       \sqrt{\frac{\dPL}{\tau}}\\
 &\leq e^{Cn}\sqrt{\frac{\dPL}{\tau}}.
 \end{aligned}
\end{equation}
The zero-deficit case follows from the qualitative equality statement.
This proves the desired estimate with one log-concave input and with the optimal order of dependence on $\lambda$. In the squared-distance formulation, the dimensional constant obtained here is at most $e^{Cn}$ as well; the absolute constant in the exponent may change.

\subsection{Epilogue: The result for general functions}\label{sec.1logconcave_general_functions}

Once again, the careful reader will have identified that many of the important aspects of the proof carried out in this section are completely independent of either function $f,g$ being log-concave. Actually, the determining factor is which quantitative Brunn-Minkowski technology is used. In the first proof with two log-concave functions we can use for instance the original elementary argument from~\cite{figalli2009refined}. For a single log-concave input, while keeping both the sharp exponent and the optimal dependence on the parameter $\lambda$, the contribution of Barchiesi and Julin~\cite{barchiesi2017robustness} gives the necessary tool to conclude. 

At the price of at least one of the sets being convex, both of these quantitative Brunn-Minkowski technologies are unconditional on the smallness of the Brunn-Minkowski deficit. This is precisely the price that must be paid to use the Brunn--Minkowski theory for two measurable sets, due to Figalli--van Hintum--Tiba~\cite{figalli2024BM}.
For instance, their Theorem~1.5 gives, when $|E|=|F|>0$,
\begin{equation}\label{eq.1logconcave_general_BM}
 \inf_{v\in\R^d}|E\sym(F-v)|^2
 \leq\frac{C_d}{\tau}|E|\dBM(E,F),
 \qquad
 \frac{\dBM(E,F)}{|E|}<\epsilon_{d,\tau},
\end{equation}
where the smallness threshold $\epsilon_{d,\tau}$ depends on the
interpolation parameter.

Combining the Brunn-Minkowski theory from~\cite{figalli2024BM} and the quantitative height transport argument using level set chunks lifted to $\mathbb{R}^{n+3}$ described in this section, we recover the sharp square-root estimate between arbitrary integrable inputs at every fixed $\lambda\in(0,1)$
from~\cite[Theorem~1.16, $p=0$]{FigalliVanHintumTibaPL}\footnote{We use the theorem numbering of the second arXiv version (6 October 2026). Corollary~1.7 there additionally supplies a common log-concave approximant; the assertion proved here concerns the distance between the two inputs.}, i.e.
\[
 \inf_{z\in\R^n}\|f-g(\cdot+z)\|_1
 \leq C_{n,\lambda}\sqrt{\dPL(\lambda)}.
\]
This claim requires more than the replacement of Theorem~\ref{thm.1logconcave_BM} by~\eqref{eq.1logconcave_general_BM}: Sections~\ref{sec.1logconcave_review_notation}--\ref{sec.1logconcave_proof} also use log-concavity in choosing the reference height, applying stability to the convex chunk $F_r$, and propagating overlap estimates through~\eqref{eq.1logconcave_growth}. Here we choose the reference height by mass and retain the longer windows needed to propagate translations through convex hulls. The additional input is Theorem~\ref{thm.sug_FvHT}, the convex-hull form of Brunn--Minkowski stability from~\cite{figalli2024BM}.

As usual, it suffices to treat
\[
 0<\dPL\le c_0,
\]
where $c_0\in(0,1]$ depends only on $n$ and $\lambda$ and is fixed in Lemma~\ref{lemma.sug_height}.

The scalar estimate in the proof of Lemma~\ref{lemma.1logconcave_distances} allows the profile coefficient $P_n=4\sqrt{en}$. We shall use repeatedly that
\[
 \|\mu-\nu\|_{L^1(\R)}\ \le\ P_n\sqrt{\frac{\dPL}{\tau}}
\]
holds for arbitrary functions (see Lemma~\ref{lemma.1logconcave_distances}), and that it does not depend on the reference height $M$ used to define $\mu$ and $\nu$, as
\[
 \|\mu-\nu\|_{L^1(\R)}=\int_0^\infty\bigl||F_t|-|G_t|\bigr|\,\dd t .
\]

We shall rely on the following form of the Brunn--Minkowski stability for measurable sets.

\begin{theorem}[{\cite[Theorem~1.3 and Remark~1.8]{figalli2024BM}}]\label{thm.sug_FvHT}
For every $d\in\mathbb N$ and $t\in(0,1/2]$ there are constants $H_d\ge1$, depending only on $d$, and $\epsilon_{d,t}>0$, depending only on $d$ and $t$ and nondecreasing in $t\in(0,1/2]$, such that the following holds. Let $A,B\subset\R^d$ be measurable sets with $|A|=|B|\in(0,\infty)$, and define $\gamma\ge0$ by
\[
 |tA+(1-t)B|=(1+\gamma)\,|A| .
\]
If $\gamma<\epsilon_{d,t}$, then there are $x\in\R^d$ and a convex set $L\supseteq A\cup(B+x)$ such that
\[
 |L\setminus A|+|L\setminus(B+x)|\ \le\ H_d\,t^{-1/2}\gamma^{1/2}\,|A| .
\]
\end{theorem}

To record the dependence inherited from this theorem, put
\begin{equation}\label{eq.sug_BM_constant}
 \epsilon_*\eqdef\min\{1,\epsilon_{d,\tau/2}\},\qquad
 \mathfrak B_{d,\tau}\eqdef
 1+\frac{H_d^2}{\tau}+\epsilon_*^{-1}.
\end{equation}
The constants denoted by $C$ below are numerical; dependence on the dimension and on the Brunn--Minkowski input will be displayed explicitly.

For general inputs, the minimum of their essential suprema can be infinite. We choose the reference height by mass instead, so that the lowest chunks are nonempty and the contribution of negative heights is controlled.

\begin{lemma}\label{lemma.sug_height}
Set $c_0=\tau/(16P_n^2)$. There is $M>0$ with the following property. Define $V,W,U$ by
\[
 f=Me^{-V},\qquad g=Me^{-W},\qquad h=Me^{-U},
\]
and the sets $A_q,B_q,C_q$, the profiles $a,b$ and the measures $\mu,\nu$ by the formulas of Section~\ref{sec.1logconcave_review_notation}, the choice of $M$ as the minimum of the two essential suprema being replaced by the present mass normalization. Then
\[
 \nu\bigl((-\infty,0]\bigr)=2P_n\sqrt{\frac{\dPL}{\tau}},\qquad a(0)>0,\qquad b(0)>0,
\]
and, for every $z\in\R^n$,
\begin{equation}\label{eq.sug_negative_heights}
 M\int_{-\infty}^0e^{-q}\,\bigl|A_q\sym(B_q-z)\bigr|\,\dd q\ \le\ 5P_n\sqrt{\frac{\dPL}{\tau}} .
\end{equation}
\end{lemma}

\begin{proof}
Choose $M$ so that the part of $g$ above $M$ has twice the mass allowed by the profile estimate:
\[
 \int_M^\infty|G_s|\,\dd s=2P_n\sqrt{\frac{\dPL}{\tau}} .
\]
 The change of variables $s=Me^{-q}$ sends the heights $q\le0$ to the levels $s\ge M$, so this plainly corresponds to the identity for $\nu((-\infty,0])$. Since this mass is positive and $b$ is nondecreasing, $b(0)>0$. If $f$ vanished above $M$, that is, if $a(0)=0$, the profiles would differ by at least this mass,
\[
 \|\mu-\nu\|_{L^1}\ \ge\ \nu\bigl((-\infty,0]\bigr)=2P_n\sqrt{\frac{\dPL}{\tau}}\ >\ P_n\sqrt{\frac{\dPL}{\tau}},
\]
which contradicts the profile estimate, as $\dPL>0$. Finally, the levels above $M$ cost at most the mass of both functions there, since $|A_q\sym(B_q-z)|\le a(q)+b(q)$:
\[
 M\int_{-\infty}^0e^{-q}\,\bigl|A_q\sym(B_q-z)\bigr|\,\dd q
 \ \le\ \mu\bigl((-\infty,0]\bigr)+\nu\bigl((-\infty,0]\bigr)
 \ \le\ 2\,\nu\bigl((-\infty,0]\bigr)+\|\mu-\nu\|_{L^1} .
\]
Using the bound on $\nu((-\infty,0])$ and Lemma \ref{lemma.1logconcave_distances}, the right-hand side is at most $5P_n\sqrt{\dPL/\tau}$, finishing the proof.
\end{proof}

From now on, we fix $M$ to be the reference height of Lemma~\ref{lemma.sug_height}. The Prékopa condition still gives $\lambda A_q+(1-\lambda)B_s\subseteq C_{\lambda q+(1-\lambda)s}$ for all $q,s\in\R$, and Lemma~\ref{lemma.1logconcave_distances} holds for this $M$, as its proof uses the reference height only through the change of variables $t=Me^{-q}$.

Moving on, we keep the chunks of~\eqref{eq.1logconcave_chunks}, on longer windows. Fix $\ell\eqdef 4n+1,$ the half-length of the part of the window on which the kernel is flat, and redefine
\[
 \begin{aligned}
 K^{}_r&\eqdef\Bigl\{(q,u)\in\R\times\R^2:\ |u|<1,\ |q-r|+|u|<\ell+1\Bigr\},\\
 \kappa(s)&\eqdef\pi e^{-s}\min\bigl\{1,(\ell+1-|s|)_+\bigr\}^2 .
 \end{aligned}
\]
For $\ell=1$ these match the previous definitions of window $K_r$ and kernel $\kappa$ used so far. The chunks are
\[
 E_r\eqdef\Bigl\{(x,q,u)\in\R^n\times\R\times\R^2:\ x\in A_q,\ (q,u)\in K_r\Bigr\},
\]
and $F_r$, $H_r$ are defined in the same way with $B_q$ and $C_q$. The computation of Lemma~\ref{lemma.1logconcave_volume_Rn+3} gives
\[
 Me^{-r}|E_r|=\mu\star\kappa (r),\qquad
 Me^{-r}|F_r|=\nu\star\kappa (r),\qquad
 Me^{-r}|H_r|=\mu_h\star\kappa (r),
\]
so the mass identities~\eqref{eq.1logconcave_chunk_masses} hold with the adjusted value of $\|\kappa\|_{L^1}$ for the current longer window kernel. As the windows are convex, $\lambda E_r+(1-\lambda)F_s\subseteq H_{\lambda r+(1-\lambda)s}$ for all $r,s\in\R$, and Lemma~\ref{lemma.smooth_height} holds for $\kappa $, whose square root is Lipschitz with compact support. 

The cost of these longer windows is explicit:
\[
 \|\kappa\|_1+\| (\sqrt\kappa)'\|_2^2\leq C e^{\ell+1}.
\]
Together with the scalar bounds in Lemma~\ref{lemma.1logconcave_distances}, this accounts for the dimensional constants in the following integrated estimates.

For measurable sets $E,F,H\subset\R^d$ of finite positive volume, the Brunn--Minkowski stability will be applied to the \emph{deficit of the triple}
\[
 \delta(E,F;H)\eqdef\Bigl[\,|H|-\bigl(\lambda|E|^{1/d}+(1-\lambda)|F|^{1/d}\bigr)^d\Bigr]+\dvol(E,F),
\]
which is at least $\dBM(E,F)+\dvol(E,F)$ when $\lambda E+(1-\lambda)F\subseteq H$. The reason why we use $H$ rather than the Minkowski combination itself is mainly because $|H_r|$ depends continuously on $r$. For the chunks, their \emph{deficit centred at $r$} is defined as
\[
 \delta_{\rm ch}(r)\eqdef\delta\bigl(E_r,F_r;H_r\bigr)=\Bigl[\,|H_r|-\lambda|E_r|-(1-\lambda)|F_r|\Bigr]+2\,\dvol(E_r,F_r),\qquad d=n+3 .
\]
The proof of Lemma~\ref{lemma.1logconcave_integrated_deficits} applies then in this context as well, and it thus yields bounds to integrals of both terms against $M\beta(r)e^{-r}\,\dd r$, where $\beta$ is the cutoff used there. Since $\delta_{\rm ch}(r)\ge0$ for $r\ge1/4$, by the Brunn--Minkowski inequality, while $\beta$ vanishes on $(-\infty,1/4]$ and equals one on $[1/2,\infty)$, this gives
\begin{equation}\label{eq.sug_deficits}
 M\int_{1/2}^\infty e^{-r}\,\delta_{\rm ch}(r)\,\dd r\leq C(n+3)^3e^{\ell+1}\dPL .
\end{equation}

Indeed, the scalar volume bound can be taken with upper coefficient $d^2\tau/2$, and the truncated transport cost with coefficient $16e(n+1)/\tau$. The factor $\tau$ cancels, while the kernel norms and the first two derivatives of the test function in~\eqref{eq.1logconcave_signed_balance} are bounded by $Ce^{\ell+1}$.

As in the previous sections, in order to use the stability of Brunn--Minkowski from Theorem~\ref{thm.sug_FvHT} we must first turn it into an assertion about sets of \emph{unequal} volume, as its original formulation only compares sets of equal volume. For convex bodies, the proof of Theorem~\ref{thm.BM_stability} rescales one body to the volume of the other, and after a translation the rescaled body is nested with the original. On the other hand, a rescaled measurable set need not be close to the original set. 

In order to circumvent this issue, we normalize both sets to have unit volume. By homogeneity of the Minkowski combination, for measurable sets $E,F\subset\R^d$ of finite positive volume,
\begin{equation}\label{eq.sug_normalization}
 \lambda E+(1-\lambda)F=\bigl(\lambda|E|^{1/d}+(1-\lambda)|F|^{1/d}\bigr)\Bigl(t\,|E|^{-1/d}E+(1-t)\,|F|^{-1/d}F\Bigr),
\end{equation}
where the new interpolation parameter is
\[
 t\eqdef\frac{\lambda|E|^{1/d}}{\lambda|E|^{1/d}+(1-\lambda)|F|^{1/d}} .
\]
Normalizing changes the interpolation parameter from $\lambda$ to $t$, but leaves the relative Brunn--Minkowski deficit unchanged. The rescaling back is then applied to the common convex set given by Theorem~\ref{thm.sug_FvHT}, whose dilates are nested. One is then able to obtain the following different volume version of the Brunn--Minkowski stability from that idea: 

\begin{lemma}\label{lemma.sug_BM}
Let $d\ge2$, and let $E,F,H\subset\R^d$ be measurable sets of finite positive volume with $\lambda E+(1-\lambda)F\subseteq H$. Then $\delta(E,F;H)\ge0$. With
\[
 c=\min\left\{\frac\tau{100},\frac{\epsilon_*}{4},
                       \frac\tau{64H_d^2}\right\},
 \qquad c^{-1}\leq100\mathfrak B_{d,\tau},
\]
the following hold for a numerical constant $C$.
\begin{enumerate}
 \item[(i)] For every such triple,
 \[
 \inf_{x\in\R^d}\bigl|E\sym(F+x)\bigr|^2\leq C\mathfrak B_{d,\tau}\bigl(|E|+|F|\bigr)\delta(E,F;H).
 \]
 \item[(ii)] If $\delta(E,F;H)\le c\,|E|$, then the convex hull $\co(F)$ has finite volume and
 \[
 |\co(F)\setminus F|\leq C\sqrt{\mathfrak B_{d,\tau}}\sqrt{\delta(E,F;H)\,|E|}\, .
 \]
 \item[(iii)] If $\delta(E,F;H)>c\,|E|$, then
 \[
 |E|+|F|\leq C\mathfrak B_{d,\tau}\delta(E,F;H).
 \]
\end{enumerate}
\end{lemma}

\begin{proof}
The bracket in $\delta(E,F;H)$ is nonnegative by the Brunn--Minkowski inequality, and by Lemma~\ref{lemma.scalar_Hellinger} we have 
\[
 \dvol(p,q)\ \ge\ \frac\tau2\bigl(\sqrt p-\sqrt q\bigr)^2,\qquad p,q\ge0 .
\]
If one volume is at least twice the other, this gives $\delta(E,F;H)\ge\dvol(E,F)\ge\frac{\tau}{48}(|E|+|F|)$, hence (iii), and (i) follows from $\inf_x|E\sym(F+x)|\le|E|+|F|$. Otherwise $\min\{t,1-t\}\ge\tau/2$, and by~\eqref{eq.sug_normalization} the normalized sets $|E|^{-1/d}E$ and $|F|^{-1/d}F$ have, at the parameter $t$, the equal-volume deficit
\[
 \gamma\ =\ \frac{|\lambda E+(1-\lambda)F|}{\bigl(\lambda|E|^{1/d}+(1-\lambda)|F|^{1/d}\bigr)^d}-1\ \le\ \frac{2\,\delta(E,F;H)}{|E|} .
\]
If $\delta(E,F;H)\le c|E|$, our choice of $c$ makes the volumes comparable and gives $\gamma\le2c<\epsilon_*$. Moreover, the convex-hull error supplied by Theorem~\ref{thm.sug_FvHT} is at most $2H_d\sqrt{c/\tau}\le1/4$, so its convex set $L$ has $|L|\le2$. Rescaling,
\[
 \bigl|\,|E|^{1/d}L\setminus E\,\bigr|+\bigl|\,|F|^{1/d}L\setminus\bigl(F+|F|^{1/d}x\bigr)\bigr|\leq\frac{4H_d}{\sqrt\tau}\sqrt{\delta(E,F;H)\,|E|}\, .
\]
Then (ii) holds because $\co(F)\subseteq|F|^{1/d}L-|F|^{1/d}x$, and (i) holds because the dilates $|E|^{1/d}L$ and $|F|^{1/d}L$ are nested after a translation and differ in volume by $\bigl||E|-|F|\bigr|\,|L|$, where
\[
 \bigl||E|-|F|\bigr|^2\ \le\ \frac4\tau\bigl(|E|+|F|\bigr)\,\dvol(E,F).
\]
If $\delta(E,F;H)>c|E|$ and the volumes are comparable, then $|E|+|F|<3|E|<\frac3c\,\delta(E,F;H)$, which is (iii) and again gives (i).
\end{proof}

From now on write $\mathfrak B=\mathfrak B_{n+3,\tau}$. We repeat the first part of Section~\ref{sec.1logconcave_proof} with longer windows. Lemma~\ref{lemma.sug_BM}(i), applied to $E_r$, $F_r$, $H_r$ with $d=n+3$, the Cauchy--Schwarz inequality, the mass identities and~\eqref{eq.sug_deficits} give
\[
 M\int_{1/2}^\infty e^{-r}\inf_{w\in\R^{n+3}}\bigl|E_r\sym(F_r-w)\bigr|\,\dd r\leq Ce^{Cn}\sqrt{\mathfrak B}\sqrt{\dPL}.
\]
The pinning argument leading to~\eqref{eq.1logconcave_pinning} applies verbatim to the longer windows, with a profile term that integrates against $Me^{-r}\,\dd r$ to $\|\kappa \|_{L^1}\|\mu-\nu\|_{L^1}$. As the kernel is flat on $[-\ell,\ell]$, the local error on the flat part of the window,
\[
 \mathcal E(r)\eqdef\inf_{z\in\R^n}\int_{r-\ell}^{r+\ell}\bigl|A_q\sym(B_q-z)\bigr|\,\dd q ,
\]
which is upper semicontinuous and hence measurable, therefore satisfies
\begin{equation}\label{eq.sug_local}
 M\int_{1/2}^\infty e^{-r}\,\mathcal E(r)\,\dd r\leq Ce^{Cn}\sqrt{\mathfrak B}\sqrt{\dPL}.
\end{equation}

Finally, one needs to put the horizontal translations together, as previously done in other parts of the proof.  In Section~\ref{sec.1logconcave_proof} the overlap estimate is propagated to all heights by the inclusions~\eqref{eq.1logconcave_growth}. For general $g$ we use instead the convex hulls provided by Lemma~\ref{lemma.sug_BM}(ii), which allow us to perform essentially the same mechanism: inside one chunk, the sections of a convex set grow at most like dilates, and the long windows make the resulting losses summable against the weights $e^{-q}$. For a set $S\subseteq\R^n\times\R\times\R^2$ and $(q,u)\in\R\times\R^2$ we write
\[
 S_{(q,u)}\eqdef\bigl\{x\in\R^n:\ (x,q,u)\in S\bigr\}
\]
for its section. For instance, the section of $F_r$ at $(q,u)\in K_r$ is $B_q$. We shall need the following elementary bounds: let $z\in\R^n$.
\begin{enumerate}
 \item[(i)] If $S\subseteq L\subseteq\R^n$ are measurable and $|L|<\infty$, then
 \begin{align*}
 |S\sym(S+z)|&\le|L\sym(L+z)|+2|L\setminus S|, \cr
 |L\sym(L+z)|&\le|S\sym(S+z)|+2|L\setminus S| .
 \end{align*}
 \item[(ii)] Let $L\subset\R^n\times\R\times\R^2$ be a bounded convex set. If $u\in\R^2$, $q_0<q_1\le q_2$ and $L_{(q_0,u)}\neq\varnothing$, then
 \[
 \bigl|L_{(q_2,u)}\sym(L_{(q_2,u)}+z)\bigr|\ \le\ \Bigl(\frac{q_2-q_0}{q_1-q_0}\Bigr)^{n}\,\bigl|L_{(q_1,u)}\sym(L_{(q_1,u)}+z)\bigr| .
 \]
\end{enumerate}

Indeed, item (i) follows from elementary inclusions, while convexity of $L$, applied to a point $x_0\in L_{(q_0,u)}$ and to the points of $L_{(q_2,u)}$, gives
\[
 L_{(q_2,u)}\ \subseteq\ x_0+\frac{q_2-q_0}{q_1-q_0}\,\bigl(L_{(q_1,u)}-x_0\bigr).
\]
For a bounded convex set $P\subset\R^n$ and $z\neq0$, the chord formula
\[
 |P\sym(P+z)|\ =\ 2\int_{z^\perp}\min\Bigl\{|z|,\ \mathcal H^1\bigl(P\cap(y+\R z)\bigr)\Bigr\}\,\dd y
\]
is nondecreasing under inclusion, invariant under translations, and gives $|\Lambda P\sym(\Lambda P+z)|\le\Lambda^n|P\sym(P+z)|$ for $\Lambda\ge1$. Together with the inclusion above, this proves (ii).\\

Our next task is to use these two estimates in order to adapt our machinery to the general case. We now claim that there is $z_0\in\R^n$ such that
\[
 M\int_0^\infty e^{-q}\,\bigl|A_q\sym(B_q-z_0)\bigr|\,\dd q\leq Ce^{Cn}\mathfrak B\sqrt{\dPL}.
\]

In order to join the horizontal translations, we need to account for the error that is obtained when one passes from the chunks to their convex hulls. To that extent, call a center $r\ge1/2$ \emph{good} if $\delta_{\rm ch}(r)\le c\,|E_r|$, with $c$ the constant of Lemma~\ref{lemma.sug_BM} for $d=n+3$, and \emph{bad} otherwise. By that lemma,
\begin{equation}\label{eq.sug_good_bad}
 \begin{aligned}
 |\co(F_r)\setminus F_r|&\ \le\ C\sqrt{\mathfrak B}\sqrt{\delta_{\rm ch}(r)\,|E_r|} &&\text{if $r$ is good},\\
 |F_r|&\ \le\ C\mathfrak B\,\delta_{\rm ch}(r) &&\text{if $r$ is bad},
 \end{aligned}
\end{equation}
Put
\[
 \mathcal D(r)\eqdef
 \sqrt{\mathfrak B}\sqrt{\delta_{\rm ch}(r)|E_r|}
       +\mathfrak B\delta_{\rm ch}(r).
\]
By Cauchy--Schwarz, the mass identities and~\eqref{eq.sug_deficits},
\[
 M\int_{1/2}^\infty e^{-r}\mathcal D(r)\dd r
 \leq Ce^{Cn}\bigl(\sqrt{\mathfrak B}\sqrt{\dPL}
                         +\mathfrak B\dPL\bigr)
 \leq Ce^{Cn}\mathfrak B\sqrt{\dPL}.
\]

For each integer $k\ge0$, choose a center $r_k\in[k+1/2,k+3/4]$ at which the sum $\mathcal E+\mathcal D$ does not exceed its average over this interval, and a translation $z_k$ attaining the infimum in $\mathcal E(r_k)$. As $e^{-k}\le e^{3/4}e^{-r}$ for $r\le k+3/4$ and these intervals are disjoint, \eqref{eq.sug_local} and the previous step give
\begin{equation}\label{eq.sug_sum}
 M\sum_{k\ge0}e^{-k}\Bigl[\mathcal E(r_k)+\mathcal D(r_k)\Bigr]\leq Ce^{Cn}\mathfrak B\sqrt{\dPL}.
\end{equation}
The flat part of the $k$-th window contains $[k+3/4-\ell,\,k+1/2+\ell]$. Hence $[0,\ell]$ lies in the window of $r_0$, each $[m,m+1]$ with $m\ge1$ lies in the windows of $r_{m-1}$ and $r_m$, and $[m-1,m+1]$ lies in the window of $r_{m-\ell+1}$ when $m\ge\ell$.

Let now $m\ge\ell$, and let $r=r_{m-\ell+1}$, whose window is flat on $[m-1,m+1]$. We claim that, for every $z\in\R^n$,
\begin{equation}\label{eq.sug_transfer}
 \begin{aligned}
 \int_m^{m+1}\bigl|B_q\sym(B_q+z)\bigr|\,\dd q
 &\ \le\ e^{1/2}\int_{m-1}^{m}\bigl|B_q\sym(B_q+z)\bigr|\,\dd q\\
 &\qquad+C\mathcal D(r).
 \end{aligned}
\end{equation}
Indeed, we divide into two cases: first, if $r$ is bad, its chunk is small, and so: 
\[
 \int_m^{m+1}\bigl|B_q\sym(B_q+z)\bigr|\,\dd q\ \le\ 2\int_m^{m+1}b(q)\,\dd q\ \le\ \frac2\pi\,|F_r|\leq C\mathfrak B\delta_{\rm ch}(r)\leq C\mathcal D(r).
\]
If $r$ is good, $\co(F_r)$ is a bounded convex subset of $\R^n\times K_r$ whose sections at $(q,u)$, for $|u|\le1/2$ and $|q-r|\le\ell$, contain $B_q$. At $q_0=\max\{r-\ell,0\}$, these sections are nonempty for every such $u$: choose $x_0\in B_0$, and observe that $(x_0,q_0,u)\in F_r$ since $q_0\ge0$ and $|q_0-r|\le\ell$. For $q_1\in[m-1,m]$ and $q_2\in[m,m+1]$ we have $q_1-q_0\ge\ell-1=4n$ and $q_2-q_1\le2$, so the growth factor in the bound (ii) above is at most
\begin{equation}\label{eq.sug_ratio}
 \Bigl(\frac{q_2-q_0}{q_1-q_0}\Bigr)^n\ \le\ \Bigl(1+\frac{1}{2n}\Bigr)^n\ \le\ e^{1/2}.
\end{equation}
Note that it is at \emph{exactly} this point that we needed to adapt the definition of $K_r, \kappa$. Using the elementary bounds (i) and (ii) above integrated in $q_1\in[m-1,m]$ and $q_2\in[m,m+1]$ and averaged over $|u|\le1/2$, gives~\eqref{eq.sug_transfer}: by Fubini's theorem the averaged section errors are at most $C\,|\co(F_r)\setminus F_r|$, which~\eqref{eq.sug_good_bad} bounds.\\

Finally, as $[m,m+1]$ lies in the windows of $r_{m-1}$ and $r_m$ for $m\ge1$, this gives
\[
 \int_m^{m+1}\bigl|B_q\sym(B_q+z_m-z_{m-1})\bigr|\,\dd q\ \le\ \mathcal E(r_{m-1})+\mathcal E(r_m)\qquad(m\ge1),
\]
and, as $[\ell-1,\ell]$ also lies in the window of $r_0$,
\[
 \int_{\ell-1}^{\ell}\bigl|B_q\sym(B_q+z_{\ell-1}-z_0)\bigr|\,\dd q\ \le\ \mathcal E(r_{\ell-1})+\mathcal E(r_0).
\]
Now let $m\ge\ell$. The integral $\int_m^{m+1}|B_q\sym(B_q+z_m-z_0)|\,\dd q$ is the cost, at the heights $[m,m+1]$, of replacing the translation $z_m$ by $z_0$. It then follows by the triangle inequality and~\eqref{eq.sug_transfer} with $z=z_{m-1}-z_0$, that
\begin{equation}\label{eq.sug_recursion}
 \begin{aligned}
 \int_m^{m+1}\bigl|B_q\sym(B_q+z_m-z_0)\bigr|\,\dd q
 &\ \le\ e^{1/2}\int_{m-1}^{m}\bigl|B_q\sym(B_q+z_{m-1}-z_0)\bigr|\,\dd q\\
 &\qquad+C\mathcal D(r_{m-\ell+1})\\
 &\qquad+\mathcal E(r_{m-1})+\mathcal E(r_m).
 \end{aligned}
\end{equation}

This implies at once that
\[
 \begin{aligned}
 \sum_{m\ge\ell}e^{-m}\int_m^{m+1}\bigl|B_q\sym(B_q+z_m-z_0)\bigr|\,\dd q
 &\ \le\ e^{-1/2}\sum_{m\ge\ell-1}e^{-m}\int_m^{m+1}\bigl|B_q\sym(B_q+z_m-z_0)\bigr|\,\dd q\\
 &\qquad+C\sum_{k\ge0}e^{-k}\Bigl[\mathcal E(r_k)+\mathcal D(r_k)\Bigr].
 \end{aligned}
\]
The sum on the right is the sum on the left plus one term with $m=\ell-1$. The sum is finite because $|B_q\sym(B_q+z)|\leq2b(q)$ and $M\int e^{-q}b(q)\dd q=1$. Since $e^{-1/2}<1$, it can be absorbed into the left-hand side. By the previous considerations the term with $m=\ell-1$ is at most $e^{-(\ell-1)}\mathcal E(r_{\ell-1})+\mathcal E(r_0)$, which is also bounded by the series. Therefore
\[
 \sum_{m\ge\ell}e^{-m}\int_m^{m+1}\bigl|B_q\sym(B_q+z_m-z_0)\bigr|\,\dd q\ \le\ C\sum_{k\ge0}e^{-k}\Bigl[\mathcal E(r_k)+\mathcal D(r_k)\Bigr].
\]
Finally, the triangle inequality and the window inclusions give $\int_0^\ell|A_q\sym(B_q-z_0)|\,\dd q\le\mathcal E(r_0)$ and, for $m\ge\ell$,
\[
 \int_m^{m+1}\bigl|A_q\sym(B_q-z_0)\bigr|\,\dd q\ \le\ \mathcal E(r_m)+\int_m^{m+1}\bigl|B_q\sym(B_q+z_m-z_0)\bigr|\,\dd q .
\]
Since $e^{-q}\le1$ on $[0,\ell]$ and $e^{-q}\le e^{-m}$ on $[m,m+1]$, multiplying by $M$ and using~\eqref{eq.sug_sum} concludes the proof.

By the modifications above, Cavalieri's principle, written in the logarithmic scale as in~\eqref{eq.log_layer_cake}, splits the distance into the heights $q\le0$, bounded by~\eqref{eq.sug_negative_heights}, and the heights $q>0$, bounded by the claim above:
\[
 \|f-g(\cdot+z_0)\|_{L^1}=M\int_\R e^{-q}\,\bigl|A_q\sym(B_q-z_0)\bigr|\,\dd q
 \ \le\ 5P_n\sqrt{\frac{\dPL}{\tau}}+Ce^{Cn}\mathfrak B\sqrt{\dPL}.
\]
Since $P_n=4\sqrt{en}$, $\ell=4n+1$ and $\mathfrak B\geq\tau^{-1}$, this proves the general result with
\begin{equation}\label{eq.sug_dimension_dependence}
 C_{n,\lambda}\leq Ce^{Cn}\left(
 1+\frac{H_{n+3}^2}{\tau}
 +\frac1{\min\{1,\epsilon_{n+3,\tau/2}\}}
 \right).
\end{equation}
The same bound covers $\dPL>c_0$ by $\|f-g(\cdot+z)\|_1\leq2$. Thus the height-transport and gluing steps have at most exponential dimensional cost; the remaining dependence is that of the coefficient and smallness threshold in Theorem~\ref{thm.sug_FvHT}. 

\section{Application to the log-Sobolev and Talagrand inequalities}
\label{sec.applications} 

In this section, once armed with the Prékopa-Leindler technology crafted so far, we fully develop the application of our theory to the potential-independent, sharp $L^1$-stability of the log-Sobolev and Talagrand inequalities developed in~\cite{kristaly2026quantitative}. As in the introduction, we emphasize the importance of the sharp exponent $1/2$ and the optimal dependence of the stability constant on $\lambda$, with a single log-concave input.

\subsection{The log-Sobolev inequality: distributions and wave functions}
\label{subsec.app_log_sobolev}
Let
\[
 \dd\mu_V=e^{-V}\dd x,\qquad \int_{\R^n}e^{-V}\dd x=1,
 \qquad V\in \mathscr{C}^2(\R^n),\qquad \nabla^2V\geq \kappa I_n,
 \quad \kappa>0.
\]
For $z\in\R^n$, write $\tau_z(x)=x+z$ and denote the translated
reference probability by $\mu_{V,z}=(\tau_z)_\sharp\mu_V$.  The family of
normalized translation profiles selected by the stability theorem is
\begin{equation}\label{eq.app_reference_translates}
 \mathcal T_V
 \eqdef
 \left\{
   \frac{\dd\mu_{V,z}}{\dd\mu_V}
   =e^{V-V(\cdot-z)}:z\in\R^n
 \right\}.
\end{equation}
Krist\'aly and P{\^\i}rvuceanu
\cite[Theorems~1.1 and~1.3(ii)]{kristaly2026quantitative}
also identify the exact normalized optimizer set for both the distributional
log-Sobolev inequality and Talagrand's inequality. In our notation these are
\[
 \mathcal O_{V,\kappa}
 \eqdef
 \left\{e^{V-V(\cdot-z)}:z\in\mathcal K_{V,\kappa}\right\},
 \qquad
 \mathcal K_{V,\kappa}
 \eqdef
 \bigcap_{x\in\R^n}\Ker\bigl(\nabla^2V(x)-\kappa I_n\bigr).
\]
Thus $\mathcal O_{V,\kappa}\subseteq\mathcal T_V$. For the standard
Gaussian, we use $\gamma_n$ both for the probability measure and its
Lebesgue density
\[
 \gamma_n(x)=(2\pi)^{-n/2}e^{-|x|^2/2}
 \text{ and }
 \mathcal O_{\gamma_n}
 \eqdef
 \left\{e^{z\cdot x-|z|^2/2}:z\in\R^n\right\}.
\]

We first formulate log-Sobolev entirely in terms of distributions. If $\varrho\ll\mu_V$ is a probability measure, define its relative entropy
and Fisher information by
\begin{equation}\label{eq.app_entropy_fisher}
 \begin{aligned}
 H(\varrho\mid\mu_V)
 =\int \log\left(\frac{\dd\varrho}{\dd\mu_V}\right)\dd\varrho,
  \quad 
 I(\varrho\mid\mu_V)
 =\int\left|\nabla\log\left(
       \frac{\dd\varrho}{\dd\mu_V}
     \right)\right|^2\dd\varrho.
 \end{aligned}
\end{equation}
The Bakry--\'Emery log-Sobolev inequality and its deficit read as follows:
\begin{equation}\label{eq.app_density_LSI_deficit}
 H(\varrho\mid\mu_V)\leq\frac1{2\kappa}I(\varrho\mid\mu_V),
 \qquad
 \dLSI(\varrho\mid\mu_V)
 \eqdef\frac1{2\kappa}I(\varrho\mid\mu_V)-H(\varrho\mid\mu_V).
\end{equation}

The same inequality has also a traditional functional formulation: every
nonnegative $f\in W^{1,2}(\mu_V)$, normalized by
$\int f^2\dd\mu_V=1$, determines the probability distribution
\begin{equation}\label{eq.app_wave_representation}
 \dd\varrho=f^2\dd\mu_V,\qquad
 f=\sqrt{\frac{\dd\varrho}{\dd\mu_V}}.
\end{equation}
In this representation,
$I(\varrho\mid\mu_V)=4\int|\nabla f|^2\dd\mu_V$, and therefore
\begin{equation}\label{eq.app_wave_LSI_deficit}
 \dLSI(\varrho\mid\mu_V)
 =\frac2\kappa\int|\nabla f|^2\dd\mu_V
  -H(\varrho\mid\mu_V).
\end{equation}
Borrowing the terminology of quantum mechanics, one may regard $f$ as a
real nonnegative wave function and $f^2\mu_V$ as its probability
distribution.  This wave-function formulation is the one studied in~\cite{dolbeault2025sharp} for the standard Gaussian measure, whereas the stability version formulated directly in terms of probability densities can be derived for a much wider class of reference log-concave measures under the same universal stability constant, as observed by~\cite{kristaly2026quantitative}.

\begin{proposition}
\label{prop.app_sharp_LS}
Let $\varrho\ll\mu_V$ be a probability measure with finite Fisher
information.  Then
\begin{equation}\label{eq.app_sharp_LS}
 \dist_{L^1(\mu_V)}\left(
   \frac{\dd\varrho}{\dd\mu_V},\mathcal T_V
 \right)
 \leq C_n{\dLSI(\varrho\mid\mu_V)}^{1/2}.
\end{equation}
The constant is independent of both $V$ and $\kappa$, and the exponent
$1/2$ is optimal.
\end{proposition}

\begin{proof}
This is the argument of Krist\'aly and P{\^\i}rvuceanu~\cite[Section~3.1]{kristaly2026quantitative}, with the sharp Prékopa--Leindler estimate described above substituted for the non-sharp Prékopa--Leindler estimate. The reader is strongly encouraged to read their work, but we sill recall the short linearization computation for completeness.  

Assume first that $g=\log(\dd\varrho/\dd\mu_V)$ is smooth and bounded, and consider the following triplet of functions satisfying the Prékopa condition:
\begin{equation}\label{eq.app_LS_PL_triple}
 u_\lambda=e^{g/(1-\lambda)-V},\qquad
 v=e^{-V},\qquad
 w_\lambda=e^{g_\lambda-V},
\end{equation}
where $w_\lambda$ is their Asplund sum.  Strong convexity of $V$ gives
\[
 g_\lambda
 \leq g+\frac{\lambda}{2\kappa(1-\lambda)}|\nabla g|^2.
\]
Since $\int v\dd x=1$, the relative Prékopa--Leindler deficit of this
triple is
\begin{equation}\label{eq.app_LS_PL_deficit}
 \begin{aligned}
 0 \le 
 \dPL(\lambda)
 &\eqdef
 \frac{
   \displaystyle\int w_\lambda\dd x
   -\left(\displaystyle\int u_\lambda\dd x\right)^{1-\lambda}
 }{
   \left(\displaystyle\int u_\lambda\dd x\right)^{1-\lambda}
 }
 \end{aligned}
\end{equation}
Its first variation at $\lambda=0$ satisfies
\begin{equation}\label{eq.app_LS_endpoint}
 \limsup_{\lambda\downarrow0}\frac{\dPL(\lambda)}{\lambda}
 \leq \frac1{2\kappa}\int |\nabla g|^2e^g\dd\mu_V
       -\int ge^g\dd\mu_V
 =\dLSI(\varrho\mid\mu_V).
\end{equation}
Writing $\tau_\lambda=\min\{\lambda,1-\lambda\}$, the sharp
Prékopa--Leindler estimate controls the input distance by
$C_n\sqrt{\dPL(\lambda)/\tau_\lambda}$.  The normalized inputs converge in
$L^1$ to the densities of $\varrho$ and $\mu_V$, and
$\lambda/\tau_\lambda\to1$.  Passing to the limit therefore
proves~\eqref{eq.app_sharp_LS}. The centered Gaussian perturbations in~\cite[Section~3.3]{kristaly2026quantitative} show that the exponent cannot be
larger than $1/2$.
\end{proof}

For the standard Gaussian $\gamma_n$, the sharp theorem of Dolbeault--Esteban--Figalli--Frank--Loss~\cite{dolbeault2025sharp} gives
\begin{equation}\label{eq.app_gaussian_Hellinger}
 \dist_{L^2(\gamma_n)}\left(f,\sqrt{\mathcal O_{\gamma_n}}\right)
 \leq C\sqrt{\dLSI(f^2\gamma_n\mid\gamma_n)}, 
\end{equation}
where 
$
  \sqrt{\mathcal O_{\gamma_n}}
  \eqdef 
  \left\{
    \sqrt{f} : f \in \mathcal O_{\gamma_n}
  \right\} 
$. For any two probability measures $\varrho_0,\varrho_1\ll\gamma_n$
such that $\sqrt{\dd\varrho_i/\dd\gamma_n}\in W^{1,2}(\gamma_n)$, a
simple application of the Cauchy-Schwarz inequality allows us to
deduce~\eqref{eq.app_sharp_LS} from~\eqref{eq.app_gaussian_Hellinger} for the standard Gaussian measure. Although their arguments are only valid for the standard Gaussian, it is clear that a $L^2$ wave-function estimate implies an $L^1$ distribution estimate. However, it cannot supply a $L^p$ control of the distributions for $p>1$. Indeed, the next Proposition shows that $L^1$ is the largest exponent for which a density estimate
of the form in Proposition~\ref{prop.app_sharp_LS} can hold. Differently from the analogous result for the Brascamp-Lieb inequality~\cite[Theorem~3.5]{machado2025quantitative}, we actually observe that the obstruction is stronger than dependence on the potential, as we show that for a fixed Gaussian reference measure, the LSI deficit cannot control any $L^p$ distance for $p>1$. 

\begin{proposition}\label{prop.app_Lp_LS_failure}
For every $n\geq1$, $1<p<\infty$, and $\epsilon>0$,
\begin{equation}\label{eq.app_Lp_LS_failure}
 \inf_{\varrho\in\mathscr P(\R^n)}
 \left\{
 \left.
 \frac{{\dLSI(\varrho\mid\gamma_n)}^{1/2}}
 {\dist_{L^p(\gamma_n)}\left(
    \dfrac{\dd\varrho}{\dd\gamma_n},\mathcal O_{\gamma_n}
  \right)}
 \ \right|\ 
 \begin{gathered}
  0<\dLSI(\varrho\mid\gamma_n)<\epsilon,\\
  0<\dist_{L^p(\gamma_n)}\left(
       \dfrac{\dd\varrho}{\dd\gamma_n},\mathcal O_{\gamma_n}
     \right)<\infty
 \end{gathered}
 \right\}=0.
\end{equation}
\end{proposition}

\begin{proof}
We work in dimension one since a tensoring argument with
$\gamma_{n-1}$ proves the claim in every dimension. For $0<\varepsilon<1$, set $ \varrho_\varepsilon \eqdef \mathcal N(\varepsilon^{-1},1-\varepsilon)$. Its density with respect to $\gamma_1$ is
\[
 \frac{\dd\varrho_\varepsilon}{\dd\gamma_1}(x)
 =(1-\varepsilon)^{-1/2}
 \exp\left(-\frac{(x-\varepsilon^{-1})^2}{2(1-\varepsilon)}
 +\frac{x^2}{2}\right).
\]
Direct Gaussian integration gives
\begin{equation}\label{eq.app_LS_gaussian_deficit}
 \dLSI(\varrho_\varepsilon\mid\gamma_1)
 =\frac12\left((1-\varepsilon)^{-1}-1+\log(1-\varepsilon)\right),
\end{equation}
which tends to zero as $\varepsilon\downarrow0$.

To estimate the distance to $\mathcal{O}_{\gamma_1}$, first notice that multiplication by $\gamma_1^{1/p}$ identifies $L^p(\gamma_1)$ with
$L^p(\dd x)$, indeed for any $G_z \eqdef e^{zx-z^2/2} \in \mathcal{O}_{\gamma_1}$ we have that 
\[
  \norm{
    \frac{\dd \varrho_\varepsilon}{\dd \gamma_1} - G_z
  }_{L^p(\gamma_1)} 
  = 
  \norm{
    \frac{\dd \varrho_\varepsilon}{\dd \gamma_1} \gamma_1^{1/p} - G_z \gamma_1^{1/p}
  }_{L^p(\mathbb{R}^n)}. 
\]
Completing the square gives
\begin{equation}\label{eq.app_LS_first_completed_square}
 \begin{aligned}
 \frac{\dd\varrho_\varepsilon}{\dd\gamma_1}(x)
 \gamma_1(x)^{1/p}
 &=
 c_{p,\varepsilon}
 \exp\left(-\frac{|x-m_\varepsilon|^2}{2\sigma_\varepsilon^2}\right),
 \end{aligned}
\end{equation}
with the constants 
\[
 m_\varepsilon
 =\frac{p}{\varepsilon[1+(p-1)\varepsilon]},
 \quad
 \sigma_\varepsilon^2
 =\frac{p(1-\varepsilon)}{1+(p-1)\varepsilon},
 \quad 
 c_{p,\varepsilon} 
 = 
  (2\pi{(1-\varepsilon)}^{p})^{-\frac{1}{2}}
 \exp
    \left(
      \frac{p-1}{2\varepsilon^2[1+(p-1)\varepsilon]}
    \right)
 .
\]
On the other hand, every element of $\mathcal O_{\gamma_1}$ satisfies
\begin{equation}\label{eq.app_LS_optimizer_completed_square}
 G_z(x)\gamma_1(x)^{1/p}
 =\frac{1}{(2\pi)^{1/(2p)}}
 \exp\left(\frac{p-1}{2}z^2\right)
 \exp\left(-\frac{|x-pz|^2}{2p}\right).
\end{equation}
Thus the Gaussian in~\eqref{eq.app_LS_first_completed_square} has center
$m_\varepsilon$ and variance $\sigma_\varepsilon^2$, whereas those in~\eqref{eq.app_LS_optimizer_completed_square} have arbitrary center and fixed variance $p$. Noticing that
$p-\sigma_\varepsilon^2 =\frac{p^2\varepsilon}{1+(p-1)\varepsilon}$, there is a constant $k_p>0$ such that, for $0<\varepsilon\leq1/2$, we have
\begin{equation}\label{eq.app_LS_gaussian_shape_gap}
 \inf_{A\geq0,\,y\in\R}
 \left\|
   \exp\left(-\frac{|x|^2}{2\sigma_\varepsilon^2}\right)
   -A\exp\left(-\frac{|x-y|^2}{2p}\right)
 \right\|_{L^p(\dd x)}
 \geq k_p\varepsilon.
\end{equation}

We now translate the first Gaussian to the origin, factor out its
multiplicative constant, and enlarge the family in
\eqref{eq.app_LS_optimizer_completed_square} by allowing its amplitude and
center to vary independently. From~\eqref{eq.app_LS_gaussian_shape_gap},
we obtain
\[
 \dist_{L^p(\gamma_1)}
 \left(\frac{\dd\varrho_\varepsilon}{\dd\gamma_1},
 \mathcal O_{\gamma_1}\right)
 \geq
 c_{p,\varepsilon} k_p \varepsilon
 \geq c_p\varepsilon e^{c_p/\varepsilon^2}.
\]
This tends to infinity as $\varepsilon\downarrow0$.
\end{proof}

\subsection{The Talagrand transport inequality}
\label{subsec.app_talagrand}

Let again $D^2V\geq \kappa I_n$, and let
$\varrho\ll\mu_V$ belong to $\mathscr{P}_2(\R^n)$. Talagrand's inequality states that the relative entropy between $\varrho$ and a reference strongly log-concave measure $\mu_V$ controls the squared Wasserstein distance between them. The inequality and its deficit are written as
\[
 W_2^2(\varrho,\mu_V)\leq\frac2\kappa H(\varrho\mid\mu_V),
 \qquad
 \dTal(\varrho\mid\mu_V)
 =H(\varrho\mid\mu_V)-\frac \kappa2W_2^2(\varrho,\mu_V).
\]
Krist\'aly and P{\^\i}rvuceanu~\cite[Section~4]{kristaly2026quantitative}
showed how to derive stability by linearizing the
Pr\'ekopa--Leindler inequality. 

To make explicit the optimal transport structure behind this linearization, consider the quadratic cost $c(x,y)\eqdef \frac{\kappa}{2}|x-y|^2$, and we recall the $c$-transform defined as
\[
 \psi^c(y)
 \eqdef\inf_x\left\{\frac{\kappa}{2}|x-y|^2-\psi(x)\right\}.
\]
Then, Kantorovich duality states that
\[
  \frac{\kappa}{2}W_2^2(\varrho,\mu_V)
  =\sup_{\psi\in C_b(\R^n)}
  \left\{
    \int\psi^c\dd\varrho+\int\psi\dd\mu_V
  \right\}.
\]
Let $\psi$ be a Kantorovich potential attaining the supremum above. The
optimal transport map $T$ from $\mu_V$ to $\varrho$ is given by
\[
 T=\id-\kappa^{-1}\nabla\psi.
\]
Thus $(\psi,\psi^c)$ is an optimal pair for the Kantorovich dual problem,
and the equality on the support of the
optimal transport plan gives
\begin{equation}\label{eq.app_Tal_duality}
 \psi(x)+\psi^c(T(x))=\frac{\kappa}{2}|x-T(x)|^2
 \quad\text{for $\mu_V$-a.e. }x,
 \qquad
 \int\psi^c\dd\varrho+\int\psi\dd\mu_V
 =\frac{\kappa}{2}W_2^2(\varrho,\mu_V).
\end{equation}
We refer to~\cite[Chap.~1]{santambrogio2015optimal} for the
Kantorovich duality, the $c$-transform, and the characterization of the
optimal map for the quadratic cost.

We shall also use the variational formula for the relative entropy
\begin{equation}\label{eq.app_entropy_variational}
 H(\varrho\mid\mu_V)
 =\sup_{\zeta\in C_b(\R^n)}
 \left\{
   \int\zeta\dd\varrho
   -\log Z_\zeta
\right\},
\text{ where }
Z_\zeta \eqdef\int e^{\zeta}\dd\mu_V,
\end{equation}
with the usual extension to admissible unbounded functions by truncation;
see~\cite{ambrosio2021lectures}. 

\begin{proposition}\label{prop.app_sharp_Tal}
There exists a dimensional constant $C_n>0$ such that every
$\varrho\ll\mu_V$ with finite entropy and second moment satisfying
$\dTal(\varrho\mid\mu_V)<C_n^{-2}$ obeys
\[
 \dist_{L^1(\mu_V)}\left(
   \frac{\dd\varrho}{\dd\mu_V},\mathcal T_V
 \right)
 \leq C_n {\dTal(\varrho\mid\mu_V)}^{1/2}.
\]
The constant is independent of $V$ and $\kappa$, and the exponent $1/2$ is
optimal.
\end{proposition}

\begin{proof}
We follow the endpoint computation in
\cite[Section~4.1]{kristaly2026quantitative}, referring there for the
transport and integrability details. The major idea is to use optimal to perform perturbations of the Prékopa--Leindler inequality with the optimal Kantorovich potentials discussed above and relate them to the relative entropy by means of the variational formula~\eqref{eq.app_entropy_variational}. With weights $1-\lambda,\lambda$, set
\[
 u_\lambda=e^{\lambda\psi-V},\qquad
 v_\lambda=e^{(1-\lambda)\psi^c-V},\qquad
 w=e^{-V}.
\]
The defining inequality for $\psi^c$ and the strong convexity of $V$ yield the
Pr\'ekopa condition. Since $\int w\dd x=1$, the relative Prékopa--Leindler deficit is
\begin{equation}\label{eq.app_Tal_PL_deficit}
 \begin{aligned}
 0 \le \dPL(\lambda)
 &\eqdef
 \frac{
   \displaystyle\int w\dd x
   -\left(\displaystyle\int u_\lambda\dd x\right)^{1-\lambda}
    \left(\displaystyle\int v_\lambda\dd x\right)^\lambda
 }{
   \left(\displaystyle\int u_\lambda\dd x\right)^{1-\lambda}
   \left(\displaystyle\int v_\lambda\dd x\right)^\lambda
 }
 \geq0.
 \end{aligned}
\end{equation}
If $\tau_\lambda=\min\{\lambda,1-\lambda\}$, then the first variation
computed in~\cite[Section~4.1]{kristaly2026quantitative} gives
\begin{equation}\label{eq.app_Tal_endpoint}
 \lim_{\lambda\downarrow0}
 \frac{\dPL(\lambda)}{\tau_\lambda}
 =\mathcal B_V(\psi),
 \qquad
 \mathcal B_V(\psi)
 \eqdef-\int\psi\dd\mu_V-\log Z_{\psi^c}.
\end{equation}

The variational formula~\eqref{eq.app_entropy_variational}, tested with
$\psi^c-\log Z_{\psi^c}$, gives
\[
 \int\psi^c\dd\varrho-\log Z_{\psi^c}\leq H(\varrho\mid\mu_V).
\]
Combining this estimate with the Kantorovich duality
identity~\eqref{eq.app_Tal_duality} yields
\begin{equation}\label{eq.app_Tal_endpoint_bound}
 \mathcal B_V(\psi)\leq\dTal(\varrho\mid\mu_V).
\end{equation}

The normalized inputs in the Prékopa--Leindler estimate converge, as
$\lambda\downarrow0$, to the densities $1$ and $Z_{\psi^c}^{-1}e^{\psi^c}$ with
respect to $\mu_V$.  Equations~\eqref{eq.app_Tal_endpoint} and
\eqref{eq.app_Tal_endpoint_bound}, together with the sharp dependence on
$\tau_\lambda$, therefore give
\[
 \dist_{L^1(\mu_V)}\left(
   Z_{\psi^c}^{-1}e^{\psi^c},\mathcal T_V
 \right)
 \leq C_n\sqrt{\mathcal B_V(\psi)}.
\]
Finally, the same duality identity gives the exact entropy remainder
\[
 \int\log\left(
   \frac{\dd\varrho/\dd\mu_V}{Z_{\psi^c}^{-1}e^{\psi^c}}
 \right)\dd\varrho
 =\dTal(\varrho\mid\mu_V)-\mathcal B_V(\psi).
\]
Pinsker's and the triangle inequalities now imply
\[
 \begin{aligned}
 \dist_{L^1(\mu_V)}\left(
   \frac{\dd\varrho}{\dd\mu_V},\mathcal T_V
 \right)
 &\leq
 \left\|\frac{\dd\varrho}{\dd\mu_V}-Z_{\psi^c}^{-1}e^{\psi^c}\right\|_{L^1(\mu_V)}
 +C_n\sqrt{\mathcal B_V(\psi)}\\
 &\leq
 \sqrt{2\left(\dTal(\varrho\mid\mu_V)-\mathcal B_V(\psi)\right)}
 +C_n\sqrt{\mathcal B_V(\psi)}
 \leq C_n\sqrt{\dTal(\varrho\mid\mu_V)}.
 \end{aligned}
\]
Centered Gaussian dilations, as in~\cite[Section~4.4]{kristaly2026quantitative}, show that the exponent is optimal.
\end{proof}

We finally test whether the potential-uniform estimate can be strengthened
to an $L^p$ density distance.  As for log-Sobolev, every power above one
fails already for the Gaussian potential.

\begin{proposition}\label{prop.app_Lp_Tal_failure}
For every $n\geq1$, $1<p<\infty$, and $\epsilon>0$,
\begin{equation}\label{eq.app_Lp_Tal_failure}
 \inf_{
    \varrho\in\mathscr{P}_2(\R^n)
 }
 \left\{
 \left.
 \frac{{\dTal(\varrho\mid\gamma_n)}^{1/2}}
 {\dist_{L^p(\gamma_n)}\left(
    \frac{\dd\varrho}{\dd\gamma_n},\mathcal O_{\gamma_n}
  \right)}
 \ \right|\ 
 \begin{gathered}
  0<\dTal(\varrho\mid\gamma_n)<\epsilon,\\
  0<\dist_{L^p(\gamma_n)}\left(
       \frac{\dd\varrho}{\dd\gamma_n},\mathcal O_{\gamma_n}
     \right)<\infty
 \end{gathered}
 \right\}=0.
\end{equation}
\end{proposition}

\begin{proof} 
Once more we work in dimension one since a tensoring argument with $\gamma_{n-1}$ proves the claim in every dimension. Fix $s>p/(p-1)$. For $0<\varepsilon<1$, consider the increasing odd map
\[
 T_\varepsilon(x)=
 \begin{cases}
 x,&|x|\leq \varepsilon^{-1},\\
 \sign(x)\bigl(s|x|-(s-1)\varepsilon^{-1}\bigr),
 &\varepsilon^{-1}<|x|\leq2\varepsilon^{-1},\\
 x+(s-1)\varepsilon^{-1}\sign(x),&|x|>2\varepsilon^{-1},
 \end{cases}
\]
and define $\varrho_\varepsilon\eqdef(T_\varepsilon)_\sharp\gamma_1$. The map $T_\varepsilon$ is a bi-Lipschitz bijection, and therefore $\varrho_\varepsilon$ is an absolutely continuous centered probability measure. Moreover, $T_\varepsilon'=s$ on $\{\varepsilon^{-1}<|x|<2\varepsilon^{-1}\}$ and $T_\varepsilon'=1$ elsewhere. Gaussian change of variables followed by integration by parts gives
\begin{equation}\label{eq.app_Tal_1d_deficit}
 \dTal(\varrho_\varepsilon\mid\gamma_1)
 =\int(T_\varepsilon'-1-\log T_\varepsilon')\dd\gamma_1.
\end{equation}
Consequently,
\begin{equation}\label{eq.app_Tal_tail_deficit}
 0<\dTal(\varrho_\varepsilon\mid\gamma_1)
 =(s-1-\log s)
 \gamma_1(\varepsilon^{-1}<|x|<2\varepsilon^{-1}),
\end{equation}
which tends to zero as $\varepsilon\downarrow0$.

Using $\varrho_\varepsilon$ also for its Lebesgue density, a Gaussian change of variables gives
\[
 \frac{\varrho_\varepsilon}{\gamma_1}(T_\varepsilon(x))
 =\frac{1}{T'_\varepsilon(x)}
   \exp\left(\frac{T_\varepsilon(x)^2-x^2}{2}\right),
\]
and one can check that 
\[
 \begin{aligned}
 \left\|\frac{\varrho_\varepsilon}{\gamma_1}\right\|_{L^p(\gamma_1)}^p
 &=\int
 \left(\frac{\varrho_\varepsilon}{\gamma_1}(T_\varepsilon(x))\right)^{p-1}
 \dd\gamma_1(x)\\
 &\geq
 \frac{s^{1-p}}{\sqrt{2\pi}}
 \int_{\varepsilon^{-1}}^{2\varepsilon^{-1}}
 \exp\left(
   \frac{(p-1)(sx-(s-1)\varepsilon^{-1})^2-px^2}{2}
 \right)\dd x.
 \end{aligned}
\]
By continuity, there are constants $c_p>0$ and $\eta_p\in(0,1)$ such that
\[
 \left\|\frac{\varrho_\varepsilon}{\gamma_1}\right\|_{L^p(\gamma_1)}^p
 \geq\frac{c_p}{\varepsilon}
       \exp\left(\frac{c_p}{\varepsilon^2}\right),
\]
and in particular this norm tends to infinity as
$\varepsilon\downarrow0$.

It remains to compute the distance to the set of Gaussian optimizers. If $z\geq0$, then $e^{zx-z^2/2}\leq1$ on the negative half-line; if $z\leq0$, the same
holds on the positive half-line. Since $\varrho_\varepsilon/\gamma_1$ is even,
\[
 \left\|\frac{\varrho_\varepsilon}{\gamma_1}-e^{zx-z^2/2}\right\|_{L^p(\gamma_1)}^p
 \geq\frac12\int
 \left(\frac{\varrho_\varepsilon}{\gamma_1}-1\right)_+^p\dd\gamma_1
 \geq 2^{-p}
 \left\|\frac{\varrho_\varepsilon}{\gamma_1}\right\|_{L^p(\gamma_1)}^p
 -\frac12.
\]
Taking the infimum over $z$ shows that
\[
 \dist_{L^p(\gamma_1)}
 \left(\frac{\varrho_\varepsilon}{\gamma_1},\mathcal O_{\gamma_1}\right)
 \longrightarrow\infty
 \qquad\text{as }\varepsilon\downarrow0.
\]
Together with~\eqref{eq.app_Tal_tail_deficit}, this gives
\[
 \frac{\dTal(\varrho_\varepsilon\mid\gamma_1)^{1/2}}
 {\dist_{L^p(\gamma_1)}
  \left(\varrho_\varepsilon/\gamma_1,\mathcal O_{\gamma_1}\right)}
 \longrightarrow0.
\]
For every prescribed $\epsilon>0$, the deficit belongs to $(0,\epsilon)$ for all sufficiently small $\varepsilon$, which proves the claim in dimension one.

\end{proof}

Propositions~\ref{prop.app_Lp_LS_failure} and
\ref{prop.app_Lp_Tal_failure} are the counterparts, for log-Sobolev and Talagrand, of the zero-infimum criterion obtained for Brascamp--Lieb in \cite[Theorem~3.5]{machado2025quantitative}. They show more than a loss of
uniformity as the potential varies: for these density distances, no square-root stability estimate in $L^p$, $p>1$, can hold even after fixing the Gaussian potential and restricting to arbitrarily small deficits. Together with Propositions~\ref{prop.app_sharp_LS} and
\ref{prop.app_sharp_Tal}, this isolates a potential-independent $L^1$ stability theory generated by Pr\'ekopa--Leindler. The wave-function estimate \eqref{eq.app_gaussian_Hellinger} belongs to a parallel $L^2$ theory and is fully compatible with this conclusion.

\bibliographystyle{alpha}
\bibliography{refs.bib}

\end{document}